\documentclass[11pt]{iopart}

\usepackage{xcolor}
\usepackage{hyperref}
\usepackage{csquotes}
\usepackage{booktabs}
\expandafter\let\csname equation*\endcsname\relax
\expandafter\let\csname endequation*\endcsname\relax
\usepackage{amsmath}
\usepackage{marginnote}
\usepackage{iopams}  
\usepackage[lined,linesnumbered]{algorithm2e}
	\RestyleAlgo{ruled}
\usepackage[pdftex]{graphicx}
\usepackage{titlesec}
\usepackage[many]{tcolorbox}
\usepackage{enumitem}

\newcommand{\new}[1]{#1}
\newcommand{\removed}[1]{}
\newcommand{\replaced}[2]{\new{#2}}{\removed{#1}}

\newtcbox{\mybox}[1][red]{on line,
arc=0pt,outer arc=0pt,colback=#1!10!white,colframe=#1!50!black,
boxsep=0pt,left=1pt,right=1pt,top=2pt,bottom=2pt,
boxrule=0pt,bottomrule=1pt,toprule=1pt}

\newtcbox{\myinbox}[1][red]{on line,
arc=0pt,outer arc=0pt,colback=#1!10!white,colframe=#1!10!white,
boxsep=0pt,left=1pt,right=1pt,top=2pt,bottom=2pt,
boxrule=0pt,bottomrule=1pt,toprule=1pt}

\newtcbox{\proxbox}{enhanced,nobeforeafter,tcbox raise base,boxrule=0pt,top=2mm,bottom=0mm,
  right=0mm,left=0mm,arc=1pt,boxsep=2pt,before upper={\vphantom{dlg}},
  colframe=green!50!black,coltext=green!25!black,colback=green!10!white,
  overlay={\begin{tcbclipinterior}\fill[green!75!blue!50!white] ([yshift=-2mm]frame.north west)
    rectangle node[text=black,font=\sffamily\bfseries\tiny] {nonsmooth} ([yshift=0mm]frame.north east);\end{tcbclipinterior}}}

\newtcbox{\gradbox}{enhanced,nobeforeafter,tcbox raise base,boxrule=0pt,top=2mm,bottom=0mm,
  right=0mm,left=0mm,arc=1pt,boxsep=1pt,before upper={\vphantom{dlg}},
  colframe=red!50!black,coltext=red!25!black,colback=red!10!white,
  overlay={\begin{tcbclipinterior}\fill[red!75!blue!50!white] ([yshift=-2mm]frame.north west)
    rectangle node[text=black,font=\sffamily\bfseries\tiny] {smooth} ([yshift=0mm]frame.north east);\end{tcbclipinterior}}}

\newtcbox{\mombox}{enhanced,nobeforeafter,tcbox raise base,boxrule=0pt,top=2mm,bottom=0mm,
  right=0mm,left=0mm,arc=1pt,boxsep=1pt,before upper={\vphantom{dlg}},
  colframe=orange!50!black,coltext=orange!25!black,colback=orange!10!white,
  overlay={\begin{tcbclipinterior}\fill[orange!75!blue!50!white] ([yshift=-2mm]frame.north west)
    rectangle node[text=black,font=\sffamily\bfseries\tiny] {momentum} ([yshift=0mm]frame.north east);\end{tcbclipinterior}}}

\newtcbox{\dualbox}{enhanced,nobeforeafter,tcbox raise base,boxrule=0pt,top=2mm,bottom=0mm,
  right=0mm,left=0mm,arc=1pt,boxsep=1pt,before upper={\vphantom{dlg}},
  colframe=blue!50!black,coltext=blue!25!black,colback=blue!10!white,
  overlay={\begin{tcbclipinterior}\fill[blue!75!green!50!white] ([yshift=-2mm]frame.north west)
    rectangle node[text=black,font=\sffamily\bfseries\tiny] {dualisation} ([yshift=0mm]frame.north east);\end{tcbclipinterior}}}

\usepackage{diagbox}
\usepackage{multirow}

\usepackage{makecell}
\usepackage{arydshln}

\titleformat{\section}{\normalfont \bfseries}{\thesection}{1em}{}
\titleformat{\subsection}{\normalfont \bfseries}
  {\thesubsection}{1em}{}
\titleformat{\subsubsection}{\normalfont \bfseries\slshape}{\thesubsubsection}{1em}{}  
\titleformat{\paragraph}[runin]
{\bfseries\slshape}{\theparagraph}{1em}{}  

\usepackage{subfig}

\usepackage{tikz}
\usetikzlibrary{arrows.meta}

\newcommand{\X}{\mathcal{X}}

\newcommand{\E}{\mathbb{E}}

\newcommand{\calL}{\mathcal{L}} 

\newcommand{\calO}{\mathcal{O}}

\newcommand{\ek}{e^{(k)}}
\newcommand{\xk}{x^{(k)}}
\newcommand{\xkm}{x^{(k-1)}}
\newcommand{\xkp}{x^{(k+1)}}
\newcommand{\yk}{y^{(k)}}

\newcommand{\ykp}{y^{(k+1)}}
\newcommand{\zk}{z^{(k)}}

\newcommand{\zkp}{z^{(k+1)}}

\newcommand{\wk}{w^{(k)}}

\newcommand{\prox}{\operatorname{prox}}

\newtheorem{remark}{Remark}
\newcommand{\argmin}{\operatornamewithlimits{argmin}}

\newcommand{\iprod}[2]{\langle #1, #2 \rangle}

\newcommand{\norm}[1]{{|\kern-1.125pt|} #1 {|\kern-1.125pt|}}

\def\Lmax{L_{\max}}

\def\reconspace{\mathcal U}
\def\reconvar{u}
\def\dataspace{\mathcal V}
\def\datavar{v}
\def\ipop{K}
\def\R{\mathbb R}
\def\RI{\R_\infty}

\def\datanum{s}

\def\recondim{d}
\def\numf{\ell}
\def\numg{m}
\def\numh{n}

\def\optspace{\mathcal X}
\def\optvar{x}
\def\dualspace{\mathcal Y}
\def\dualvar{y}
\def\auxspace{\mathcal Z}
\def\auxvar{z}
\def\optop{A}

\def\optfuna{f}
\def\optfunb{g}
\def\optfunc{h}

\def\optfun{\Phi}

\newcounter{dummy} \numberwithin{dummy}{section}

\newtheorem{myex}[dummy]{Example}

\begin{document}

\title[A Guide to Stochastic Optimisation for Large-Scale Inverse Problems]{A Guide to Stochastic Optimisation for \\ Large-Scale Inverse Problems}

\author{Matthias J. Ehrhardt$^1$, \v Zeljko Kereta$^2$, Jingwei Liang$^3$, \\ Junqi Tang$^4$}

\address{$^1$University of Bath, Bath BA2 7JU, UK\\
$^2$University College London, London, UK\\
$^3$Shanghai Jiao Tong University, Shanghai, China\\
$^4$University of Birmingham, Birmingham B15 2TT, UK}
\ead{m.ehrhardt@bath.ac.uk}
\vspace{10pt}
\begin{indented}
\item[]\today 
\end{indented}

\begin{abstract}
Stochastic optimisation algorithms are the de facto standard for machine learning with large amounts of data. 
Handling only a subset of available data in each optimisation step dramatically reduces the per-iteration computational costs, while still ensuring significant progress towards the solution.
Driven by the need to solve large-scale optimisation problems as efficiently as possible, the last decade has witnessed an explosion of research in this area. 
Leveraging the parallels between machine learning and inverse problems has allowed harnessing the power of this research wave for solving inverse problems. 
In this survey, we provide a comprehensive account of the state-of-the-art in stochastic optimisation from the viewpoint of \new{variational regularisation for} inverse problems \new{where the solution is modelled as minimising an objective function}. 
We present algorithms with diverse modalities of problem randomisation and discuss the roles of variance reduction, acceleration, higher-order methods, and other algorithmic modifications, and compare theoretical results with practical behaviour.
We focus on the potential and the challenges for stochastic optimisation that are unique to \new{variational regularisation for} inverse imaging problems and are not commonly encountered in machine learning.  
We conclude the survey with illustrative examples from imaging \removed{problems} \new{on linear inverse problems} to examine the advantages and disadvantages that this new generation of algorithms bring to the field of inverse problems.  
\end{abstract}

\vspace{2pc}
\noindent{\it Keywords}: Inverse Problems, Imaging, Variational Regularisation, First-Order Algorithms, Large-Scale Optimisation, Stochastic Optimisation

\submitto{\IP}

\newpage

\begin{minipage}{0.85\textwidth}
\begin{small}
    \setcounter{tocdepth}{2}
    \tableofcontents
\end{small}    
\end{minipage}

\section{Inverse Problems Meets Optimisation} \label{sec:ip}

This review is concerned with the efficient solution of optimisation problems in the context of inverse problems. This is particularly relevant when data sets are large or forward operators are costly to evaluate.
In order to set the scene, we are interested in recovering an object $\reconvar^\dagger \in \reconspace$ from indirect noisy measurements $\datavar \in \dataspace$.
We assume that the relationship between $\reconvar$ and $\datavar$ can be approximated by the forward model
\begin{align}
    \ipop (\reconvar^\dagger) \approx \datavar, \label{eq:ip}
\end{align}
where $\ipop : \reconspace \to \dataspace$ is a mapping between  Hilbert spaces\footnote{While some mathematical modelling for inverse problems is formulated on Banach spaces, see e.g.~\cite{Benning2018actanumerica}, for optimisation, it is vastly more convenient to consider Hilbert spaces \cite{bauschke2011convex, Hohage2014}. The dimension of such spaces is not as important and many results are dimension independent.} $\reconspace$ and $\dataspace$. 
Many tasks can be described through this expression, such as deblurring in image processing, tomographic reconstruction in imaging or regression in statistics, to name a few. Concrete applications in tomography include computed tomography (CT) \cite{Hansen2021}, positron emission tomography (PET) \cite{Ollinger1997a} and magnetic resonance imaging (MRI) \cite{Fessler2010a}.

In this review we primarily consider linear inverse problems, where the mapping $\ipop$ is a linear operator and we denoting its action on elements of $\reconspace$ by $\ipop\reconvar$. 
Even though nonlinearity introduces additional mathematical and computational challenges \new{(like nonconvexity)}, many techniques presented in this paper can be readily extended to nonlinear problems.

In applications both the solution $\reconvar$ and data $\datavar$ are typically discretised, and we seek to recover an element in $\R^\recondim$ from measurements in $\R^\datanum$. Complex numbers can be readily incorporated by identifying $\mathbb C$ with $\R^2$. An entry $\datavar_i$ represents the $i$-th measurement, e.g.\ of a digital detector, and $\reconvar_i$ corresponds to the $i$-th parameter of the solution, e.g.\ the value of a particular image voxel. 

Many modern applications have in common that either $\recondim$ and $\datanum$ are both large, or that the forward operator $\ipop$ is expensive to evaluate.
Consequently, computing an approximate solution to the inverse problem \eqref{eq:ip} is computationally challenging. 
These two issues are generally highly related as e.g.\ the cost of matrix-vector products scales linearly with both the number of input and output dimensions. To illustrate the point of high dimensions, consider the PET-MR scanner Siemens Biograph mMR where the uncompressed data (so called \enquote{span 1}) has more than 350 million elements and images are commonly modelled with around 20 million voxels. It is important to highlight, that the words \enquote{large} and \enquote{expensive} are relative, and have to be compared to the available resources in a specific application.
 
Instead of directly solving \eqref{eq:ip}, it is common to approximate the solution via an optimisation problem. The typical choice is
\begin{align}
    \min_{\reconvar \in \reconspace} D(\ipop \reconvar, \datavar), \label{eq:opt}
\end{align}
where $D : \dataspace \times \dataspace \to \RI := \R \cup \{\infty\}$ measures how closely $\ipop \reconvar$ approximates $\datavar$, and is commonly called the \emph{data-fidelity} or \emph{data\new{-}fit term}. The data-fidelity term $D$ is typically a \replaced{non-negative}{nonnegative} function minimised by identical arguments, $\inf_{\datavar'} D(\datavar', \datavar) = D(\datavar,\datavar) = 0$. 
Allowing $D$ to take on the value $\infty$ is sometimes needed to extend the definition of $D$ to a vector space (e.g.\ see the Kullback--Leibler divergence below) or to model constraints. See \cite{Bredies2018book} for a discussion on extended real numbers~$\RI$. 

An alternative viewpoint is to model, and then maximise, the likelihood of observing $\datavar$ given $\ipop \reconvar$. The data-fidelity can then be defined as the negative log-likelihood and the sought estimator is the one maximising the likelihood, see e.g.~\cite{Bertero1998book, Benning2018actanumerica} for more dedicated discussions.

Most commonly used data fidelities are \emph{separable}.
In the finite-dimensional data setting, $\dataspace = \R^\datanum$, this means that they can be written as
\begin{align}
    D(\datavar', \datavar) = \sum_{i=1}^\datanum d(\datavar_i', \datavar_i). \label{eq:separable}
\end{align}
The function $d:\R\times\R\rightarrow\R_\infty$ can be determined using statistical properties of the data or the measurement noise. 
For example, using the Gaussian, Laplace and Poisson distributions leads to least squares, absolute differences and Kullback--Leibler divergence, respectively, defined as follows
\begin{itemize}
\item \textbf{Least squares}
\begin{align*}
    d(\datavar_i', \datavar_i) = \frac1{2\sigma_i^2} (\datavar_i' - \datavar_i)^2.
\end{align*}
\item \textbf{Absolute differences}
\begin{align*}
    d(\datavar_i', \datavar_i) = |\datavar_i' - \datavar_i|.
\end{align*}
\item \textbf{Kullback--Leibler divergence}
\begin{align*}
    d(\datavar_i', \datavar_i) = \operatorname{KL}(\datavar_i | \datavar_i') = 
    \begin{cases}
        \datavar_i' - \datavar_i + \datavar_i \log(\datavar_i/\datavar_i') & \text{if $\datavar_i' > 0$ and $\datavar_i > 0$}, \\
        \datavar_i' & \text{if $\datavar_i' > 0$ and $\datavar_i = 0$}, \\
        \infty & \text{else}.
    \end{cases}
\end{align*}
\end{itemize}
While the Gaussian assumption is most frequently applied, Poisson distributed data has also gained a lot of attention, see e.g.\ \cite{Setzer2010, Bardsley2000, Hohage2016}.

\subsection{Variational Regularisation}
The original inverse problem \eqref{eq:ip} is ill-posed, in the sense that it can have no solutions, infinitely many solutions or even if the solution exists and is unique, it is highly sensitive to the noise in the measurements.
Solving the optimisation problem \eqref{eq:opt}, which minimises the fidelity of the reconstruction to the measurements, instead of \eqref{eq:ip}, can help overcome the nonexistence of the solution.
However, it inherits other issues, such as \replaced{non-uniqueness}{nonuniqueness} of solutions and noise amplification. For example, if $K$ is ill-conditioned, small changes in the observations (for example due to noise) can lead to large changes in the reconstructed solution.
Moreover, only minimising the data fidelity will result in reconstructions that fit to noise as well.
A common approach to address this is to add a \emph{regulariser} $R : \reconspace \to \RI$ to the optimisation objective and instead solve
\begin{equation}\label{eq:varreg}
    \min_{\reconvar \in \reconspace} \left\{D(\ipop \reconvar, \datavar) + \lambda R(\reconvar)\right\},
\end{equation}
where $\lambda > 0$ is called the \emph{regularisation parameter}. This concept is known as \emph{variational regularisation} and has been extensively studied extensively over the past decades. Many regularisers $R$ have been used over the years, see \cite{Scherzer2008book, Benning2018actanumerica} and the references therein. Instead of providing a comprehensive overview, we will discuss representative regularisers with properties relevant for inverse problems.

If $\mathcal{U} = \R^\recondim$, the simplest regulariser is the squared 2-norm: $R(\reconvar) = \frac12 \|\reconvar\|^2_2 = \frac12 \sum_i \reconvar_i^2$. This is often referred to as Tikhonov or Miller regularisation \cite{Ito2014book, miller1970regularisation} and can be used to stabilise the reconstruction and make the modelled solution unique. Alternatively, using the 1-norm, $R(\reconvar) = \|\reconvar\|_1 = \sum_i |\reconvar_i|$, promotes a specific property in the solution: sparsity, i.e.\ only few coefficients of the solution are nonzero. 

Regularising with 1- and squared 2-norm does not capture well the structure of desired solutions generally encountered in imaging applications. They promote reconstructions with small entries which is not desirable for images.
Moreover, it does not take into account local neighbourhood or smoothness features within the image. 
Therefore, regularisers that penalise differences between neighbouring pixels are better suited for imaging tasks.
This can be achieved by not penalising $\reconvar$ directly but rather its gradient: $R(\reconvar) = \frac12 \|\nabla\reconvar\|^2_2$. 
However, this is in strict terms well-defined only for smooth functions $\reconvar$. An alternative perspective is to interpret $\nabla$ as a linear operator that computes (spatial) finite differences of a discretised vector $\reconvar$. We do not to want delve deeper into the details, but instead refer to e.g.~\cite{Burger2013} or \cite{chambolle2016introduction} for more details in the continuous or discrete setting, respectively. 

Using the squared $2$-norm promotes similarity between neighbouring pixels and can cause blurring on the edges of an object. 
In order to preserve the edges it is common to replace the squared $2$-norm of the gradient with the $1$-norm, as
\begin{equation}\label{eq:tv}
    R(\reconvar) = \norm{\nabla \reconvar}_1.
\end{equation}
This defines the total variation~(TV), which is a particularly good model if we believe that the desired solution can be well described as a piece-wise constant function, a common model for natural images. \new{This formula is well-defined for smooth functions $u$ but can be generalised to other functions with less regularity, see e.g.~\cite{Burger2013}.}
 
It is known that TV introduces staircasing artefacts when an image has piece-wise linear structures, see e.g.\ \cite{Caselles2007}. To circumvent this problem, several high-order variants of TV have been proposed, see e.g.~\cite{CL97, bredies2010total}. 
Perhaps the most popular choice is the total generalised variation (TGV) \cite{bredies2010total}, defined as 
\begin{equation}\label{eq:tgv}
    R(\reconvar) = \inf_w \left\{\norm{\nabla \reconvar - w}_1 + \beta \norm{\nabla w}_1\right\},
\end{equation} 
for some $\beta > 0$. Here $w$ is the vector field of the piece-wise linear component of $u$, while $\nabla u - w$ is the vector field of the piece-wise constant component of $u$. \new{A comparison between TV and TGV is provided in Figure \ref{fig:cmp_tv_tgv}.}

An alternative to the \replaced{non-smooth}{nonsmooth} TV is the so called Huber TV, see e.g.\ \cite{chambolle2016introduction}, which replaces the $1$-norm with a smooth approximation. In the discretised setting it can be written as
\begin{align*}
    R(\reconvar) = \sum_i h_\gamma(\|(\nabla \reconvar)_i\|_2),
\end{align*}
where the sum is over all components of the gradient field $\nabla \reconvar$ and 
\begin{align*}
    h_\gamma(t) = \begin{cases} t & \text{if } |t| > \gamma , \\ 
    \frac1{2\gamma} t^2 + \frac{\gamma}{2} & \text{if }|t| \leq \gamma .
    \end{cases}
\end{align*}
Huber TV is differentiable, with a $\|\nabla\|^2/\gamma$-Lipschitz continuous gradient, and in several applications it overcomes the staircasing artefact that TV introduces.

Another source of regularisation can be formulated as constraints. In contrast to the previously discussed models, which encourage the solution to have certain properties, constraints can be used to instead enforce the desired properties. 
The most popular constraints in imaging are box constraints, where we only seek solutions within a $\recondim$-dimensional box 
\[C = \{x \in \R^\recondim: \ell_i \leq x_i \leq r_i, i=1,2, \ldots, \recondim\}.\]
The most relevant constraints of this type are the \replaced{non-negativity}{nonnegativity} constraints, where one takes $\ell_i=0$ and $r_i = \infty$. 

\begin{remark}[Statistical Perspective]\label{rem:bayesian_variational}
Variational regularisation \eqref{eq:varreg} can be articulated from a Bayesian perspective in terms of maximum a posteriori (MAP) estimation~\new{\cite{stuart2010bayesian, dashti2013map}}.
Given a density $p(\reconvar)$ for the probabilistic model of $u$, MAP considers maximising the posterior probability $p(\reconvar\vert\datavar)\propto p(\datavar\vert\reconvar) p(\reconvar)$, which is equivalent to minimising its negative logarithm, i.e.\ $$\min_{\reconvar} \left\{ -\log p(\datavar\vert\reconvar) - \log p(\reconvar)\right\}.$$

For example, if we assume that the data is Gaussian distributed with mean $\ipop\reconvar$ and covariance $\sigma^2 I$ and that the solution is Gaussian distributed with zero mean and covariance $\lambda^{-1}I$, the MAP can be written as $\argmin_\reconvar \left\{ \sigma^{-2} \|\ipop\reconvar-\datavar\|^2+\lambda\|\reconvar\|^2\right\}$.

The separability of the data fit \eqref{eq:separable} can be motivated from a statistical perspective, and follows for instance from assuming that the data points are conditionally independent given the estimated solution. 

\new{In this review, we focus on the optimization perspective which concerns efficiently computing a MAP estimator. While other statistical aspects of stochastic approximation might also be interesting for some readers, such as MMSE estimation, it is out-of-scope for this optimization-focused review.}
\end{remark}

\subsection{Optimisation Template}\label{subsec:optim_template}

As we just discussed, inverse problems can be solved via optimisation in the context of variational regularisation \eqref{eq:varreg}. All aforementioned optimisation problems are special cases of the \emph{optimisation template}
\begin{align} \label{eq:opttemplate}
    \min_{\optvar \in \optspace} \Bigl\{\optfun(\optvar) = \optfuna(\optop\optvar) + \optfunb(\optvar) + \optfunc(\optvar)\Bigr\}.
\end{align}
\new{Throughout the paper, we consider the following assumptions to the objective
\begin{enumerate}[label={(A.\arabic{*})}, ref={(A.\arabic{*})}, leftmargin=1.5cm]
    \item\label{assump_1} Functions $\optfuna : \dualspace \to \RI$ and $\optfunb : \optspace \to \RI$ are proper closed and convex. Moreover, both are \emph{prox-friendly}, i.e.\ that their proximal operators, see \eqref{eqn:prox_operator} for a definition, can be either evaluated in closed-form or efficiently approximated. 

    \item\label{assump_2} The function $\optfunc : \optspace \to \R$ is convex and $L$-smooth, i.e.~it is continuously differentiable and its gradient has Lipschitz constant $L>0$, see \eqref{eqn:L_smooth} for a definition. 

    \item\label{assump_3} $\optop : \optspace \to \dualspace$ is a bounded linear operator. 
\end{enumerate}
Note that we do not want to (or can) compute the proximal operator of $f\circ \optop$, but instead want to decouple the function $f$ from the operator $\optop$, and treat them separately. Techniques to do so are discussed in Section~\ref{sec:optalgs}.
}
\removed{
Here the functions $\optfuna : \dualspace \to \RI$ and $\optfunb : \optspace \to \RI$ are convex and \emph{prox-friendly}, i.e.\ they have the property that their proximal operators, see \eqref{eqn:prox_operator} for a definition, can be either evaluated in closed-form or efficiently approximated. The function $\optfunc : \optspace \to \R$ is convex and $L$-smooth, i.e.~it is continuously differentiable and its gradient has Lipschitz constant $L>0$, see \eqref{eqn:lipschitz_cont}. 
Lastly, $\optop : \optspace \to \dualspace$ is a bounded linear operator. 
Note that we do not want to (or can) compute the proximal operator of $f\circ \optop$, but instead want to decouple the function $f$ from the operator $\optop$, and treat them separately.
Techniques to do so are discussed in Section~\ref{sec:optalgs}.
}

The conversion from a given minimisation problem for variational regularisation to the optimisation template is not unique as we show in the following example.
\begin{myex} \label{ex:TV}
We consider computing a TV-regularised solution as
\begin{align*}
    \min_{\reconvar \in \reconspace} \left\{\frac 12 \|\ipop\reconvar-\datavar\|_2^2 + \lambda \|\nabla \reconvar\|_1\right\}.
\end{align*}
This can be written via our template \eqref{eq:opttemplate} with $\optvar =\reconvar$, $\optspace=\reconspace$, and 
\begin{itemize}
    \item[{\rm (i)}] $f(y) = \lambda\|y\|_1$, $A = \nabla$, $g(x) = 0$, $h(x) = \frac 12 \|Kx-v\|_2^2$;
    \item[{\rm (ii)}] $f(y) = 0$, $A = 0$, $g(x) = \lambda \|\nabla x\|_1$, $h(x) = \frac 12 \|Kx-v\|_2^2$;
    \item[{\rm (iii)}] $f(y) = 0$, $A = 0$, $g(x) = \|\nabla x\|_1$, $h(x) = \frac 1{2\lambda } \|Kx-v\|_2^2$.
\end{itemize}
The three formulations above (and many more) are all mathematically valid, and used in many applications but lead to algorithms with very different properties. For most of these choices the proximal operators of $f$ and $g$ can be computed in closed-form. For (ii) and (iii), it is common to compute the proximal operator of the total variation $\|\nabla x\|_1$ via warm-started fast gradient projection (FGP) iterations \cite{beck2009fast}, see e.g.\ \cite{Ehrhardt2016mri} for an MRI application.
\end{myex}

The next example illustrates that the optimisation variable $x$ does not need to coincide with the inverse problems variable $u$.
\begin{myex}
    We consider the same inverse problem as in Example~\ref{ex:TV}, but replace the regulariser with TGV \eqref{eq:tgv}, i.e.\ we aim to solve
    \begin{align*}
        \min_{u, w} \left\{\frac 12 \|\ipop\reconvar-v\|_2^2 + \lambda \norm{\nabla \reconvar - w}_1 + \lambda \beta \norm{\nabla w}_1\right\}.
    \end{align*}
    The minimiser can then be computed via template \eqref{eq:opttemplate} with $x = (u,w)$ and
    \begin{align*}
    f(y_1, y_2) = \lambda\|y_1\|_1 + \lambda \beta \|y_2\|_1, \quad \optop = \begin{pmatrix} \nabla & -I \\ 0 & \nabla \end{pmatrix}, \quad g(x) = 0, \quad  h(u, w) = \frac 12 \|\ipop u-v\|_2^2.
    \end{align*}
    Most notably, the optimisation variable $x$ consists of the image to be reconstructed $u$ and the vector-field $w$ in the definition of TGV. 
\end{myex}

\new{To illustrate the difference between TV and TGV regularisations, we consider an image denoising problem, which corresponds to $\ipop$ being identity operator in the above two examples. The outcome is provided in Figure \ref{fig:cmp_tv_tgv}. As shown in Figure \ref{fig:cmp_tv_tgv}(a), the clean image contains piecewise linear structure, the noisy observation is obtained by adding additive white Gaussian noise. The second row shows the denoising outcome using TV and TGV. It can be observed that TGV manages to preserve the linear regions while TV introduces staircase artefacts. 

\begin{figure}[htbp]
    \centering
    \subfloat[Original clean image]{\includegraphics[width=0.475\textwidth]{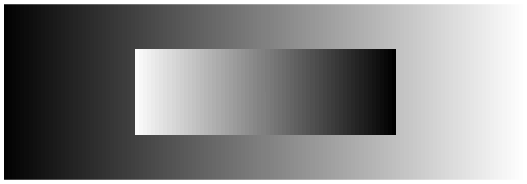}\label{fig:img_original}} 
    \hfill
    \subfloat[Noisy image]{\includegraphics[width=0.475\textwidth]{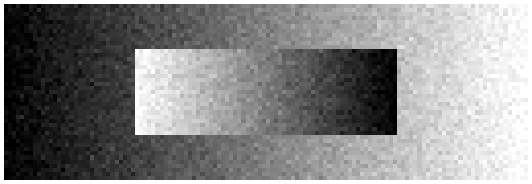}\label{fig:img_denoise}} \\
    \subfloat[Denoising via TV]{\includegraphics[width=0.475\textwidth]{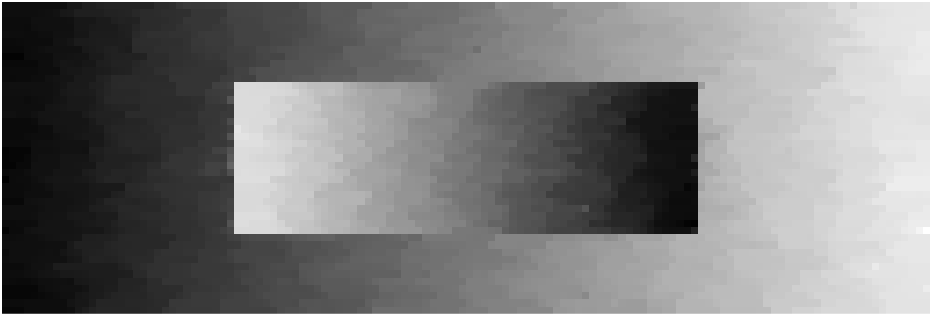}}\label{fig:img_tv} 
    \hfill
    \subfloat[Denoising via TGV]{\includegraphics[width=0.475\textwidth]{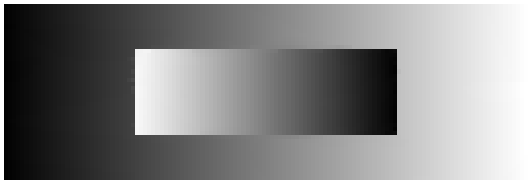}\label{fig:img_tgv}} \\
    
    \caption{Comparison of TV and TGV on image denoising problem. }    
    \label{fig:cmp_tv_tgv}
\end{figure}
}

For the third example, we consider how can constraints be incorporated into the optimisation template.
\begin{myex}
    We consider the same problem as in Example~\ref{ex:TV}, but this time we add the \replaced{non-negativity}{nonnegativity} constraint
    \begin{align*}
        \min_{\reconvar \geq 0} \left\{\frac 12 \|\ipop\reconvar-\datavar\|_2^2 + \lambda \|\nabla \reconvar\|_1\right\}.
    \end{align*}
This can be set into our template \eqref{eq:opttemplate} with $\reconvar = \optvar$, $f(\dualvar) = \|\dualvar\|_1, \optop = \nabla, h(\optvar) = \frac{1}{2} \|\ipop \optvar-\datavar\|_2^2$ and $g(\optvar) = \iota_{[0,\infty)^d}(\optvar)$ is the indicator function of the non-negative orthant (see \eqref{eqn:indicator_function} for the definition). 
Since all $x$ that violate the constraint take the value $\infty$, this choice for $g$ guarantees that the optimal solution is indeed \replaced{non-negative}{nonnegative}. 
\end{myex}

Though the above examples are mostly concerned TV regularisation, they are readily extended to other types of regularisation such as Huber TV, or Wavelet regularisation (by replacing the gradient operator $\nabla$ with a wavelet transform $W$).

\subsection{Further Reading}

\paragraph{Advanced Modelling} The present survey focuses on the most commonly used mathematical models for variational regularisation.
However, this only scratches the surface of what can be done in this context. For example, a drawback of variational regularisation is the bias introduced by the regulariser, particularly when the data is not very informative. The bias can be reduced by the two-step technique in \cite{Brinkmann2016} or iterative regularisation \cite{Benning2018actanumerica}. Another interesting direction is the coupling of multiple measurements, which can correspond to different modalities or energy channels \cite{Ehrhardt2023, Arridge2021review} or by considering the problem as a dynamic inverse problem over time \cite{Schuster2018dynamicip, Hauptmann2021}. Such forward models require tailored mathematical modelling that can be done in the context of variational regularisation.

\paragraph{Machine Learning}
Recent years have seen a rise of learned regularisers for inverse problems \cite{Arridge2019}, for example dictionary learning \cite{tovsic2011dictionary}, fields-of-expert model \cite{roth2005fields}, gradient-step denoisers \cite{hurault2021gradient,tan2023provably}, regularisation-by-denoising \cite{romano2017little}, deep image prior with TV \cite{baguer2020diptv,barbano2024image}, generative regularisers \cite{duff2024regularising}, approaches based on the input-convex neural networks (ICNN) \cite{amos2017input} or network Tikhonov (NETT) \cite{li2020nett}. 
Among them, a focus is on learning informative yet convex regularisation, which 
utilise specially parameterised neural networks, such as the ICNN \cite{mukherjee2023learned} and convex ridge regularisers \cite{goujon2022neural}.

\subsection{Other Reviews}
To the best of our knowledge this is the first review that focuses on stochastic optimisation algorithms for inverse problems. However, there have been other reviews concerned with the efficient computation of inverse problems, some of which also focus on optimisation approaches. In \cite{Cullen2013} and \cite{Chung2015}, the authors discuss the efficient solution of large-scale inverse problems focusing on deterministic algorithms from numerical linear algebra. The reviews \cite{Burger2016firstorder},\cite{chambolle2016introduction}, and \cite{valkonen2016optimisingbig} are concerned with efficient first-order optimisation algorithms for inverse problems and imaging. While we also discuss first-order algorithms in detail, our particular focus is on stochastic approaches to overcome the large-scale nature of the problem. The paper \cite{Bellavia2023} surveys stochastic methods for training neural networks for image classification but does not consider inverse problems. The connection of convex optimisation to machine learning and AI for medical imaging is discussed in \cite{Xu2022}. They do touch upon stochastic optimisation but do not go into much detail or specifics that are of interest for inverse problems.

\section{Deterministic Algorithms for Variational Regularisation} \label{sec:optalgs}

\new{From the discussion in Section \ref{subsec:optim_template}, in the generic form we need to deal with the following minimisation problem
\begin{equation*}
\min_{\optvar \in \optspace} \Big\{  \Phi(\optvar) = f(\optop\optvar) + g(\optvar) + h(\optvar) \Big\},
\end{equation*}
under the Assumptions \ref{assump_1}-\ref{assump_3}. 
Developing efficient solvers to tackle variational optimisation problems has been an incredibly active and productive research area, especially since the new millennium. 
Among the existing methodologies, first-order methods, which utilise only the first-order derivative of the objective function, have been the most popular choice. Moreover, first-order methods are the backbone of the majority of the stochastic optimisation algorithms. 

Hence in this section we first present an overview of first-order methods, from foundational principles to recent algorithmic advancements. 
In Section~\ref{subsec:building_blocks} we first discuss the difficulties in solve generic variational optimisation problems to motivate the introduction of basic building blocks and design principles for first-order methods. 
In Section~\ref{subsec:advanced_1storder} we present representative first-order algorithms for minimising different variants of \eqref{eq:opttemplate}. We conclude in Section~\ref{sec:determinstic_further} with a discussion of further topics:  higher-order schemes, non-Euclidean methods, nonconvex optimisation and nonlinear problems. 

We emphasise that first-order optimisation is a rich and rapidly evolving field of research, consequently, our presentation is not comprehensive. 
We refer to monographs \cite{bauschke2011convex, beck2017first} and reviews \cite{LaurentCondat2023, Bachmayr2009, Benning2018actanumerica, Benning2021choose} for more details.

\subsection{Building Blocks of First-Order Algorithms}\label{subsec:building_blocks}

The challenges of solving \eqref{eq:opttemplate} in general are considered from the following sources: nonsmoothness, function-operator composition, the summation of several terms\new{, and even non-convexity}. 
The development of first-order optimization algorithms has been a reflection of the growing complexities of problems in inverse problems, computer science and engineering, as well as the increasing demand for computational efficiency. 

First-order optimization algorithms rely on gradient information to guide the iterative process of minimizing an objective function. This family of methods has evolved significantly over the decades, beginning with the simple yet foundational gradient descent algorithm and extending to more advanced techniques like momentum-based acceleration, proximal algorithms, and modern operator splitting techniques. Each stage in their evolution has been motivated by the need to address limitations of earlier approaches, adapt to new problem structures, and improve computational efficiency. 

Review the development history of first-order optimisation algorithms, the idea of 
{\it divide and conquer} has been consistently applied throughout. It has become a guiding  principle for designing algorithms to address the challenges.
It consists of decomposing the original problem into simpler subproblems, tackling each subproblem separately, and assembling all the components into a provable and implementable numerical scheme.}

\subsubsection{Basic Elements}\label{sec:fom+basic_element}

\new{
Over decades of development, techniques for decomposing the optimisation problem such as \eqref{eq:opttemplate} are well developed. In general, we can use the following rules
\begin{itemize}
    \item Directly decompose the objective: take \eqref{eq:opttemplate} for example, the problem can be decomposed into terms $h$, $g$ and $f\circ \optop$.
    \item Decompose based on optimality condition of optimal solution, this approach is similar to the first one, as the optimality in general is consists of the (sub)derivatives of each function. We refer to a recent review \cite{LaurentCondat2023} for dedicated discussion. 
    \item Decompose in terms of optimisation variables, this is particularly useful when the optimisation problem is separable in in the variables. 
\end{itemize}
Whichever decomposition approach, in the ``conquer'' step the procedure essentially boils down to managing $h$, $g$ and $f\circ \optop$. 
Therefore, in this part, based on the properties (Assumptions \ref{assump_1}-\ref{assump_3}) of the terms assigned in Section \ref{subsec:optim_template}, we discuss the basic approaches to tackle them individually. 
}

\paragraph{Gradient Descent} 
Gradient descent (GD) is one of the oldest iterative methods for unconstrained optimisation, dating back to at least the middle of the 19th century. When minimising a smooth differentiable objective $h$, GD iteratively performs a step in the direction of the negative gradient at the current iterate
\begin{equation}\label{eq:gradientdescent}
\xkp = \xk - \tau_k \nabla h(\xk),    
\end{equation}
where $\tau_k > 0$ is the step{\removed{-}}size. 
As the name suggests, GD is a descent method that reduces function values in each iteration $h(\xkp) \leq h(\xk)$ assuming that the step{\removed{-}}sizes $\tau_k$ are well-chosen\footnote{The optimisation literature often focuses on linear convergence rate for the strongly convex case. Here, we put more emphasis on the merely convex case, since in practice the optimisation problems associated to inverse problems are often only strongly convex with a very small convexity parameter.}. 
\new{When the function $h$ is $L$-smooth, i.e.\ $\nabla h$ is $L$-Lipschitz continuous, then under proper constant stepsize choice, GD obtains $\mathcal{O}(1/k)$ convergence rate for the objective function value \cite{beck2017first}.}

\paragraph{Proximal Point Algorithm}
\new{When minimising the \replaced{non-smooth}{nonsmooth} function $g$, one can still apply GD by replacing the gradient $\nabla h$ with a sub{\removed{-}}gradient (see \eqref{eqn:subdifferential} for a definition). 
However, the resulting sub{\removed{-}}gradient method is no longer a descent algorithm and does not necessarily converge under constant step{\removed{-}}size choice \cite{beck2017first}, moreover it can only obtain $\mathcal{O}(1/\sqrt{k})$ convergence rate which is worse than GD. As a consequence, a much more popular approach to deal with minimising nonsmooth objective is to use the proximal operator of $g$, see definition in \eqref{eqn:prox_operator}, which results in proximal point algorithm (PPA) \cite{rockafellar1976monotone}. 
}

For a proper, closed and convex function $g$ the corresponding proximal point algorithm, or the proximal minimisation algorithm, takes the form
\begin{equation}\label{eq:proximaldescent}
    \xkp 
    = \mathrm{prox}_{\tau_k g} (\xk)
    := \argmin_{\optvar\in \optspace} \left\{ \tau_k g(\optvar) + \frac{1}{2} \norm{\optvar-\xk}^2 \right\},
\end{equation}
where $\tau_k>0$ is the step{\removed{-}}size.
Convergence of PPA to a minimiser is ensured provided the step{\removed{-}}size does not vanish too quickly, for example $\sum_k \tau_k=\infty$. 
\new{PPA is also a descent method which is the same as GD \cite{beck2017first}, under constant stepsize choice it admits $\mathcal{O}(1/k)$ convergence rate for the objective function value.}

Note that each iteration requires solving an optimisation problem to compute the proximal operator.
However, the proximal operator of many widely used (regularisation) functions is known in closed-form \cite{combettes2011proximal, proximityoperatorProxRepository}. 

\begin{remark}\label{rmk:ode}
    GD and PPA can be interpreted as different discretisations of the gradient flow $x'(t) = - \nabla g(x(t)), t\geq 0$. Applying the explicit Euler scheme to the gradient flow yields GD \eqref{eq:gradientdescent}: $(\xkp-\xk)/\tau = -\nabla g(\xk)$.
    On the other hand, applying the implicit Euler scheme $(\xkp-\xk)/\tau = -\nabla g(\xkp)$, requires solving an equation $\xkp + \tau \nabla g (\xkp)=\xk$, to compute $\xkp$. 
    This in turn is equivalent to $\xkp=\mathrm{prox}_{\tau g} (\xk)$, which is PPA \eqref{eq:proximaldescent}.
    \new{When $g$ is nonsmooth, the implicit Euler scheme leads to solving the proximal operator of $g$.}
\end{remark}

\new{
\begin{remark}
Over the past two decades of development, especially the advances in nonsmooth optimisation, when dealing with convex optimisation problem, nonsmoothness nowadays is not necessarily much more challenging to deal with than smoothness: 1) Algorithm-wise, as long as the proximal operator of the nonsmooth objective is easy to compute, there is no difference compare to the smooth case. 2) Theory-wise, nonsmoothness introduces technicalities such as set-valued inclusion for optimality condition and argument in proving subgradient convergence. 
3) As both GD and PPA are descent algorithms, they have the same $\mathcal{O}(1/k)$ convergence rate for the objective function value. 
\end{remark}
}

\paragraph{Acceleration}
Acceleration is an essential component of any effective first-order method.
The most successful acceleration technique is \enquote{inertial acceleration} or \enquote{momentum}, dating back to the heavy-ball method \cite{polyak1964some} which can significantly speed up GD, but lacks global acceleration guarantees. 
Nesterov's accelerated gradient~(NAG) method was the first inertial acceleration that provably improves the objective function convergence rate of GD from \new{$\calO(1/k)$ to $\calO(1/k^2)$} 
for smooth convex objectives \cite{nesterov1983method}.
Fast iterative soft-thresholding algorithm (FISTA) extends this speed up to nonsmooth convex objectives \cite{beck2009fast}.
We refer to \cite{JEMS2023} for a more comprehensive treatment of the development of inertial acceleration over the past decade. 
Other acceleration techniques include preconditioning \cite{pock2011diagonal,liu2021acceleration}, vector extrapolation \cite{Anderson,scieur2016regularized} and so on. We refer to the monograph \cite{Aspremont2021Acceleration} on acceleration for further details.

\begin{remark}
Similar to GD, NAG also has an ODE interpretation. In \cite{sujmlr16}, the authors showed that the limit of NAG iterations is a second-order differential equation: $x''(t) + ({3}/{t}) x'(t) = - \nabla g (x(t))$.  
\end{remark}

    \begin{figure}[htbp]
    \centering
    \includegraphics[width=0.6\textwidth]{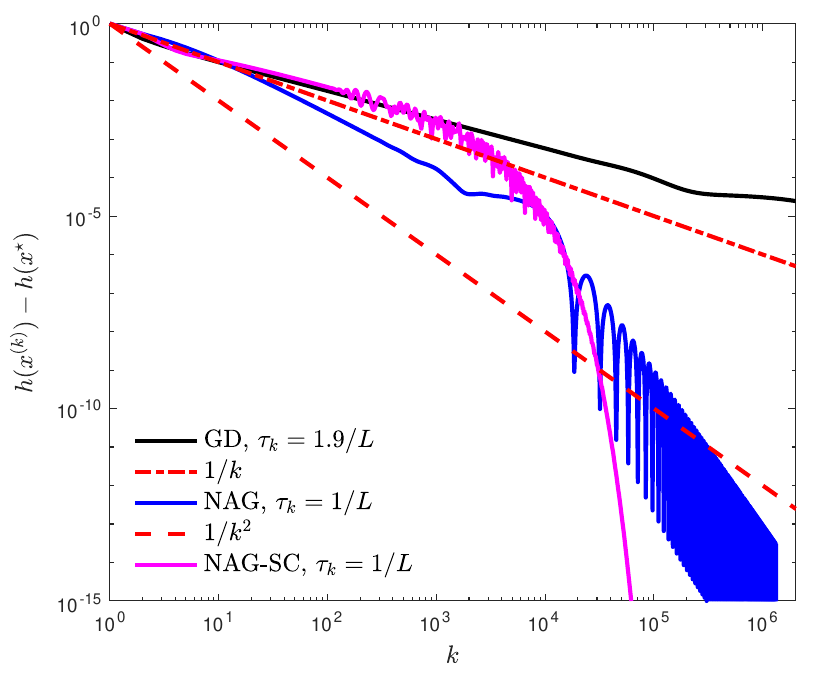}
    \caption{Comparison of GD, NAG on least-squares problem. }    
    \label{fig:cmp_gd_nag}
    \end{figure}
    
\new{
\begin{myex}\label{ex:gd_nag_ppa}
In this example, by considering least-squares $ \min_{\optvar}h(\optvar)=\frac12\norm{\ipop\optvar-\datavar}^2$, we provide a comparison on GD, NAG in standard form and NAG adapt to strong convexity (NAG-SC). We have $\ipop \in \mathbb{R}^{100\times 100}$ being a tri-diagonal matrix with diagonal elements equal to $2$ and off-diagonal elements equal to $-1$, $\datavar\in\mathbb{R}^{100}$ is a random vector. 
Under this setting, the problem is strongly convex with condition number of the order $10^7$. 
In Figure \ref{fig:cmp_gd_nag} we demonstrate the objective function value convergence: let $L=\norm{\ipop}^2$
\begin{itemize}
    \item GD with stepsize $\tau_k=1.9/L$ is the slowest, while NAG ($\tau_k=1/L$) achieves faster convergence rate. Note that the convergence rate of NAG is actually faster $1/k^2$ which is due to the fact the strong convexity. 

    \item 
    In NAG-SC we adapt the acceleration parameter according to the condition number which achieves even faster convergence rate. 

\end{itemize}
From the comparison, we can see that accelerated gradient method is less affected by the condition number of the objective function than gradient descent, and can achieve much faster convergence behavior. 
\end{myex}

\begin{remark}
In their standard form, GD, NAG and PPA are widely used optimisation algorithms with well-understood convergence properties. 
Although GD and PPA address objective functions with different smoothness properties, the two algorithms share the same convergence rate.
When using constant stepsizes for a convexobjective function, all algorithms achieve sublinear convergence rate. NAG improves the $\mathcal{O}(1/{k})$ objective function convergence rate of GD to $\mathcal{O}(1/k^2)$ which is the significance of the method. For strongly convex objective functions, all methods achieve a linear convergence rate. Between GD and NAG, NAG provides a faster linear rate under the same step size.
\end{remark}
}

\paragraph{Function-Operator Decomposition} 
Now we turn to the composition term $f\circ A$, whose proximal operator in general is difficult to compute, \new{we refer to \ref{app:prelims} for some exceptions}. 
\removed{A notable exception is when $\optop\optop^\ast=\alpha I,~\alpha>0$, in which case the proximal operator of $f \circ \optop$ can be computed as $\prox_{f\circ \optop}(x) = x + \tfrac{1}{\alpha} \optop^\ast \big( \prox_{\alpha f}(\optop x) - \optop x \big)$ \cite{beck2017first}.
For the general case} This \removed{issue can be} \replaced{circumvented}{is typically computed} by one of two popular approaches: duality or method of multiplier. 
The former uses the Fenchel--Moreau theorem, which states that the bi-conjugate (see \eqref{eqn:convex_conjugate} and \eqref{eqn:bi_conjugate} for definition) of a proper, closed convex function is equal to itself.
This gives
\begin{align}
f(\optop\optvar) = f^{\ast\ast} (\optop\optvar) := \sup_{\dualvar\in\dualspace} \big\{  \iprod{\optop\optvar}{\dualvar} -  f^\ast(\dualvar) \big\}, \label{eq:fenchelmoreau}
\end{align}
and thus the optimisation problem can be recast as a saddle-point problem, and then solved by the primal-dual algorithm PDHG (see Algorithm \ref{alg:pdhg}). 
An alternative is to decompose $f\circ \optop$ via an auxiliary variable into a constrained problem
\begin{align*}
\min_{\optvar \in \optspace,\ \auxvar\in \auxspace } \big\{ f(\auxvar) \quad \text{ such that }\quad \auxvar = \optop\optvar \big\}.
\end{align*}
By adding the constraint back to the objective via a Lagrange multiplier, the problem can be solved by dual ascent, method of multipliers \cite{boyd2011distributed}, or the alternating direction method of multipliers (ADMM) (see Algorithm \ref{alg:admm}). PDHG and ADMM are deeply connected, and equivalent to each under some conditions, see \cite{chambolle11pdhg, chambolle2016introduction, o2020equivalence} for discussions.

\paragraph{Coordinate Descent} 
Instead of updating all the entries of the optimisation variable $\optvar$ in each iteration, coordinate descent (CD) methods update only one or a small subset of coordinates:
\begin{equation}\label{eq:descentscheme}
\begin{aligned}
\xkp 
= \xk - \tau_k \nabla_{B_k} h(\xk)e_{B_k}
= \xk - \tau_k \sum_{i\in B_k} \nabla_{i} h(\xk)e_{i}.
\end{aligned}
\end{equation} 
Here $B_k$ is a subset of the coordinates of $\optvar$, $e_i$ is the $i$-th canonical basis vector and $\nabla_{i} h(\xk) = \frac{\partial h(\xk)}{\partial \optvar_i} $ denotes the directional derivative along $i$-th coordinate. 
When $B_k$ is a singleton, \eqref{eq:descentscheme} is the basic CD scheme, otherwise it is called {\it block} coordinate descent. CD-based methods are commonly used for problems that have a separable structure or that are coordinate-friendly \cite{shi2016primer}. 
Popular deterministic choices for the order in which the coordinates are updated include cyclic and fixed permutation orders. Randomised orders are discussed in Section~\ref{subsec:scd}. 
The introduction of CD to optimisation is due to \cite{bertsekas2015}\footnote{The first edition of the book was published in 1989.}. 
Convergence results were first shown for differentiable objectives \cite{luo1992convergence}, and later extended to \replaced{non-differentiable}{nondifferentiable} objectives \cite{tseng2001convergence}. 
We refer to \cite{nesterov2012efficiency,beck2013convergence,richtarik2014iteration} for theoretical advances of CD algorithms.

\subsubsection{Designing a First-Order Algorithm} 

With the above basic elements at hand, designing a first-order methods consists of the following steps: 
\begin{itemize}
    \item {\bf Problem decomposition} The original optimisation problem is decomposed into simpler subproblems, according to properties of the individual functions. 

    \item {\bf Subproblem processing and assembly} 
    Each subproblem is assigned an approach to tackle it (using e.g.\ Table \ref{tab:elements} as a guideline). The resulting components are then assembled to form an iterative scheme.
    In principle, optimality conditions can serve as a reference for the assembly of the algorithm.

    \item {\bf Refinement} The iterative scheme might require further improvements to address an existing issue. 
    First, one or several of its components might require solving an additional optimisation problem, e.g.\ for the computation of the proximal operator.
    Second, theoretical convergences guarantees might be either lacking or too restrictive, and thus the approach ought to be modified.
    Third, theoretical or practical convergence may be slow, and thus require acceleration or other techniques to improve the performance.

\end{itemize}

\begin{table}[htbp]\renewcommand\arraystretch{1.35}
\centering
\begin{tabular}{c@{\hspace{4mm}}c@{\hspace{4mm}}c}
\toprule
Property & Method & References \\ \midrule \midrule
Smoothness & $\xkp=\xk-\tau_k\nabla h(\xk)$ & \cite{kantorovich1960,bauschke2011convex} \\ \hline
Prox-friendly &  $\xkp = \mathrm{prox}_{\tau_k g}(\xk)$ & \cite{moreau1965,parikh2014} \\ \hline
\multirow{2}{*}{Composition} & Duality: $f(\optvar) = \sup_{\dualvar} \left\{\iprod{\dualvar}{\optop \optvar} - f^\ast(\dualvar)\right\} $ & \cite{Chambolle2004, chambolle11pdhg} \\ 
& Constraint: $\min_{\optvar,\auxvar} \big\{ f(\auxvar) ~~~\textrm{s.t.}~~ \auxvar = \optop\optvar \big\}$ & \cite{gabay1976dual,boyd2011distributed} \\ \hline 
Coordinate friendly & $\xkp = \xk - \tau_k \nabla_{i_k} h (\xk) e_{i_k}$  &      \cite{bertsekas2015,wright2015coordinate} \\ \bottomrule
\end{tabular}
\caption{Foundational elements in designing first-order algorithms} \label{tab:elements}
\end{table} 

It is important to emphasise that the above steps are not necessarily unique.
In the following we present classical first-order algorithms that can be applied to solve instances of \eqref{eq:opttemplate}. We start with simpler cases where only two of the terms in \eqref{eq:opttemplate} are nontrivial, and then move to the more general three-term case.

\subsection{Advanced First-Order Algorithms}\label{subsec:advanced_1storder} 

In this section we present representative first-order algorithms for minimising two-term and three-term optimisation problems. 

\subsubsection{Approaches for a Two-Term Sum}\label{sec:twotermsum} 

Minimisation of a two-term sum of functions is widely studied with several popular approaches in the literature. 
Moreover, approaches for two-term problems often serve as building blocks for more advanced methods. 

\paragraph{Proximal Gradient Descent}  
When \replaced{$f$ vanishes}{$f=0$}, \eqref{eq:opttemplate} reduces to a two-term sum
\begin{equation}
\label{eq:var_model_pgd}
\min_{\optvar \in \optspace} \left\{ \Phi(\optvar) = g(\optvar) + h(\optvar)\right\},
\end{equation}
where $h$ is smooth and $g$ is prox-friendly.
The most popular scheme to solve \eqref{eq:var_model_pgd} is the proximal gradient descent (PGD), see Algorithm~\ref{alg:pgd}, also known as the forward-backward splitting method~\cite{2005_Combettes_Signal}.  
PGD performs a gradient descent step on the smooth term $h$, and then evaluates the proximal operator of $g$.

\begin{center}
\begin{minipage}{0.95\linewidth} 
\begin{algorithm}[H] 
\caption{Proximal Gradient Descent (PGD)}\label{alg:pgd}
\DontPrintSemicolon
\KwIn{$x^{(0)} \in \optspace$, $\tau > 0$.}
\For{$k=0,1,...$}{
\(\tilde x^{(k)}=\) \mybox[red]{\(\xk - \tau \nabla h(\xk)\)} \hfill{\scriptsize \texttt{gradient descent on smooth term}}\\
\(\xkp =\) \mybox[green]{\(\mathrm{prox}_{\tau g}(\tilde x^{(k)})\)} \hfill{\scriptsize\texttt{proximal descent on nonsmooth term}}
}
\end{algorithm}
\end{minipage}
\end{center}

PGD is the simplest extension of GD to the nonsmooth case, and it shares similar convergence properties:    
provided that $h$ is $L$-smooth and the step{\removed{-}}size satisfies $\tau < 2/L$, PGD converges at a rate $\calO(1/k)$ in the objective \cite{beck2009fast}\footnote{The optimisation literature often considers the strongly convex case, too. We omit it here since in practice the optimisation problems associated to inverse problems are merely strongly convex}.   

\paragraph{Accelerated PGD}
In spite of its simplicity, PGD is in practice hindered by its slow $\calO(1/k)$ convergence rate. 
Fast iterative soft-shrinkage algorithm (FISTA) \cite{beck2009fast}, extends the acceleration methodology of NAG to the \replaced{non-smooth}{nonsmooth} convex case.

\begin{center}
\begin{minipage}{0.95\linewidth} 
\begin{algorithm}[H] 
\caption{Fast Iterative Soft-shrinkage Algorithm (FISTA)}\label{alg:fista}
\DontPrintSemicolon
\KwIn{$x^{(0)} = x^{(-1)} \in \optspace$, $\tau >0$, $t_0=1$.}
\For{$k=0,1,...$}{
$t_{k+1} = \frac{1+\sqrt{1+4t_{k}^2}}{2},~ a_k = \frac{t_{k}-1}{t_{k+1}}$\\
$\tilde x^{(k)} =$  \mybox[orange]{\(\xk + a_k( \xk - \xkm)\)} \hfill{\scriptsize\texttt{momentum}}\\
\(\xkp =\)  \mybox[green]{\(\mathrm{prox}_{\tau g}\big( \myinbox[red]{\( \tilde x^{(k)} - \tau \nabla h(\tilde x^{(k)})\)}\big)\)}\hfill{\scriptsize\texttt{proximal descent on nonsmooth term}}
}
\end{algorithm}
\end{minipage}
\end{center}

Provided that $h$ is convex and $L$-smooth and  $\tau\leq1/L$, FISTA improves the convergence rate in the objective value to  $\calO(1/k^2)$. 
It can be shown that the iterates of FISTA converge to a global minimiser with a modified update rule $a_k = (k-1)/(k+d), d>3$ \cite{chambolle2015convergence}. 
Moreover, the authors also show that the objective function convergence speed is indeed $o(1/k^2)$. 

Despite the theoretical advantage, in practice FISTA is not necessarily faster than PGD, especially when high precision solution is sought. The reason is that FISTA is not a descent method and shows (periodic) oscillatory behaviour, typically when the problem is (locally) strongly convex. To circumvent this problem, restarting schemes can significantly improve the performance of FISTA \cite{o2015adaptive,sujmlr16}. 

\paragraph{Primal-Dual Splitting Methods} 
When \replaced{$g$ vanishes}{$h=0$}, \eqref{eq:opttemplate} reduces to
\begin{equation}\label{eq:var_model_fAh}
    \min_{\optvar\in \optspace} \left\{ \Phi(\optvar) = f(\optop\optvar) + \replaced{h}{g}(\optvar)\right\} , 
\end{equation}
where $f$ and $g$ are both nonsmooth but prox-friendly. 
As previously discussed, the proximal operator of $f\circ \optop$ is in general not available in closed form and alternative approaches are needed to handle $f\circ \optop$.

Using the concept of duality, plugging \eqref{eq:fenchelmoreau} into \eqref{eq:var_model_fAh} gives the saddle point problem
\begin{equation}\label{eqn:saddle_point_problem}
\min_{\optvar\in\optspace} \max_{\dualvar\in\dualspace} \left\{ \replaced{h}{g}(\optvar) + \iprod{\optop\optvar}{\dualvar} - f^\ast(\dualvar)\right\} .
\end{equation}
Chambolle and Pock \cite{chambolle11pdhg} propose to solve \eqref{eqn:saddle_point_problem} via an algorithm that alternates between a descent in the primal and an ascent in the dual on an extrapolated variable, see Algorithm \ref{alg:pdhg} below for details.

\begin{center}
\begin{minipage}{0.95\linewidth}
\begin{algorithm}[H]         
    \caption{Primal-Dual Hybrid Gradient Algorithm (PDHG)}\label{alg:pdhg}
        \DontPrintSemicolon
        \KwIn{$\optvar^{(0)} \in \optspace, \dualvar^{(0)} \in \dualspace$, $\sigma, \tau > 0$.}
        \For{$k=0,1,\ldots$}{
            $\xkp = $ \mybox[green]{\(\mathrm{prox}_{\tau g} (             \textrm{\myinbox[red]{$ \xk - \tau \optop^\ast \yk $}})\)} \hfill{\tiny\texttt{proximal descent on nonsmooth term}}\\
            $\bar{x}^{(k+1)} =$ \mybox[orange]{\(2\xkp - \xk\)} \hfill{\tiny\texttt{extrapolation}}\\
            $\ykp =$ \mybox[blue]{$\mathrm{prox}_{\sigma f^\ast} \big( \textrm{\myinbox[red]{\(\yk + \sigma \optop \bar{\optvar}^{(k+1)}\)}} \big)$} \hfill{\tiny\texttt{proximal ascent and dualisation}}
        }
\end{algorithm}
\end{minipage}
\end{center}

This algorithm is known as the primal-dual hybrid gradient (PDHG) or the Chambolle--Pock algorithm. 
PDHG can be equivalently written as PPA on $\optspace \times \dualspace$ \cite{he2012}. Hence, one can obtain the convergence of $(\xk,\yk)$ to a saddle point $(x^\ast, y^\ast)$ where  $x^\ast$ is a solution of \eqref{eq:var_model_pgd}. 
The most known convergence result for PDHG is that the primal-dual gap converges as $\calO(1/k)$ 
provided $\sigma\tau\|\optop\|^2 < 1$ \cite{chambolle11pdhg}. 
In a recent work \cite{Ma2023} improved the step{\removed{-}}size condition to $\sigma\tau\|\optop\|^2<4/3$. 
We also refer to~\cite{goldstein2013adaptive,goldstein2015adaptive} for alternative adaptive step{\removed{-}}size schemes. 

\paragraph{Alternating Direction Method of Multipliers} 
Addressing $f\circ \optop$ through the method of multipliers
converts \eqref{eq:var_model_fAh} into a constrained problem
\begin{equation}\label{eq:problem-admm}
\min_{\optvar \in \optspace,\ \auxvar\in \auxspace } \left\{\Phi(\optvar,\auxvar):= f(\auxvar) + g(\optvar)\right\} \quad 
\textrm{ such that} \quad \optop\optvar - \auxvar = 0.
\end{equation}
The objective $\Phi$ is separable in $\optvar$ and $\auxvar$, and the constraint enforces their coupling.
The associated augmented Lagrangian reads 
\begin{equation}\label{eq:Lagrangian}
\calL(\optvar,\auxvar; \dualvar) = f(\auxvar) + g(\optvar) + \iprod{\dualvar}{\optop\optvar-\auxvar} + \frac{\tau}{2}\|{\optop\optvar-\auxvar}\|^2
\end{equation}
where $\tau>0$ and $w$ is the Lagrangian multiplier. 
To find a saddle-point of $\calL(x,z; w)$, the Alternating Direction Method of Multipliers (ADMM, Algorithm \ref{alg:admm}) updates $z, x$, and $w$ one at a time.  

\begin{center}
\begin{minipage}{0.95\linewidth}
\begin{algorithm}[H]         
    \caption{Alternating Direction Method of Multipliers (ADMM)}\label{alg:admm}
        \DontPrintSemicolon
        \KwIn{$\optvar^{(0)} \in \optspace, \dualvar^{(0)} \in \dualspace, \tau > 0$.}
        \For{$k=0,1,...$}{
            $\zkp \in \argmin_{\auxvar\in \auxspace}~ \calL(\xk,\auxvar,\yk)$ \hfill{\tiny\texttt{updating the dual variable}}\\
            $\xkp \in \argmin_{\optvar\in\optspace}~ \calL(\optvar,\zkp,\yk)$ \hfill{\tiny\texttt{updating the primal variable}}\\
            $\ykp = \yk + \tau ({\optop\xkp  - \zkp})$
        }
    \end{algorithm}
    \end{minipage}
    \end{center}
We refer to \cite{glowinski1975approximation,gabay1976dual} for the early development of ADMM and \cite{boyd2011distributed} for a more modern survey. 
ADMM is equivalent to applying the Douglas--Rachford splitting method to the dual of \eqref{eq:var_model_fAh}, see \cite{eckstein1992douglas,combettes2007douglas}.

The update for the auxiliary variable can be written as $\zkp=\prox_{f/\tau}( \optop\xk - \auxvar + \tfrac{1}{\tau}\yk)$, and is thus uniquely determined.
On the other hand, the minimisers of $\calL(\optvar,\zkp,\yk)$, and thus $\xkp$, need not to be unique, unless $\optop$ has full column-rank. 
The convergence analysis of ADMM is in general complicated.
However, when $f$ and $h$ are closed and convex, and $\calL$ admits a saddle point, convergence in the objective and the constraint can be shown \cite{boyd2011distributed}.
The convergence of the primal sequences $\xk$ and $\zk$ is more difficult to establish but when $\optop$ is of full column-rank, the convergence of $(\xk,\zk)$ can be shown \cite{eckstein1992douglas}. 

\begin{remark}\label{remark:admm2}
Instead of updating $\zkp$ and $\xkp$ sequentially, as in Algorithm \ref{alg:admm}, we can update them simultaneously
$$(\zkp,\xkp) \in \argmin_{\optvar\in\optspace,\auxvar\in \auxspace}~ \calL(\optvar,\auxvar,\yk). $$ 
The resulting method is called the augmented Lagrangian method \cite{beck2017first} or method of multipliers \cite{boyd2011distributed}. However, the simultaneous update has in general no closed-form expression. 
Hence, sequential (block) updates are preferred. This can be readily extended to problems with more than two blocks; see Section~\ref{sec:threetermsum} for a discussion.     
\end{remark}

\begin{remark}\label{remark:admm}
    ADMM is traditionally derived to solve a variation of the problem \eqref{eq:problem-admm} with a more general constraint: $Ax+Bz=u$.  This requires only mild modifications of the Algorithm \ref{alg:admm}. However, the first step is no longer a proximal operator, but requires minimising $f(z) + \frac{\tau}{2} \|{A\xk +Bz-u + \tfrac{1}{\tau}\wk }\|^2$, if $B^\ast B\neq I$. To solve this subproblem, there are two standard strategies. First strategy is to add a proximal term $\iprod{z-\zk}{Q(z-\zk)}$ with $Q = \frac{1}{\beta} I - B^\ast B$ where $\beta>0$ is such that $Q$ is symmetric positive definite. In doing so, the problem reduces to computing the proximal operator of~$f$.
    The second strategy is to linearise the quadratic term at point $\zk$, which corresponds to applying one-step PGD to solve the corresponding subproblem. 
\end{remark}

\subsubsection{Approaches for a Three-Term Sum}\label{sec:threetermsum}
In previous section we discussed standard approaches for an objective decomposed into a sum of two terms. 
Most of the approaches presented therein can be extended to the three-term optimisation template \eqref{eq:opttemplate}. However, sometimes more flexibility is needed. 

\paragraph{Condat--V\~u Algorithm}
Condat and V\~u independently developed a primal-dual forward-backward algorithm for \eqref{eq:opttemplate} as a generalisation of the PDHG, see Algorithm~\ref{alg:condatvu}. 

\begin{center}
\begin{minipage}{0.95\linewidth} 
    \begin{algorithm}[H] 
        \caption{Condat--V\~u Algorithm}\label{alg:condatvu}
        \DontPrintSemicolon
        \KwIn{$\optvar^{(0)} \in \optspace, \dualvar^{(0)} \in \dualspace$, $\tau, \sigma > 0$}
        \For{$k=0,1,\ldots$}{
            $\xkp = $ \mybox[green]{$\prox_{\tau g}\big(\text{\mybox[red]{$\xk - \tau \nabla h(\xk) - \tau \optop^* \dualvar^{(k)}$}}\big)$} \hfill{\tiny\texttt{proximal descent on nonsmooth term}}\\
            $\overline{x}^{(k+1)} = $ \mybox[orange]{$2\xkp - \xk$} \hfill{\tiny\texttt{extrapolation}}\\
            $\dualvar^{(k+1)} = $ \mybox[blue]{$\prox_{\sigma f^\ast}\left( \text{\myinbox[red]{\(\dualvar^{(k)} + \sigma \optop \overline{x}^{(k+1)}\)}} \right)$} \hfill{\tiny\texttt{proximal ascent and dualisation}}
        }
    \end{algorithm}
\end{minipage}
\end{center}

The iterates of the Condat--V\~u algorithm converge if $\tau(\sigma\|\optop\|^2 + L/2)<1$ where $L$ is the Lipschitz constant of the gradient of $h$ \cite{condat2013primal,vu2013splitting}. Note that each step requires one evaluation of $\optop, \optop^*, \prox_{\tau g}, \prox_{\sigma f^\ast}$ and the gradient $\nabla h$.

\paragraph{PD3O}
The primal-dual three operator algorithm (PD3O) is another generalisation of PDHG and other two- and three-term algorithms. The difference between PD3O and Condat--V\~u is in the final term in the update for $y^{(k+1)}$, where PD3O adds the gradient momentum term $\tau(\nabla h(\xk)-\nabla h(\xkp))$.

\begin{center}
\begin{minipage}{0.95\linewidth} 
    \begin{algorithm}[H] 
        \caption{Primal-Dual Three Operator Algorithm (PD3O)}\label{alg:pd3o}
        \DontPrintSemicolon
        \KwIn{$\optvar^{(0)} \in \optspace, \dualvar^{(0)} \in \dualspace$, $\sigma,\tau >0$ }
        \For{$k=0,1,\ldots$}{
                    $\xkp = $ \mybox[green]{$\prox_{\tau g}\big(\text{\mybox[red]{$\xk - \tau \nabla h(\xk) - \tau \optop^* \dualvar^{(k)}$}}\big)$} \hfill{\tiny\texttt{proximal descent on nonsmooth term}}\\
        $\overline{\optvar}^{(k+1)}=$ \mybox[orange]{$2\xkp-\xk+\tau(\nabla h(\xk)-\nabla h(\xkp))$}\hfill{\tiny\texttt{extrapolation}} \\
        $\dualvar^{(k+1)} = $ \mybox[blue]{$\prox_{\sigma f^\ast}\big(\text{\myinbox[red]{$\dualvar^{(k)} + \sigma \optop \overline{x}^{(k+1)}$}} \big)$} \hfill{\tiny\texttt{proximal ascent and dualisation}}\\
        }
    \end{algorithm}
\end{minipage}
\end{center}
PD3O requires one more evaluation of the gradient $\nabla h$ than Condat--V\~u, however, these costs can be amortised by storing $\nabla h(\xkp)$ and then using it in the next iteration. 

It has been shown to converge for $\sigma\tau\|\optop\|^2 < 1$ and $\tau L<2$ \cite{yan2018pd3o}.
Importantly, the convergence condition decouples the conditions coming from PDHG and GD. This is in contrast to the Condat--V\~u algorithm where the convergence criteria couples the parameters $\tau$ and $\sigma$, with both $\|\optop\|^2$ and $L$.

\paragraph{Further Splitting Approaches}
Methods presented in previous sections are only a small selection of  first-order schemes for solving various instances of the optimisation template \eqref{eq:opttemplate} but listing all the existing algorithms is beyond the scope of this paper. 
Other important algorithms for the two-term sum are the Douglas--Rachford splitting algorithm \cite{lions1979splitting,combettes2007douglas} and the Loris--Verhoeven algorithm \cite{loris2011lv} (also known as proximal alternating predictor corrector algorithm).
For the three-term sum they include primal-dual fixed point algorithm (PDFP) \cite{chen2016pdfp}, primal-dual Davis--Yin algorithm~(PDDY)~\cite{davis2016davisyin} and the asymmetric forward-backward adjoint splitting (AFBA)~\cite{latafat2017afba}.
The existing two- and three-term algorithm are often intimately connected. For example, with specific parameter and function choices PD3O is a generalisation of Loris--Verhoeven, Davis--Yin, PDHG and Douglas--Rachford algorithms.
We refer to \cite{Jiang2023, LaurentCondat2023} for further connection between these algorithms. 
ADMM can also be extended to solve three-term sum problems, however, as shown in \cite{chen2016direct}, the resulting algorithm is not necessarily convergent, and thus requires modifications.

Numerical approaches for optimisation problems with more than three terms (sometimes called the multi-block problem) have also been studied in the literature. In \cite{raguet2013generalized} the authors proposed a so-called ``generalised Forward-Backward splitting'' algorithm for the minimisation of a smooth function and several prox-friendly terms. 
There are several primal-dual splitting methods for solving multi-block problems, e.g.~\cite{combettes2012primal,vu2013splitting,combettes2014variable}.  
A widely adopted strategy is block-wise updates, which shares the same spirit as coordinate descent; see also the discussion in Remark \ref{remark:admm2}.

\subsection{Further Reading}\label{sec:determinstic_further}
Here we discuss advanced topics related to presented approaches, including: higher-order schemes, Bregman methods, incremental algorithms and so on.

\paragraph{Second-Order Methods}
First-order methods have several drawbacks.
For example, they often exhibit slow convergence speed, are not scale invariant, and often require tedious parameter tuning. 
Higher-order methods such as Newton's method address some of the issues by using higher-order derivatives to estimate the objective \cite{nocedal2006numerical}. 
In many imaging applications, computing the Hessian can be prohibitively expensive.  
Quasi-Newton methods try to overcome this problem by approximating the Hessian using only gradient information.
Perhaps the most-known quasi-Newton method is the Broyden--Fletcher--Goldfarb--Shanno (BFGS) algorithm \cite{broyden1970convergence,fletcher1970new,goldfarb1970family,shanno1970conditioning}.
BFGS achieves superlinear convergence under suitable convexity and smoothness conditions, and it has been applied to nonsmooth problems  \cite{lewis2013nonsmooth,yu2008quasi}.
However, its memory requirements can be high, which has led to the development of the limited-memory BFGS (L-BFGS) \cite{liu1989limited}, that has been applied to several inverse problems~\cite{vogel2000lbfgsip, tsai2018fqn}.

The superlinear convergence of Newton methods can be extended to nonsmooth optimisation problems for the class of so called semi-smooth functions \cite{qi1993semismooth}. Semi-smooth Newton methods have been applied for Lasso problems \cite{li2018highly}, and have been extended to PGD and Douglas--Rachford splitting type algorithms \cite{xiao2018regularized, hu2022local,liang2023squared}, primal-dual algorithms \cite{wang2023quasi}, to name a few, with applications in TV minimisation \cite{ng2007ssqntv} and seismic tomography \cite{ulbrich2015ssqnseismic}.

\paragraph{Bregman Proximal Methods} 
Most of this review is restricted to the Hilbert space / Euclidean setting. Generalisations are an active field of research, often based on the Bregman distance \cite{bregman1967relaxation}. 
Let $\phi: \optspace \to \mathbb{R}$ be a strictly convex and continuously differentiable function. The Bregman distance  associated to $\phi$ is defined as $D_\phi(x, y) = \phi(x) - \phi(y) - \langle \nabla \phi(y), x - y \rangle$. 
The Bregman distance is in general not symmetric and does not satisfy the triangle inequality. However, it is \replaced{non-negative}{nonnegative} and only equals to $0$ if and only if $x=y$. 
The notion of $L$-smoothness can be generalised to ``relatively smooth'' associated to a Bregman distance \cite{bauschke2017descent,teboulle2018simplified}.
Various Bregman first-order methods can be designed in this framework, such as Bregman proximal algorithms \cite{bauschke2017descent}, Bregman ADMM \cite{wang2014bregman}, Bregman PGD \cite{mukkamala2019beyond,hanzely2021accelerated}, and Bregman Primal-Dual splitting \cite{chambolle2022accelerated,jiang2022bregman}.  

\paragraph{\replaced{Nonconvex Optimisation}{Nonlinear Inverse Problems}}
While many applications can be well modelled with a linear forward model, there are also many applications where a linear model is insufficient, e.g.\ electrical impedance tomography (EIT) \cite{Cheney1999}, phase retrieval for imaging~\cite{shechtman2015phase}, photoacoustic tomography \cite{wang2015photoacoustic}, full waveform inversion~\cite{tarantola1984inversion} with applications in geophysics, \cite{virieux2009fullwaveform}, optical imaging generally \cite{arridge1997optical} and diffuse optimal tomography specifically \cite{boas2001imaging}, to name just a few.
From an optimisation perspective this is challenging as the resulting minimisation problem is generally nonconvex, its regularity properties are heterogeneous and most algorithms are only guaranteed to converge to locally optimal solutions.
Some convex optimisation algorithms have been extended to nonconvex optimisation and applied to inverse problems \cite{forsgren1998primal,conn2000primal,valkonen2014primal, valkonen2021first,chen2021non}.
\new{Recent years have seen much progress on the theoretical foundations of optimisation algorithms for nonconvex problems by exploiting properties like the Polyak/Kurdyka-\L{}ojasiewicz condition \cite{attouch2010proximal, bolte2014proximal, karimi2016linear, fatkhullin2022sharp}. Identifying which nonlinear inverse problems meet this condition remains an open research question.}

\paragraph{Learning-to-Optimise}
With the emergence of deep learning, learning-to-optimise (L2O) is becoming a promising direction in optimisation \cite{li2016learning, chen2022learning}. This line of work was in its early days aimed at accelerating the training of machine learning models, focusing mostly on improving the network training, compared to classical hand-crafted optimisers \removed{such as SGD and Adam} using reinforcement learning \cite{andrychowicz2016learning}. However, most such methods lack convergence guarantees which motivated research on provably convergent L2O schemes, which aim to learn parts of the algorithm such as the step{\removed{-}}sizes for PDHG \cite{banert2020data} and the mirror maps of the mirror descent algorithms~\cite{tan2023data,tan2023boosting}, with applications in both machine learning and inverse imaging problems. In~\cite{ehrhardt2024learning}, the authors developed a framework for greedy L2O with guaranteed convergence on the training data. Theoretical results of the generalisation guarantees for L2O schemes were recently developed \cite{sucker2024learning}, leading a new dimension of research in this area.

\new{\paragraph{Iterative Regularisation}
Thus far we have discussed solving inverse problems via variational regularisation \eqref{eq:varreg} with optimisation methods. An alternative is to not include a regulariser directly and instead stop the iterations early. For example, the Landweber iterations apply gradient descent directly on the data fidelity. With well-chosen stepsizes, number of iterations and initialisation, the Landweber method can be shown to converge to the minimal-norm solution of \eqref{eq:ip}, if one exists \cite{landweber1951iteration}.
Choosing a good stopping criterion is an important and task dependent problem.
Landweber iterations for a nonlinear forward operator $K$ were analysed in \cite{hanke1995convergence}.  
Turning an optimisation method into a regularisation method is commonly referred to as iterative regularisation, which has been used for linearised Bregman iterations \cite{Yin2008, Yin2010} and their nonconvex variants \cite{Benning2021choose}, dual algorithms \cite{garrigos2018iterative}, primal-dual algorithms \cite{Molinari2020}, and analysed in Banach spaces \cite{schopfer2006banachlandweber}, to name a few.}

\section{Stochastic Algorithms for Variational Regularisation} \label{sec:stoch}
The first-order numerical schemes discussed in Section \ref{sec:optalgs} have been the predominant methods for solving inverse problems over the past decades, due to their efficiency, simplicity, and relatively fast convergence speed. However, due to rapid increases in problem complexity and data sizes over the recent years, the computational costs associated with first-order methods is becoming prohibitively high. This is the main driving force for the development of algorithms with a low per-iteration computational cost, which are often stochastic in nature. In this section we present an overview of existing stochastic optimisation approaches particularly relevant for inverse problems. The layout of the section is similar to that of Section~\ref{sec:optalgs}. 
In Section \ref{subsec:stoch_design} we discuss the basic ingredients of randomisation for first-order optimisation algorithms, and in Sections \ref{subsec:SGM} and  \ref{subsec:scd} we introduce several representative algorithms.
We conclude in Section \ref{subsec:further_stochastic}.

\subsection{Designing First-Order Stochastic Algorithms}\label{subsec:stoch_design}

When designing stochastic optimisation algorithms we aim to reduce the per-iteration computational costs while not degrading the performance too much, thus creating an overall faster and more efficient approach. 
In order to discuss the design principles of stochastic algorithms, we first rewrite problem \eqref{eq:opttemplate}, revealing a finite-sum structure,
\new{\begin{equation}\label{eq:opttemplate_finite_sum}
\min_{\optvar \in \optspace} \left\{\Phi(\optvar) = \sum_{i=1}^\numf f_i(\optop_i \optvar) + \sum_{i=1}^\numg g_{\new{i}}(\optvar) + \sum_{i=1}^\numh h_i(\optvar)\right\}.
\end{equation}}
For all $i$ the operators $\optop_i:\optspace\rightarrow\dualspace_i$ are linear and bounded, $f_i:\dualspace_i\rightarrow\R_\infty$ and $g_i :\optspace\rightarrow\R_\infty$ are proper, closed, convex, and prox-friendly, and $h_i:\optspace\rightarrow\R$ are convex and $L_i$-smooth. 
Moreover, we denote $\dualspace=\dualspace_1\times\ldots\times\dualspace_\numf$. 
The analysis of stochastic algorithms for this problem relies on several key quantities, the most used one is the largest Lipschitz constant of the smooth terms, $\Lmax=\max\{L_1,\ldots,L_\numh\}$.  

\begin{remark}
    The finite-sum structure in \eqref{eq:opttemplate_finite_sum} naturally arises in inverse problems when considering MAP estimators, see Remark \ref{rem:bayesian_variational}. In contrast, in machine learning, the loss function often represents an \emph{empirical risk} which consists of the average loss, i.e.\ each sum is divided by the number of summands. While mathematically equivalent, this rescaling has consequences for some formulas that follow.

    \new{The empirical risk is an approximation to the population risk where the objective is the expected value of a random process. This situation naturally arises in machine learning but is less common in inverse imaging problems. A few exceptions are list-mode PET imaging and Cryo-EM where the data can be seen as samples from a random direction. To the best of our knowledge this connection has not been exploited yet.}
\end{remark}

\begin{myex}\label{eg:stochastic_CT}
    Let $\ipop$ be a forward operator, modelled by the discrete X-Ray transform with a parallel-beam geometry and $180$ angles from the $[0, \pi)$ range, $\reconvar^\dagger$ the object of interest, and $\datavar$ the observed sinogram (encoding the noise as well). 
    The forward operator~$\ipop$, object $\reconvar^\dagger$, and data $\datavar$ are approximately related via the forward model $\ipop\reconvar^\dagger=\datavar$.
    To recover $\reconvar^\dagger$ from $\datavar$ we consider the variational problem
    \begin{align}\label{eqn:radon_varreg}
\argmin_{\optvar\in\X} \Big\{\Phi(\optvar)=\frac{1}{2} \|\ipop\optvar-\datavar\|_2^2 +\alpha\|\nabla\optvar\|_1 + \iota_{[0,\infty)^d}(\optvar)\Big\},
\end{align}
where $\alpha>0$. 
The objective in \eqref{eqn:radon_varreg} can be written in several ways to deal with the costs of evaluating $\ipop$ and $\ipop^\ast$.
For the following, let $\ipop_i$, for $1\leq i\leq 30$ be the CT Radon operator that uses $6$ angles from the range $[\pi i/180, \pi)$ with increments of size $\pi/30$.
The data $\datavar$ is split in a corresponding fashion, giving subset data vectors $\datavar_i$ that use $1/30$-th of the data in $\datavar$.
We can then identify functions in \eqref{eqn:radon_varreg} with the template \eqref{eq:opttemplate_finite_sum} as
\begin{itemize}
    \item $g(\optvar)=\alpha \|\nabla \optvar\|_1+\iota_{[0,\infty)^d}(\optvar)$, $f\equiv0$, and $h_i(\optvar)= \frac{1}{2}\|\ipop_i\optvar-\datavar_i\|_2^2$, $\numh=30$;
    \item $g(\optvar)=\alpha\|\nabla \optvar\|_1+\iota_{[0,\infty)^d}(\optvar)$, $\optop_i=\ipop_i$ with $f_i(\optvar)=\frac{1}{2}\|\optvar-\datavar_i\|^2$, $\numf=30$ and $h\equiv0$.
\end{itemize}
As in Section~\ref{sec:ip}, matching the inverse problem to the template is not unique and an algorithm design question. In Section \ref{sec:numerics} we will revisit this example numerically.
\end{myex}

Similar to Section \ref{subsec:building_blocks}, we first discuss basic approaches for addressing these challenges and ideas to design stochastic algorithms. 

\subsubsection{Basic Ingredients}

\paragraph{Stochastic Gradient Approximation} 
The underlying idea of stochastic gradient methods is to approximate the full gradient of the smooth term by a stochastic estimator.
The prototypical example is the stochastic gradient descent (SGD)~\cite{robbins1951stochastic} which in each step approximates the full gradient $\nabla h(\xk) = \sum_{i=1}^\numh \nabla h_i(\xk)$ by a gradient of a single function $\tilde \nabla h(\xk) \approx \numh \nabla h_{i_k}(\xk)$, using an random index $i_k$, \new{selected uniformly with replacement}\footnote{\new{Other forms of sampling of indices have been studied, such as uniform sampling without replacement and different forms of importance sampling. However, in practice the performance is often indistinguishable, but the price to pay is a more complicated convergence analysis}}
As a result, the per-iteration computational effort is reduced by a factor of $1/\numh$, and thus one iteration of GD corresponds to executing $\numh$ iterations of SGD. This is referred to as one \emph{epoch} or \emph{data pass}. The drawback of SGD is that the error introduced by the variance of the random process is not reduced along the iterations. 
A zoo of algorithms has been developed to address this issue over the years,
which together with SGD, will be discussed in more detail in subsection \ref{subsec:SGM}.

\paragraph{Stochastic Proximal Point Algorithms}
The stochastic proximal point algorithm is the stochastic counterpart of the PPA for \replaced{non-smooth}{nonsmooth} problems, by evaluating a single randomly chosen proximal operator only, i.e. $\xkp = \operatorname{prox}_{\tau_k g_{i_k}}(\xk)$ \cite{bianchi2016ergodic}. 
Similarly to the smooth setting, variance reduced variants can be derived \cite{Traore2023}.

\paragraph{Randomised Dual Coordinate Descent} 
As is the case for deterministic algorithms, stochastic algorithms can also be applied to the dual formulation of the problem.
Namely, the dual problem to $\min_x \sum_{i=1}^\numf f_i(\optop_i \optvar)$ is computed as $\min_y \sum_{i=1}^\numf f_i^*(y_i)$ subject to $\sum_{i=1}^\numf \optop_i y_i = 0$ \cite{shalev2013stochastic}. Notice that the finite-sum structure in the primal problem leads to a dual problem which is separable in $y$ and thus naturally amendable to dual coordinate descent strategies \cite{richtarik2014iteration}. This idea has been exploited in primal-dual algorithms for \eqref{eq:opttemplate_finite_sum} like QUARTZ \cite{qu2015quartz} and SPDHG \cite{Chambolle2018spdhg}. 
\new{
\begin{remark}
    It is worth noting that stochastic gradient descent and stochastic coordinate descent are deeply connected via duality, and we refer to \cite{shalev2013stochastic,gower2015randomized} for discussions.
\end{remark}
}

Table \ref{tab:elements-sto} summarises the basic ingredients for the design of first-order stochastic optimisation algorithms.  
Next we provide detailed explanations of these methods and discuss specific state-of-the-art algorithms that use these concepts. 

\begin{table}[htbp]\renewcommand\arraystretch{1.35}
\centering
\begin{tabular}{c@{\hspace{4mm}}c@{\hspace{4mm}}c}
\toprule
Target & Tool & References \\ 
\midrule \midrule
$\sum_{i=1}^\numh h_i(x)$ & Stochastic gradient approximation & \cite{robbins1951stochastic, bottou2018optimisation} \\
\hline
$\sum_{i=1}^\numh h_i(x)$ & Better stochastic gradients with history & \cite{Defazio2014, johnson2013svrg} \\
\hline
$\sum_{i=1}^\numg g_i(x)$ & Stochastic proximal operators & \cite{bianchi2016ergodic, Traore2023} \\
\hline
$\sum_{i=1}^\numf f_i(A_i x)$ & Coordinate descent for dual problem & \cite{shalev2013stochastic, Chambolle2018spdhg} \\ \bottomrule
\end{tabular}
\caption{Elements in designing stochastic first-order algorithms.
} 
\label{tab:elements-sto}
\end{table}

\subsubsection{Designing Stochastic Algorithms}
Once the computational expensive components of a deterministic optimisation algorithm have been identified, they can be replaced with a stochastic alternative. 
However, there are still several delicate tasks, for example, decomposing \eqref{eq:opttemplate} into the finite-sum form~\eqref{eq:opttemplate_finite_sum} is not unique,  choosing $\numf, \numg, \numh$ and partitioning $\optop$ such that the performance of the algorithms are optimised is not a straightforward task in general. 

In the rest of the section we introduce several representative stochastic optimisation algorithms and discuss related topics.

\subsection{Stochastic Gradient Methods} \label{subsec:SGM}

In this section we review in detail algorithms based on stochastic approximation of the gradient, which dates back to 1950s. 
In their seminal work \cite{robbins1951stochastic}, Robbins and Monro \new{proposed a stochastic approximation procedure for finding the root of a function that can be expressed as an expected value. 
Another important insight is that the data in many large-scale real-world applications show a substantial amount of true or approximate redundancy \cite{bottou2018optimisation}. 
Thus, using all available data in each computation is inefficient. 
Moreover, }GD is quite robust in the sense that convergence towards the minimum can be ensured even when there are small errors in the descent direction, as long as the errors decrease in the long run.
This has significant practical benefits, allowing to efficiently solve large scale optimisation problems.

Since the 2010s, stochastic gradient schemes have seen a resurgence of interest, due to the machine learning renaissance, and the library of existing methods has been significantly enriched. 
The new algorithmic modifications include adaptive steps-sizes, variance reduction schemes, and many more.

\subsubsection{Vanilla Stochastic Gradient Descent}\label{subssec:SGD}

\replaced{When the composite terms vanish and we consider a single prox-friendly function, i.e.\ $\numf=0$, $\numg=1$, template \eqref{eq:opttemplate_finite_sum} reduces to}{When we consider a single prox-friendly function and without the composite terms, i.e.\ $\numf=0$, $\numg=1$, template \eqref{eq:opttemplate_finite_sum} reduces to}
\begin{equation}\label{eq:two_finite_sum}
\min_{\optvar \in \optspace} \left\{ \Phi(\optvar) = g(\optvar) + \sum_{i=1}^\numh h_i(\optvar) \right\}.
\end{equation}
\removed{Stochastic gradient descent (}SGD\removed{)} \cite{robbins1951stochastic} is a classical stochastic optimisation method for efficiently solving this type of a problem, see Algorithm \ref{alg:sgd}.
In each iteration $k$ it performs a step in the (negative) direction of $\nabla h_{i_k}(\optvar)$ for a randomly sampled index $i_k\in\{1,\ldots,\numh\}$ using a pre-defined step{\removed{-}}size schedule $\tau_k > 0$ (a.k.a.\ learning rate).

\begin{center}
\begin{minipage}{0.95\linewidth} 
    \begin{algorithm}[H] 
        \caption{Stochastic Gradient Descent (SGD)}\label{alg:sgd}
        \DontPrintSemicolon
        \KwIn{$x^{(0)} \in \optspace$, $\tau_k > 0$.}
        \For{$k=0,1,\ldots$}{
        Sample $i_k \in\{1,\ldots,\numh\}$ uniformly at random \\
        $\tilde\nabla h(\xk) =$ \mybox[yellow]{$\numh\nabla h_{i_k}(\xk)$} \hfill{\tiny\texttt{stochastic gradient}}\\
        $\xkp =$ \mybox[green]{$\mathrm{prox}_{\tau_k g}(\text{\myinbox[red]{$\xk - \tau_k\tilde\nabla h(\xk)$}})$}  \hfill{\tiny\texttt{proximal descent on nonsmooth term}}
        }
    \end{algorithm}
\end{minipage}
\end{center}

Instead of using a single function $h_i$ in each iteration, in machine learning it is common to use several randomly selected functions, resulting in the so-called mini-batch SGD.
In contrast, for inverse imaging problems it is more common to pre-batch the data by partitioning the forward operator $\ipop$, and define functions $h_i$ from the partition. In this context, batches are usually referred to as \emph{subsets}, and thus the number of batches and batch size become \emph{number of subsets} and \emph{subset size}. We will discussed subset selection in more detail in Section~\ref{sec:challenges}. SGD has been used in image reconstruction, e.g.\ for PET~\cite{Twyman2023rdp}, CT \cite{Gao2018}, optical tomography \cite{chen2018opticaltomography} and cryo-EM \cite{singer20220singleparticlecryoem,PunjaniNatureMethods2017cryoem}.

The descent direction in SGD is an unbiased estimator of the full gradient, that is, $\E[ \numh\nabla h_i(\optvar)]=\nabla h(\optvar)$. 
\new{Such a stochastic approximation allows us to rewrite stochastic gradient descent into the following form
\begin{equation}\label{eq:stochastic_perturbation}
    \xkp = \mathrm{prox}_{\tau_k g} \big( \xk - \tau_k ( \nabla h(\xk) + \ek ) \big)
\end{equation}
where $\ek = \tau_k\tilde\nabla h(\xk) - \nabla h(\xk)$ is the approximation/perturbation error. 
From this perturbation perspective, to prove convergence, we need to control the approximation error. 
If the error $\ek$ is deterministic, according to \cite{schmidt2011convergence} which studies the convergence property of PGD with error in gradient, \eqref{eq:stochastic_perturbation} is convergent if $$ \sum_{i=0}^{k} \tau_i\norm{e^{(i)}} = o(\sqrt{k}) .$$ This means that: 1) If $\norm{\ek}$ is bounded away from $0$, the stepsize $\tau_k$ should vanish to $0$ fast enough; 2) If $\norm{\ek}$ vanishes, then constant stepsize can be allowed.
The first scenario is exactly the case of SGD, that the error $\ek$, though stochastic, is not vanishing. This is why in many convergence analysis of SGD, one either needs to impose the assumption that the variance is bounded,  
$\E [\|\ek\|^2] = \E [\|\numh\nabla h_i(\optvar) - \nabla h(\optvar)\|^2] \leq C ,$ for some $C>0$, which essentially means that the error $\ek$ is bounded;
Or impose the so-called ``expected smoothness'' assumptions as studied in \cite{khaled2020better} for weaker conditions than bounded variance. 
}

\new{Nevertheless, either assumption above only shows $\ek$ is controllable, one needs to use vanishing stepsizes to achieve convergence, and SGD does not converge for constant step{\removed{-}}sizes.}
As an example, typical vanishing step{\removed{-}}sizes takes the following form according to \cite{robbins1951stochastic}
\begin{align}\label{eqn:sgd_stepsize_condition}
\sum_{k=1}^{\infty} \tau_k = \infty, \quad \sum_{k=1}^{\infty} \tau_k^2 < \infty .
\end{align}
For example $\tau_k = k^{-\delta}, \delta \in (1/2, 1]$ satisfies the above conditions, though other choices are also popular, e.g.\ \cite{bottou2018optimisation}.
However, decaying step{\removed{-}}size regimes eventually lead to a very slow speed of convergence. 
To circumvent the slowdown in convergence, there are in general two types of modifications: variance reduction schemes and adaptive step{\removed{-}}sizes.

\begin{remark}
    In Algorithm \ref{alg:sgd} indices $i_k$ are sampled at random, which is typically done according to a uniform distribution in the discrete set $\{1,\ldots,n\}$, without replacement.
    Alternatively, a probability $p_i$ is assigned to each $i\in\{1,\ldots,n\}$.
    Choice of the sampling probability distribution can reduce the variance in SGD.
    However, for certain problems it can be shown that the optimal probability distribution is adaptive and requires solving a minimisation problem (in each iteration) that can be harder than the original problem.
    For $L$-smooth functions setting $p_i=L_i/\sum_{j=1}^n L_j$ is a good proxy \cite{zhao2015importance}. 
\end{remark}

\new{
\begin{remark}\label{rmk:sde}
Following Remark \ref{rmk:ode}, SGD can also be connected with continuous dynamical systems. In particular, SGD can be seen as the discretisation of the stochastic differential equation  $ x'(t) = -\nabla f(x(t)) + \sqrt{\eta} \Sigma(x(t)) {\rm d} W_t$, where $W(t)$ is a Wiener process (standard Brownian motion) capturing the stochastic noise and $\Sigma(x(t))$ is a noise scaling matrix that depends on the variance of the stochastic gradients. We refer to \cite{li2019stochastic,latz2021analysis} and the references therein for more discussions.  
\end{remark}
}

\subsubsection{Variance Reduction Schemes}\label{subsec:variance_reduction}
Stochastic variance reduction algorithms aim to reduce the variance present in SGD through historical gradient information.
In the traditional statistics literature such quantities are sometimes referred to as control variates. 
These approaches combine the low-per iteration complexity of SGD with constant step{\removed{-}}sizes (and faster convergence rates) of PGD, at a cost of increased storage requirements.
The study of variance reduction algorithms began with Stochastic Average Gradient (SAG)~\cite{2012_Roux_Stochastic}, Stochastic Averaged Gradient Am\'elior\'e\footnote{The term am\'elior\'e is used in \cite{gower2020variance}. It is French and translates to English as \enquote{improved}.} (SAGA) \cite{Defazio2014}, and Stochastic Variance-Reduced Gradient (SVRG) \cite{johnson2013svrg}, but has since grown significantly.
We refer to \cite{driggs2022biased, gower2020variance} for an overview. 
Variance reduction algorithms have been used for many image reconstruction problems, such as PET \cite{Twyman2023rdp, kereta2021svrem}, CT \cite{davood2017vrct}, optical tomography \cite{macdonald2020efficient}, seismic tomography \cite{hu2023seismic} and blind deblurring~\cite{driggs2021stochastic}. 
Below we describe SAGA and SVRG, two representative variance-reduced variants of SGD. 

\paragraph{SAGA and SAG} 
In order the reduce the variance, SAG and SAGA imitate PGD by storing a table of the latest computed gradients for each sub-function $h_i$. 
The entire table is used in each iteration as the approximation of the full gradient, and thus information from all the gradients $\nabla h_i$, which lowers the variance.
However, the table is updated only at the position of the most recently used (i.e.\ randomly sampled) function index, which reduces the computational effort.

\begin{center}
\begin{minipage}{0.95\linewidth} 
\begin{algorithm}[H] 
    \caption{Stochastic Average Gradient Am\'elior\'e (SAGA)}\label{alg:saga}
        \DontPrintSemicolon
        \KwIn{$x^{(0)} \in \optspace$, $\tau>0$}
        \For{$k=0,1,\ldots$}{ 
            Sample $i_k \in\{1,\ldots,\numh\}$ uniformly at random\\
            $\tilde\nabla h(\xk)=$ \mybox[yellow]{$\numh \big(\nabla h_{i_k}(\xk) - \gamma_{i_k}^{(k)}\big) + \sum_{i=1}^\numh \gamma_i^{(k)}$} \hfill{\tiny\texttt{variance reduction}}\\
        $\xkp =$ \mybox[green]{$\mathrm{prox}_{\tau_k g}(\text{\myinbox[red]{$\xk - \tau_k\tilde\nabla h(\xk)$}})$}  \hfill{\tiny\texttt{proximal descent on nonsmooth term}}
            $\gamma_{j}^{(k+1)} = \begin{cases} \nabla h_{j}(\xk) & \textrm{if~~} j = i_k , \\ \gamma_{j}^{(k)} & \textrm{otherwise} \end{cases}$\hfill{\tiny\texttt{updating the stored gradients}}
        }
\end{algorithm}
\end{minipage}
\end{center}
When $h_i$ are convex, the average iterate of SAGA converges to a minimiser of the objective at a sublinear rate $\mathcal{O}(1/k)$ provided
$\tau = 1/(3 \numh \Lmax)$
\cite{defazio_phd}. 
Note that the convergence of SAGA is given in terms
of the average iterate, and not the final iterate, as common for deterministic algorithms in Section \ref{sec:optalgs}. 
This is a fairly common issue for stochastic optimisation algorithms, and arises due to proof techniques. However, in practice it is common to use the final iterate instead of the average one, regardless of the missing theoretical guarantees.

SAG follows nearly identical iterations, except that there is no $\numh$ factor in front of $\nabla h_{i_k} (\xk) - \gamma_{i_k}^{(k)}$, which lowers the variance and makes SAG {\it biased}.

The original convergence analysis of SAG required a computer aided approach to complete the proof, though the convergence was recently proven without the need for computer aid \cite{driggs2022biased}.

Storing the table of historical gradients can become infeasible when the number of subsets is large, or when data is of a high dimensionality. 
However, when applied to linear models (in the sense of \cite{gower2020variance}), the memory requirement can be significantly improved in certain applications, see Section \ref{sec:challenges} for details.

\new{
\begin{remark}
    Momentum based acceleration methods can be seen as a form of variance reduction, as they use a weighted average of previous gradients which gives more weight to recent gradients. To see this, consider the following accelerated SGD using heavy-ball momentum,
\[
\begin{aligned}
\yk &= \xk + a(\xk - \xkm) \\
\xkp &= \yk - \tau v^{(k)}   \\
\end{aligned}
\]
where $v^{(k)}$ represents the stochastic gradient at step $k$. We consider constant inertial parameter $a$ and stepsize $\tau$ for simplicity. Under this setting, the update of $\xkp$ can be written into a recursion which reads
\[
\begin{aligned}
\xkp 
&= \xk - \tau \sum_{i=0}^{k} a^i v^{(k-i)}  
= \xk - \tau \Big( v^{(k)} + \sum_{i=1}^{k} a^i v^{(k-i)} \Big)  .
\end{aligned}
\]
The term $\sum_{i=1}^{k} a^i v^{(k-i)}$ can be seen as variance reduction term, with exponentially decaying weight. 
While in contrast, SAG(A) uses a plain average of the recent gradients. 
\end{remark}
}

\paragraph{SVRG}   
Instead of storing a table of gradients, SVRG stores a reference iterate, often called an anchor point, or a snapshot, as well as the full gradient evaluated at the reference iterate. 
Both the reference iterate and the full gradient are updated at pre-defined intervals during the iterations.
To this end, SVRG requires two loops: in the outer loop we update the reference iterate and compute the full gradient, and in the inner loop we perform stochastic updates using a combination of the current gradient and the full gradient from the outer loop. 
The length  $t\in\mathbb{N}$ of the inner loop is by convention set to $2\numh$ for convex and $5\numh$ for nonconvex problems.

\begin{center}
\begin{minipage}{0.95\linewidth} 
\begin{algorithm}[H] 
    \caption{Stochastic Variance Reduced Gradient (SVRG)}\label{alg:svrg}
        \DontPrintSemicolon
        \KwIn{$\optvar^{(0)} \in \optspace$, $\tau>0$,
        inner loop length $t\in\mathbb{N}$}
        Set $\overline{\optvar}^{(0)}=\optvar^{(0)}$\\
        \For{$j=0,1,\ldots$}{  
            Compute $\overline{\gamma}^{(j)} = \sum\nolimits_{i=1}^{\ell} \nabla h_{i}(\overline\optvar^{(j)})$\hfill{\tiny\texttt{full gradient at reference iterate}}\\
            \For{$k=0,\ldots,t-1$}{
                 Sample $i_k \in\{1,\ldots,\numh\}$ uniformly at random \\
                 $\tilde \nabla h(\optvar^{(k)})=$ \mybox[yellow]{$\numh \big(\nabla h_{i_{k}}(\optvar^{(k)}) - \nabla h_{i_{k}}(\overline{\optvar}^{(j)})\big) + \overline{\gamma}^{(j)}$} \hfill{\tiny\texttt{variance reduction}}\\
        $\xkp =$ \mybox[green]{$\mathrm{prox}_{\tau_k g}(\text{\myinbox[red]{$\xk - \tau_k\tilde\nabla h(\xk)$}})$}  \hfill{\tiny\texttt{proximal descent on nonsmooth term}}
            }
            Set $\overline{\optvar}^{(j+1)} = \optvar^{(t)}$\hfill{\tiny\texttt{store reference iterate}}
        }
\end{algorithm}
\end{minipage}
\end{center}
SVRG was initially designed only for smooth objectives \cite{johnson2013svrg}.
The proximal version for composite objectives of the type \eqref{eq:two_finite_sum} was developed and analysed in \cite{xiao2014}.
They show for strongly convex problems that the average iterate of SVRG converges in the objective value at a linear rate for $\tau < 1/(4 \numh \Lmax)$
and sufficiently large $t$~\cite{xiao2014}.

The advantage of SVRG over SAG(A) is that it requires storing only two pieces of data, i.e. $\overline{\optvar}^{(j)}$ and $\overline{\gamma}^{(j)}$, hence its memory requirements are independent of $\numh$.
However, its disadvantages are the need for two loops, which introduces an additional tuning parameter, and the need to evaluate the full gradient every time the reference iterate is updated (which is not needed by SAGA).

\paragraph{Adaptive Subset Sizes}
The number $\numh$ of smooth terms, or equivalently the batch size, has a direct impact on the variance of SGD and can be derived for general cases of unbiased stochastic estimators with bounded variance \cite{bottou2018optimisation}.
Specifically, the variance of the stochastic estimator increases as the number of subsets increases, at a roughly linear rate.
Therefore, variance can be reduced not by reducing the step{\removed{-}}size but rather by adaptively or dynamically selecting the batch or subset size. 
This topic has been predominantly studied from the deep learning perspective, e.g.~\cite{samuel2018dontbatch,devarakonda2017adabatch, balles2017adabatch}.
However, it is important to note that increasing $\numh$ also increases the per-iteration computational cost and can thus run into same scalability issues that PGD encounters.
Thus, timely increases of the batch size are critical in practice, and often specific criteria are designed to address this issue, e.g.~\cite{byrd2012sample}.
While different variations of this methodology have been used in certain inverse problems, such as diffuse optical tomography \cite{macdonald2020efficient} and PET \cite{twyman2019iterative}, the topic has not yet seen widespread adoption.

\subsubsection{Modes of Convergence and Convergence rates}

\new{The study of modes of convergence and convergence rates is a topic of critical importance for stochastic iterative methods, and has been investigated under a wide range of assumptions, though this has been predominantly done from the machine learning community and for related optimisation tasks.
Due to the stochastic nature of the iterations, the convergence of stochastic iterative method is studied from a probabilistic perspective.
That is, iterates (or their functions) are considered as random variables and one of several modes of convergence can be considered: 
convergence in mean (i.e. in expectation), almost sure convergence (with probability $1$), or convergence in probability (given as high probability guarantees). For stochastic methods, these are studied either for the last iterate, or for the average iterate.

Convergence rates, provided they exist, are typically as a function (power) of the iterate number, which heavily depends on the properties of the objective. 
Historically, convergence rates for stochastic iterative methods have mostly been developed for the last iterate of strongly convex objectives, i.e. when $h_i$ and $g$ are both strongly convex and Lipschitz smooth.
In case of non-strongly convex objectives, the convergence analysis is often conducted for the average iterate instead.
However, this sets an unfortunate predicament since for inverse problems the most informative mode of convergence of iterative methods is the convergence of the last (or running) iterate for non-strongly convex functions $h_i$ and convex, proper, lower semi-continuous $g$.
Such convergence rates have unfortunately not been established for all methods. 

An important consideration for the convergence of stochastic gradient methods is the effect of gradient variance on the step-size regime. Specifically, in order to control the gradient noise SGD requires decaying stepsizes, which impede its convergence speed. This issue is resolved with variance reduction methods, which are able to decrease the gradient noise fast enough, while maintaining constant stepsizes to provide a faster rate of convergence.
In Table \ref{tab:stoch_convergence_rates} we summarise the current state of the art regarding convergence rates and mode of convergence for Lipschitz smooth $h_i$ and proper, convex, lower semi-continuous $g$, with respect to average and last iterate, and convexity properties.

\begin{table}[htbp]\renewcommand\arraystretch{1.35}
\centering
\begin{tabular}{c@{\hspace{4mm}}c@{\hspace{4mm}}c@{\hspace{4mm}}cc@{\hspace{4mm}}c}
\toprule
& \multicolumn{2}{c}{convex}  & & strongly convex \\ \cline{2-3} \cline{5-5}
Method & Average Iterate & Last Iterate  & &  Last Iterate & References\\ 
\midrule \midrule
SGD & $\mathcal{O}(1/\sqrt{k})$  & $\mathcal{O}(1/\sqrt{k})$ & & $\mathcal{O}(1/k)$ & \cite{garrigos2023handbook}\\
\hline
SVRG &  $\mathcal{O}(1/k)$ & almost sure  & & $\mathcal{O}(\rho^k)$ & \cite{poon2018local, xiao2014, zhang2022linear}\\
\hline
 SAGA & $\mathcal{O}(1/k)$ &  almost sure   & & $\mathcal{O}(\rho^k)$ & \cite{poon2018local, Defazio2014}\\
\bottomrule
\end{tabular}
\caption{
\new{Convergence rates of stochastic first-order algorithms. In all instances we consider proper, lower semi-continuous $g$ and $L$-smooth $h_i$. Term $k$ refers to the current iterate number and $0<\rho<1$ denotes an algorithm- and function-dependent constant.}
} 
\label{tab:stoch_convergence_rates}
\end{table}
}

\subsubsection{Further Stochastic Gradient-Based Methods}
SGD and variance reduction algorithms have in recent years seen many modifications and adaptations to specific scenarios of interest. 
For example, Loopless-SVRG algorithm~\cite{kovalev2020lsvrg} is a modification of SVRG that uses only one loop and instead of updating the reference iterate and full gradient at pre-defined intervals, in each iteration it takes a random sample from $[0, 1]$. If the sampled value is below a predefined probability, reference and full gradient are updated, and otherwise a standard update is perform.
SARAH \cite{nguyen2017sarah} is a biased variant of SVRG that updates the full gradient not only in the outer loop but also in the inner loop. Similar biased variance-reduction gradient schemes include SPIDER \cite{fang2018,wang2019}. 
See also \cite{driggs2022biased} for a gradient estimator named SARGE which combines the ideas of SAGA and SVRG. 

Algorithms for the three-term sum from Section \ref{sec:threetermsum} also admit stochastic variants. 
For example, stochastic variants of PD3O or Condat--V\~u algorithms can be combined with SGD or SVRG in a mix-and-match fashion by replacing the gradient $\nabla h(\optvar)$, with a corresponding stochastic estimator, see \cite{dualize2022salim}.
Similarly, for other first-order methods that contain gradient component, one can replace the gradient with a stochastic estimator, such as the primal-dual fixed-pint algorithm \cite{zhu2021stochastic}, the stochastic Davis--Yin three-operator splitting method \cite{yurtsever2016stochastic}.  
Notice that these stochastic methods use stochasticity to the primal variable, e.g. $\nabla h(\optvar)$. While for dual variables, we leave the discussion in Section \ref{subsec:scd} as it leads to different type of stochastic optimisation methods. 

Principles used for SGD, SAGA, and SVRG have also been applied to maximum likelihood expectation maximisation, and MAP, see \cite{cappe2009sem, chen2018vrsem, karimi2019sagaem}, and have been applied to tomographic image reconstruction \cite{kereta2021svrem}. 

\subsubsection{Accelerated Stochastic Variance-Reduced Gradient Methods}
As variance-reduced gradient methods can use constant step{\removed{-}}size and have the same convergence rate as their deterministic counterparts, a natural question to ask is whether NAG/FISTA can be extended to them. 
This is a rather challenging task as NAG/FISTA are very sensitive to errors \cite{aujol2015stability}.

The first attempt along this direction is \cite{nitanda2014stochastic}, in which FISTA is combined with SVRG gradient estimator, however no theoretical guarantee when the problem is not strongly convex. 
Later on in~\cite{lin2015universal}, based on an accelerated proximal point algorithm, a so-called Catalyst acceleration is proposed to speeding up the solving of~\eqref{eq:two_finite_sum}. Catalyst contains two loops, the outer loop is accelerated proximal point algorithm, while the inner loop is solved with a stochastic gradient scheme to certain accuracy. 
In~\cite{allen2017katyusha}, an acceleration scheme named Katyusha is proposed for accelerating SAGA and SVRG. By incorporating a ``negative momentum'', the author manages to prove an $\calO(1/k^2)$ convergence rate without strong convexity. In~\cite{driggs2022acceleration}, a more elegant acceleration scheme is proposed, see Algorithm~\ref{alg:acc-vr} below, and the authors show that the negative momentum is not needed to achieve the $\calO(1/k^2)$ convergence. 

\begin{center}
\begin{minipage}{0.95\linewidth} 
\begin{algorithm}[H] 
    \caption{Universal Scheme for Accelerating Variance-Reduced Methods}\label{alg:acc-vr}
        \DontPrintSemicolon
        \KwIn{$x^{(0)}, y^{(0)}, z^{(0)} \in \optspace$, $\tau>0$, $\eta_k>0$}
        \For{$k=0,1,\ldots$}{
            $x^{(k+1)} = \eta_k z^{(k)} + (1 - \eta_k)y^{(k)}$\\
            Compute a variance-reduced stochastic gradient estimate $\tilde\nabla h(\xkp)$ \\
            $\zkp=\prox_{\tau g}\big(\zk-\tau \tilde\nabla h(\xkp)\big)$\\
            $\ykp = \eta_k \zkp + (1 - \eta_k)\yk$
        }
\end{algorithm}
\end{minipage}
\end{center}

In recent years, vector extrapolation schemes are also adopted in the stochastic optimisation, for example the nonlinear acceleration \cite{scieur2017nonlinear} schemes (e.g. the Anderson Acceleration and its variants) are applied to accelerate the stochastic gradient methods. This type of schemes are very closely related to quasi-Newton schemes such as L-BFGS, as showed in \cite{scieur2019generalized}.

\subsubsection{Incremental Gradient Methods} 
Before wave of stochastic optimisation algorithms for machine learning in the 2010s, another related family of optimisation algorithms received increased interest: incremental gradient methods. Incremental gradient methods are also designed for solving objective functions with the finite-sum structure \eqref{eq:two_finite_sum} as for SGD (see Algorithm~\ref{alg:sgd}). The core idea is that in each step of the iteration, the algorithm chooses one function from the finite-sum via a predetermined rule (cyclic or fixed permutation), and apply GD to the chosen function \cite{Bertsekas2011}. 
Similar to sub{\removed{-}}gradient descent and SGD, to ensure convergence, a carefully selected step{\removed{-}}size is required: $\sum_{k} \tau_k = \infty$, $\sum_{k} \tau_k^2 < \infty$. This idea is used in the modified block sequential regularised expectation-maximisation (BSREM) and the relaxed OS separable paraboloidal surrogates (OS-SPS)~\cite{DePierro2001, Ahn2003}, the former is used currently in clinical PET imaging \cite{Sah2017}.

To circumvent the need for a diminishing step{\removed{-}}size, an incremental averaged gradient (IAG) scheme was developed in \cite{Blatt2007} which can afford constant step{\removed{-}}size. The basic idea of the method is averaging the history gradients, which is also the key behind stochastic averaged gradient methods  SAG and SAGA (see e.g.\ Algorithm~\ref{alg:saga})~\cite{schmidt2017sag,Defazio2014}. 

This methodology can be extended to \replaced{non-smooth}{nonsmooth} optimisation (e.g.\ for PGD).
For instance when the \replaced{non-smooth}{nonsmooth} function is decomposable into a finite-sum, one can choose a single term in the sum and compute its proximal operator. We refer to the dedicated survey \cite{Bertsekas2011} for more details.

Lastly, compared to the stochastic counterpart, incremental gradient methods have weaker theoretical guarantees, such as more restrictive step{\removed{-}}size criteria \cite{Bertsekas2011, Blatt2007, schmidt2017sag}. The numerical investigation of \cite{Twyman2023rdp} in the context of PET imaging showed that BSREM was inferior to stochastic variance reduced algorithms in every considered metric. See Section \ref{sec:numerics} for a comparison of incremental gradient methods to stochastic methods.
We add that incremental gradient algorithms that access gradient indices in a cyclic manner can show exponentially bad behaviour for certain pathological use-cases \cite{herman1993hmorder, recht2012access, defazio_phd}.

\subsection{Stochastic (Dual) Coordinate Descent}\label{subsec:scd}
In the context of inverse problems, stochastic gradient methods often reduce the complexity in measurement space. For example, when the smooth term $h$ corresponds to data fidelity, e.g.\ $\frac 12\|\ipop\optvar-\datavar\|^2$.
An alternative is to reduce the complexity in the space of the optimisation variable. Recall the CD from Section \ref{sec:fom+basic_element}, which performs the iterations 
\[ \xkp=\xk-\tau_k \nabla_i h(\xk),\]
where $\nabla_i h(\optvar) = \frac{\partial h (\optvar)}{\partial \optvar_i} e_i$ and $e_i$ is the \replaced{$j$-th}{$i$-th} canonical basis vector.
As in the deterministic case, instead of using a single coordinate we can use a subset of (randomly selected) coordinates. Note that the direction $\nabla_i h(\optvar)$ is an unbiased estimation of $\nabla h(\optvar)$. 
Moreover, it is also variance-reduced. If $\xk\rightarrow\optvar^\ast$, then due to the optimality condition of $\optvar^\ast$, we get $\nabla h(\xk) \to 0$ and $\nabla_i h(\xk) \to 0$, thus $\|\nabla_i h(\xk)-\nabla h(\optvar^\ast)\|^2\rightarrow0$. 

To provide a reduction in computational costs it usually requires coordinate-friendly structures (such as separable functions), that is, coordinate updates should be significantly cheaper than gradient computations. Because the optimisation problems associated with inverse problems usually do not exhibit this structure, coordinate descent methods are largely not used in this context. An exception is for instance~\cite{Gao2018} which uses coordinate descent in large-scale CT image reconstruction to overcome memory bottlenecks on GPUs. It turns out that coordinate descent for inverse problems is more useful when applied on dual variables as we will see below.

Whereas stochasticity has been observed to improve the performance for gradient estimators, as compared to using a deterministic (e.g.\ cyclic) sampling of the gradients, the situation for coordinate descent is more nuanced.
On one hand, there are classes of \replaced{non-convex}{nonconvex} problems for which cyclic CD fails to converge to a stationary point~\cite{powell1971convergence}, but stochastic CS (SCD) is expected to converge within only a few steps, and in other cases SCD can be significantly faster than cyclic CD \cite{wright2015coordinate}.
On the other hand, one can design classes of forward operators such that cyclic CD for the corresponding quadratic problem converges several times faster than SCD \cite{gurbuzbalaban2017scdvsccd}. 
For a further comparison between cyclic, stochastic, and adaptive choices for block CD we refer to~\cite{nutini2022bcd}.

\subsubsection{Stochastic Primal-Dual Splitting}
For simplicity we will restrict ourselves to the case when $\numg = 1$ and $\numh = 0$, i.e.\ our template \eqref{eq:opttemplate_finite_sum} takes the form
\begin{equation}\label{eq:opttemplate_finite_sum_noh}
\min_{\optvar \in \optspace} \left\{ \Phi(\optvar) = \sum_{i=1}^\numf f(\optop_i \optvar) + g(\optvar)\right\} , 
\end{equation} 
with associated saddle-point problem
\begin{equation}\label{eq:opttemplate_finite_sum_noh:saddle}
\min_{\optvar \in \optspace} \max_{\dualvar \in \dualspace} \left\{ \sum_{i=1}^\numf \left( \langle \optop_i \optvar, \dualvar_i\rangle - f_i^*(\dualvar_i)\right) + g(\optvar)\right\}.
\end{equation} 

When applying PDHG (Algorithm \ref{alg:pdhg}) to this problem every iteration requires evaluating matrix-vector products for all $\optop_i^*$ to update the primal variable and $\optop_i^\ast$ to update the dual variable which can be computationally costly. Thus, in \cite{Chambolle2018spdhg} it was proposed to perform a dual coordinate ascent instead. The resulting algorithm is called 
Stochastic Primal-Dual Hybrid Gradient (SPDHG) and provided in in Algorithm~\ref{alg:spdhg}.

\begin{center}
\begin{minipage}{0.95\linewidth} 
\begin{algorithm}[H] 
    \caption{Stochastic Primal-Dual Hybrid Gradient (SPDHG)}\label{alg:spdhg}
        \DontPrintSemicolon
        \KwIn{$x^{(0)} \in \optspace, y^{(0)} = \overline{\dualvar}^{(0)} \in \mathcal Y$, $\tau, \sigma > 0$}
        \For{$k=0,1,\ldots$}{ 
        $\xkp =$ \mybox[green]{$\prox_{\tau g}(\text{\myinbox[red]{$\xk - \tau \optop^* \overline{\dualvar}^{(k)}$}})$}\hfill{\tiny\texttt{proximal descent on nonsmooth term}}\\
        Sample $i_k \in\{1,\ldots,\numf\}$ uniformly at random \\
        $\dualvar^{(k+1)}_i = \begin{cases} \mybox[blue]{$\prox_{\sigma f_i^*}\big(\text{\myinbox[red]{$\dualvar^{(k)}_i + \sigma \optop_i \xk$}}\big)$} & \text{if $i = i_k$} \\
        \dualvar^{(k)}_i & \text{else}
        \end{cases}$ 
        \hfill{\tiny\texttt{proximal ascent and stochastic dualisation}}\\
        $\overline{\dualvar}^{(k+1)}_i = \begin{cases} \mybox[orange]{$\dualvar^{(k+1)}_i + \numf  (\dualvar^{(k+1)}_i - \dualvar^{(k)}_i)$} & \text{if $i = i_k$} \\
        \dualvar^{(k)}_i & \text{else}
        \end{cases}$\hfill{\tiny\texttt{extrapolation}}
        }
\end{algorithm}
\end{minipage}
\end{center}
SPDHG has been shown to converge in Bregman distance if $\tau\sigma < 1/(\numf \max_i \|\optop_i\|^2)$. This convergence result can be extended and the algorithm be generalised to other samplings and step{\removed{-}}sizes, see \cite{Chambolle2018spdhg, alacaoglu2022convergence, Gutierrez2021ssvm} for more details.

Note that due to the linearity of $\optop$, we can rewrite the primal update as
\begin{align*}
    \optop^* \overline{\dualvar}^{(k+1)} 
    = (1 + \numf) \optop^*_{i_k} (\dualvar^{(k+1)}_{i_k} - \dualvar^{(k)}_{i_k}) + \sum_{i=1}^\numf \optop^*_i \dualvar^{(k)}_i.
\end{align*}
Note that this bears a close resemblance to a SAGA estimator of the gradient of the bilinear term in \eqref{eq:opttemplate_finite_sum_noh:saddle} with respect to $\optvar$. However, due to the factor $1 + \numf$ instead of $\numf$, this is a biased gradient estimator.
Its effectiveness on large-scale problems has been shown in various applications such as PET with sinogram data~\cite{Ehrhardt2019pmb}, PET with listmode data~\cite{ListmodePET}, {PAT \cite{SPDHGPATthesis2020},} CT \cite{Chambolle2024, riis2021computed, CIL_II}, multi-spectral CT \cite{CIL_II}, parallel MRI~\cite{Gutierrez2021ssvm, Gutierrez2022} and magnetic particle imaging \cite{Zdun2021}.

As for any primal-dual algorithm, it is difficult to tune the ratio of the primal and dual step{\removed{-}}sizes. This is addressed empirically in \cite{Zdun2021} and with theoretical foundations provided in \cite{Chambolle2024}.
Note that the problem can also be approached similarly by other primal-dual algorithms, e.g.\ \cite{shalev2013stochastic, fercoq2015accelerated, Fercoq2019}.

\subsection{Further Reading}\label{subsec:further_stochastic}

\paragraph{Stochastic Second-Order Methods}
Standard second-order methods have seen limited application in large scale problems, due to the costs and intractability in computing the Hessian of the objective.  
As discussed in Section \ref{sec:determinstic_further}, deterministic quasi-Newton methods use information about the curvature without computing the Hessians. 
Instead, approximation techniques are used to approximate the Hessian or the corresponding matrix-vector products.
A direct extension of quasi-Newton methods to the stochastic setting can be challenging, since noisy estimates of the curvature (coupled with stochastic estimators of the gradient) can destabilise the algorithm. 
A number of methods have been proposed to solve this issue, and leverage the strengths of quasi-Newton and stochastic gradient methods. 

A pioneering work in this field is the stochastic quasi-Newton method \cite{schraudolph2007stochastic}, which introduced stochastic modifications of BFGS and L-BFGS, showing a sub-linear convergence rate.
This limitation was addressed by \cite{moritz2016linearly}, achieving a linear convergence rate (under mild assumptions) through a combination of the stochastic variant of L-BFGS \cite{byrd2016sqn} with variance reduction. 

Sub-sampled quasi-Newton methods, see e.g. \cite{Roosta2014, erdogdu2015ssm, peng2016ssm, roosta2019ssm}, randomly sub-sample the Hessian to further reduce the costs, and use matrix concentration inequalities to ensure that curvature information is well preserved.
Despite these advances, stochastic quasi-Newton methods are still not widely adopted neither in the inverse imaging community, with notable counterexamples~\cite{milzarek2019sqpat, hanninen2022sqgn, Perelli2021} and they continue to be an active research area.

\paragraph{Forward Operator Sketching}

The underlying idea of sketching methods is to replace the original inverse problem with an inverse problem that is much smaller and thus more manageable. 
This is often studied for the least-squares problem $\min_{\optvar \in \mathbb{R}^\recondim} \|\ipop\optvar - \datavar\|^2$. 
The idea is to find a sketching matrix $S \in \mathbb{R}^{m \times \datanum}$, where $m \ll \datanum$ such that the solutions of $$\min_{\optvar \in \mathbb{R}^\recondim} \|S(\ipop\optvar - \datavar)\|^2$$ are, in a certain sense, similar to the solutions of the original problem but are drastically cheaper to compute.

Common sketching techniques include random projections, sub-sampling, and leverage score sampling~\cite{woodruff2014sketching,drineas2006sampling,mahoney2011randomised}. However, this type of sketching comes at a cost since there is an irreversible loss of information \cite{pilanci2016iterative}. 

\paragraph{Iterative Hessian Sketching}
An alternative approach for solving inverse problems via sketching is through iterative Hessian sketching (IHS). 
Similarly to Newton's method, IHS constructs quadratic approximations to the original problem using Hessians, which are then sketched.
In more detail, consider the minimisation of a smooth convex objective $h(\optvar)$. 
At each iteration, a sketching matrix $S \in \mathbb{R}^{\datanum \times \recondim}$ is used to construct a sketched Hessian matrix $H_S = (H^{1/2})^* S^* S H^{1/2}$, where $H$ is the Hessian matrix at the current iterate $\xk$. IHS then follows the iterations
$$\xkp = \xk - \tau_k H_S^{-1} \nabla h(\xk),$$
where $\tau_k>0$ is the step{\removed{-}}size. 
For certain convex optimisation problems IHS is faster than baseline methods
~\cite{pilanci2016iterative,pilanci2017newton}. It can also be combined with (accelerated) PGD \cite{tang2017gradient}, and it was used to accelerate multi-coil MRI reconstruction \cite{oscanoa2024coil}. 

A particularly common choice for the sketching matrix are Gaussian sketches, which have the strongest theoretical guarantees (with respect to compression quality) but are computationally inefficient. 
They can be accelerated through the Johnson--Lindenstrauss transform~(FJLT)~\cite{ailon2009fast,2008_Ailon_Fast}, and the sparse Johnson--Lindenstrauss transform (SJLT)~\cite{clarkson2013low}, to name a few techniques. 
In particular, SJLT trades off computational efficiency for theoretical guarantees, but it is often used in practice.

Although they have some clear advantages, sketching-based methods are rarely applied to inverse imaging problems.
This can primarily be attributed to the fact that sketching may destroy any beneficial structure present in the forward operator. 
For example, in CT the Radon transform is a sparse matrix for which we can only use subsampling as the sketching, leading to suboptimal performance. 
The MRI operator includes FFTs which again prevents naive application of sketching. 
Recently in \cite{oscanoa2024coil}, a special sketching operator which can preserve the fast computational structure of FFT is developed and applied to multi-coil MRI reconstruction task.

\paragraph{Stochastic Mirror/Bregman Descent}
Stochastic gradient descent has in recent years been extended to the problems that do not rely on the Euclidean structure. 
Stochastic mirror descent is one powerful such generalisation which in each step can be written in the proximal point form as 
\[ \xkp = \argmin_{\optvar\in\reconspace} \left\{\tau_k \tilde\gamma^k + D_\phi(\optvar, \xk)\right\},\]
where $\phi$ is the mirror map and $\tilde \gamma^k$ is an estimator of the gradient of the objective at $\xk$.
SMD is often investigated through the lens of stochastic schemes for nonsmooth Lipschitz continuous convex functions \cite{bianchi2004smd}, though its convergence can be shown for nonsmooth convex objectives \cite{lan2012smd}, smooth objectives \cite{bubeck2015convex} and other conditions on the objectives, and it has been applied to ill-posed inverse problems \cite{jin2023mirror}.
Adaptive step{\removed{-}}size schemes have also been studied in recent years \cite{dorazio2023smd}.

SGD has also been designed to when both $\dataspace$ and $\reconspace$ are Banach spaces by applying a proximal point like algorithm \cite{jin2023banach}, with applications in CT.
In this case selecting $\dataspace$ and $\reconspace$ as certain $\mathcal{L}^p$, or Sobolev, spaces can act as an implicit regulariser of the solution and allows adapting to different noise assumptions.
However, choosing $\dataspace$ and $\reconspace$ to optimise the reconstruction is a challenging problem. 
This can to a degree be alleviated by using modular (or variable exponent) Banach spaces \cite{lazaretti2023variable}, which can adapt to cases of spatially variable noise types.

\paragraph{Adaptive step{\removed{-}}size Methods for Training Neural Networks} 
A significant advancement of SGD in the context of training (deep) neural networks are  adaptive step{\removed{-}}size stochastic gradient schemes, such as AdaGrad \cite{duchi2011adaptive}, AdaDelta \cite{zeiler2012adadelta}, RMSProp \cite{tieleman2012lecture}, \replaced{ADAM}{Adam} \cite{kingma2014adam}, to name a few. We refer to \cite{ruder2016overview} for a brief overview of these methods.

The main motivation behind this type of methods is the observation that some variables, or components of the solution, are updated significantly less during the training than other variables.
More precisely, let $d^{(k)} = \tilde{\nabla} h (\xk)$ denote the stochastic gradient computed in $k$-th iteration, for $\xk\in\R^\recondim$.
Then assume that for an iteration $T$ there exists an $1\leq i\leq\recondim$, such that the accumulation of update directions, defined as $D_i^{(T)}=\sum\nolimits_{k=1}^{T} | d_i^{(k)} |$ is small compared to $D_j^{(T)}$ for $i\neq j$. In this case the corresponding entry $\optvar_i$ of the solution will require a significantly longer training time. 
Adaptive steps-size methods overcome this issue by using step{\removed{-}}sizes tailored to each entry of $\optvar$. 
To this end, let $D^{(T)}$ be the accumulation of the stochastic gradient $d^{(k)}$.
For example, as above let
$
D^{(T)}_i = \epsilon + \sum_{k=1}^{T} | d_i^{(k)} |
$, for $i=1\ldots,n$, and a small positive constant $\epsilon > 0$. For a diagonal matrix $H_T = \mathrm{diag}(1/D_i^{(T)})$ the adaptive step{\removed{-}}size SGD can follow the iterations
\begin{equation}\label{eq:adaptive-sgd}
\xkp = \xk - \tau H_k d^{(k)} .
\end{equation}
Therefore, if for some $i$ the value $D^{(k)}_i$ is very small, then the corresponding entry of $H_k$ will be large, and hence the entry $\xkp_i$ will use a larger effective step{\removed{-}}size and receive a larger update. 

Compared to the standard SGD scheme, in general adaptive step{\removed{-}}size schemes have two advantages.
First, tuning of step{\removed{-}}size becomes automatic such that each entries converge at relatively same speed. Second, step{\removed{-}}size choices are tuning-free, for instance in \eqref{eq:adaptive-sgd} we only need to choose $\tau>0$. 
Over the past decade, adaptive step{\removed{-}}size SGD schemes have been more popular in training neural networks, while deeper theoretical understandings are obtained for the standard SGD.

\begin{remark}
We would like to highlight that adaptive step{\removed{-}}size stochastic gradient schemes are the dominant optimisation schemes in machine learning, especially the ADAM algorithm. In the field of inverse (imaging) problems, even though the adoption of stochastic algorithms is quite trendy nowadays, adaptive step{\removed{-}}size schemes have received rather limited attention. In comparison, variance reduction schemes are more popular. 
\end{remark}

\paragraph{Randomised Prox} 
Algorithms in Section \ref{sec:stoch} reduce the iteration complexity by using a low-cost estimator of the gradient.
Similarly, algorithms in Section \ref{subsec:scd} use a lower-cost update in the dual space.
This speeds up the computations in inverse problems since it reduces the computational costs inherent in applying the operator $\ipop$ or $\optop$.
However, the second source of computational costs is in the evaluation of the proximal operator, particularly for algorithms that compute $\prox_{\tau g}$ or $\prox_{\tau f}$ in every iteration. 
Since stochastic algorithms require more iterations to conduct the same number of data passes (and utilise all of the measurement data), this can incur further computational costs if $\prox_{\tau g}$ is expensive to evaluate.

To address this, new approaches have been developed in recent years that take the stochastic approach to the proximal domain.
In \cite{condat2023randprox} RandProx algorithm is introduced which uses a stochastic estimator for the proximal step in the primal-dual Davis--Yin algorithm. 
The stochastic estimator is in its simplest form defined by a probability parameter $p$, i.e.\ depending on a coin flip the proximal operator is either applied or not applied.
ProxSkip \cite{mishchenko2022proxskip} applies this methodology directly to proximal gradient-style algorithms. 
However, these methods are predominantly applied to federated or decentralised learning, and their use in inverse imaging problems has not yet been thoroughly investigated.

\paragraph{Stochastic Iterative Regularisation}
Regularising property of stochastic optimisers has been a well observed and studied property in the machine learning communities.
For linear inverse problems and a least squares type objectives, it has been shown that SGD \cite{jin2019sgdreg,jin21sgdsaturation} and SVRG \cite{jin2022svrgreg} for well-chosen a priori stopping criteria have a regularising property.
Similar results have been shown for mirror and Banach space variants of SGD \cite{jin2023mirror,jin2023banach}, for randomised Kaczmarz methods \cite{needell2014paved,tondji2023faster}, and other stochastic optimisation schemes.

\paragraph{Nonlinear Inverse Problems}
Stochastic gradient methods have been successfully applied to some suitably-structured nonlinear inverse problems. For example the randomised Kaczmarz method and its SVRG-inspired variance reduced variants have been successfully applied to phase-retrieval problems \cite{Xian2022}, EIT~\cite{van2012adaptive,zhou2023deep}, optical tomography \cite{macdonald2020efficient}, and other problems. However, in general, the application of stochastic gradient methods for \replaced{non-linear}{nonlinear} problems poses additional challenges, due to the complexities introduced by the forward model and its properties. \new{Moreover, nonlinear inverse problems usually require iteration dependent stepsizes which are commonly found via backtracking. Backtracking for random algorithms is an open research question and no general solution is available; see also Chapter \ref{sec:challenges}.
These problems often require a problem-specific solution and more research is needed in this area.}

\section{Challenges for Stochastic Optimisation for Inverse Problems} \label{sec:challenges} 
In the previous two sections, we provided an overview of both deterministic and stochastic optimisation algorithms. In this section, we present a thorough discussion on various theoretical and practical aspects of stochastic optimisation methods.  

\subsection{When and Why Are Stochastic Gradient Methods Fast?}
As we mentioned at the beginning of Section \ref{sec:stoch}, the aim of developing stochastic optimisation algorithms is to deal with large-scale optimisation problems. One may naturally wonder whether stochastic optimisation methods are guaranteed to be faster than their deterministic counterparts. Below are some key points. 

\paragraph{Stochastic Acceleration Factor}
Whether an inverse problem is suitable or not for the application of a stochastic gradient method largely depends on the structure of the optimisation problem. which can assist us in better understanding when and why stochastic gradient methods can be efficient or inefficient in different inverse problems scenarios.

In the works of \cite{tang2019limitation,tang2020practicality}, a metric for estimating whether the inverse problem is suitable for applying stochastic gradient methods is proposed.
This metric, named the \emph{stochastic acceleration factor}, is defined as the ratio of the two Lipschitz constants $\Upsilon
:= L / \Lmax$. In the ideal cases for stochastic methods, the ratio between the gradient Lipschitz constant is close to $\numh$. Many machine learning problems with common datasets fall into this category, as well as some important inverse imaging problems such as CT reconstruction. However, there are many negative examples with $\Lmax \approx L$, making $\Upsilon \approx 1$, in which case stochastic gradient methods are not recommended. These include compressed sensing inverse problems, image denoising or image inpainting. For image deblurring it strongly depends on the size of the kernel, see also the numerical example in Section~\ref{sec:numerics}.

In the context of linear inverse problems where the data-fit is given by the least-squares, a specialised analysis reveals that the success of stochastic acceleration depends on the structure of the forward operator. Theoretical analysis provided by~\cite{tang2020practicality} shows that an inverse problem is suitable for stochastic gradient methods if and only if the Hessian of the inverse problem has a fast-decaying spectrum. In the optimisation community, it is often overlooked that stochastic gradient methods are exploiting the intrinsic low-dimensional structure of the data matrices. If such a structure is weak in certain problems, the use of stochastic algorithms may not result in any benefits.

\paragraph{When are Variance-Reduced Methods Beneficial?} 
The importance of variance-reduced methods depends highly on the amount of noise in the measurements. It has been shown in \cite{Vaswani2019} that SGD converges linearly in the so-called interpolation regime, which in the context of inverse problems essentially equates to a noise-free measurement regime. In this case SGD can even be faster than variance-reduced methods as it can use a constant step{\removed{-}}size of $1/(\numh\Lmax)$\cite{gower2019sgd}. In contrast, variance-reduced methods typically require smaller constant step{\removed{-}}sizes. For instance, the largest known provably-convergent step{\removed{-}}size is $1/(3\numh\Lmax)$ for SAGA~\cite{Defazio2014}, and $1/(4\numh\Lmax)$ for SVRG~\cite{xiao2014}. When the measurement noise is not negligible vanilla SGD requires using a decaying step{\removed{-}}size, in order to ensure convergence, and becomes slow for high accuracy solutions. In such cases, variance-reduced stochastic gradient methods are recommended due to their fast convergence in this regime which is prevalent for inverse problems. See also the numerical experiments in Section~\ref{sec:numerics}.

\paragraph{Efficient Partial Computations}
As we discussed above, testing the stochastic acceleration factor is a necessary condition for the theoretical utility of SGD for specific inverse problems. Here we describe another test related to parallelisability. Practically, stochastic methods for inverse problems all partition the forward operator $\ipop$ into smaller blocks $\ipop_i$ such that $(\ipop x)_i = \ipop_i x$. A fundamentally important question is whether $\ipop_i x$,for $i=1,\dots,s$ can be computed as fast, or faster than $\ipop x$. Stochastic methods essentially serialise the problem which adds a computational overhead. The practical feasibility of a stochastic method then depends on both the stochastic acceleration factor and the parallelisability of the forward operator. For example, most implementations of the X-ray transform are parallelisable over its angles and the forward operator in parallel MRI is parallelisable over the coils.

\paragraph{Subset Sizes and Computational Overhead}
The number of subsets allows us interpolating between deterministic and high-variance stochastic algorithms, both of which have their benefits and drawbacks. To illustrate this point we compare PGD with SGD. In every iteration of PGD, the most costly operations are the evaluation of the gradient and of the proximal operator. In imaging inverse problems it is common to use proximal operators with a \replaced{non-negligible}{nonnegligible} computational cost, e.g.\ the proximal operator of the total variation, which can be approximated via iterations \cite{beck2009fast, rasch2020inexact} or by applying a deep neural network \cite{hurault2021gradient,tan2023provably}. 
Thus, a single data pass of SGD has the same computational cost for the gradient (as it is evaluated only on a subset), but it requires evaluating the proximal operator $\numh$ times, which often leads to a considerable overhead. Therefore, for most applications, the number of subsets $\numh$ needs to be chosen conservatively to avoid this computational overhead.

\subsection{Hyperparameter Tuning}
Parameter tuning has been an important yet long-standing challenge in stochastic optimisation. Due to their random nature, most algorithmic parameters are usually much harder to be either selected a-priori or tuned on-the-fly, compared to their deterministic counterparts. The typical algorithmic parameters which commonly need to be tuned include:
\paragraph{Step{\removed{-}}size Selection}
step{\removed{-}}size selection is an important problem in stochastic optimisation, and has been a driving   force behind the development of advanced stochastic algorithms, e.g.\ variance-reduced stochastic algorithms. 
On the one hand, stochastic algorithms inherit the step{\removed{-}}size issues that deterministic algorithms exhibit, such as requiring to know the Lipschitz constant, which may not be available, too expensive to compute accurately, or highly spatially varying. 
On the other hand, step{\removed{-}}sizes with theoretical guarantees are often too conservative and in practice, step{\removed{-}}sizes larger than the theoretical limit often work well and provide an even faster performance. For example, in the first experiment of Section \ref{sec:numerics}, step{\removed{-}}size $1/(1.75nL_{\max})$ is used for SAGA under which the algorithm converges. 
Similarly, primal-dual based algorithms require proper balancing of the primal and dual step{\removed{-}}sizes to yield good performance. 
Overall, the step{\removed{-}}size selection for stochastic optimisation algorithms needs more dedicated research. 

\paragraph{Backtracking} While the line-search schemes can usually be implemented for deterministic gradient methods such as PGD or FISTA \cite{beck2009fast}, applying these schemes for stochastic gradient methods is in general difficult, both in theory and in practice. Traditional line search methods that require evaluating the full objective are undesirable, since this would require using all the data. Using the corresponding stochastic variants can perform well in practice but lacks theoretical guarantees \cite{gower2020variance}. Recently, a number of breakthroughs have been made for on-the-fly tuning of the step{\removed{-}}size for stochastic gradient methods when applied in machine learning \cite{defazio2023learning}. However, this has not yet been fully adapted to the proximal setting which is prevalent in inverse problems. 
    
\paragraph{Subset Size} The selection of the subset size is often a tricky task in practice. Although in theory the computational complexity is optimised for stochastic methods when tiny subset sizes are chosen, in practice one has to also consider the fact that most modern computational devices are highly parallel which means that the hardware would strongly favour larger subsets, leading to a trade-off. Currently there seems to be little study on providing practical guideline on selecting subset sizes in the literature.

\paragraph{Subset Pattern} Selecting the subset pattern, i.e.\ partitioning of the forward operator, is essential for optimally applying stochastic gradient methods. As discussed in \cite{tang2020practicality}, the subset pattern determines the gradient Lipschitz constants on subset, and hence affects the admissible step{\removed{-}}size and convergence speed. The influence of the subset selection has also been discussed in~\cite{Gutierrez2022} for parallel MRI. The choice of a suitable subset pattern is reasonably well-understood for some inverse problems like full-angle 2D CT. However, for many other problems like limited-angle CT, 3D PET etc, this problem has not been resolved.

\paragraph{Restarting} Accelerated stochastic gradient methods do not automatically adapt to the strong-convexity or hidden restricted strong-convexity which is prevalent when using sparsity-inducing regularisers. Without restarting, these methods would still show fast convergence but only initially. When approaching high-accuracy solutions, the algorithms start to oscillate around the solution, slowing down the convergence. To address this, adaptive restart schemes have been developed to accelerate momentum-based methods~\cite{tang2018rest,fercoq2020restarting}.

\paragraph{Warm-starting} In practice, variance-reduced methods such as SVRG and SAGA, are slower than vanilla SGD at the start. Replacing the initial iterations of variance-reduced methods with a warm-start phase by running vanilla SGD, a faster convergence at the beginning of iterations has been observed \cite{Twyman2023rdp,konevcny2017semi}. However, it is unclear how long this phase should be a-priori.

\subsection{Memory Considerations}
The stochastic algorithms presented here have very different memory footprints which can play a significant role in the way they are used in applications. 
In the following we will discuss the impact of memory requirements, and some mitigation strategies for SGD, SAG(A), SVRG, and SPDHG, though similar considerations can be applied to other algorithms.

Vanilla SGD does not require storing any additional pieces of information, and thus has the same memory footprint as GD: \new{$\mathcal O(\dim \optspace)$} required to store the current iterate.

SAG and SAGA, in their standard form, require storing the table of computed gradients. The sum of the elements in the table then either needs to be stored or computed anew in each iteration.
    Thus, the memory overhead of the standard implementation of SAG and SAGA is \replaced{$\mathcal O(\numh |\optspace|)$}{$\mathcal O(\numh \dim \optspace)$}, which can be prohibitively large if $\numh$ and \replaced{$|\optspace|$}{$\dim \optspace$} are large. In particular, in large scale CT the dimensionality \replaced{$|\optspace|$}{$\dim \optspace$} can be very large. 
    However, memory requirements can be significantly lowered in certain mitigating scenarios. Consider linear models \cite{Defazio2014, gower2020variance}, i.e.\ when $h_i=\overline{h}_i\circ \optop_i$ for $i=1,\ldots,\numh$ and linear operators $\optop_i$. 
    This is for example the case for the least-squares data fidelity or the Kullback--Leibler divergence (see Section \ref{sec:ip}).
    The gradient of $h_i$ then satisfies $\nabla h_i(\optvar)=\optop_i^\ast \nabla\overline{h}_i(\optop_i\optvar)$.
    In other words, instead of storing $\nabla h_i(\optvar)$ we can store $\nabla\overline{h}_i(\optop_i\optvar)$. As in the original SAGA, each iteration needs one application of $\optop_i$ and $\optop_i^*$, thus no extra computational cost is introduced.
    The modified SAGA algorithm is written in Algorithm \ref{alg:modified_saga}.
    \begin{center}
\begin{minipage}{0.95\linewidth} 
\begin{algorithm}[H] 
    \caption{Modified SAGA}\label{alg:modified_saga}
        \DontPrintSemicolon
        \KwIn{$x^{(0)} \in \optspace$, $\tau>0$}
        Initialise $y^{(0)}_i = \nabla \overline{h}_i(\optop_i\optvar^{(0)})$ and $\gamma^{(0)}=\sum_i \optop_i^\ast y^{(0)}_i$\\
        \For{$k=0,1,\ldots$}{ 
            Sample $i_k \in\{1,\ldots,\numh\}$ uniformly at random\\
            $\overline y=\nabla\overline{h}_{i_k}(\optop_{i_k}\xk)$\\
            $\tilde\nabla=\optop_{i_k}^\ast(\overline y - y_{i_k}^{(k)}\big)$ and $\gamma^{(k+1)}=\gamma^{(k)}+\tilde\nabla$\\
            $\tilde\nabla h(\xk)=\numh\tilde\nabla+\gamma^{(k)}$ \\
            $\xkp=\prox_{\tau g}\big(\xk-\tau\tilde\nabla h(\xk)\big)$\\
            $y_{j}^{(k+1)} = \begin{cases} \overline y & \textrm{if~~} j = i_k , \\ y_{j}^{(k)} & \textrm{otherwise} \end{cases}$
        }
\end{algorithm}
\end{minipage}
\end{center}

\new{The storage requirements of the modified SAGA are \replaced{$\mathcal O(|\optspace|+ |\dualspace|)$}{$\mathcal O(\dim \optspace + \dim \dualspace)$}. Note that the stored sum of the gradients $\boldsymbol{\gamma}^{(k)}$ is in $\optspace$, whereas dual variables $y_i^{(k)}$ are in $\dualspace_i$. Thus, storing all $y_i^{(k)}$, for $i=1,\ldots,n$ can be done with memory $\dim \dualspace$. Moreover, we can exploit the linearity of the model so that, as in the original SAGA, evaluation of the adjoint of only one of the operators $\optop_i$ is needed in each iteration.} 

The choice of the preferred variant of SAGA depends on the imaging modality and the number of subsets chosen. For example, in CT images and data are often of similar sizes so the modified SAGA algorithm would be preferable already for small $\numh$. For PET, the image volumes are often a magnitude small than the data, so for the setting as in Section~\ref{sec:ip}, the modified algorithm starts to become preferable for subset numbers exceeding $15$.

The standard implementation of SVRG (see Algorithm \ref{alg:svrg}) requires storing the reference iterate $\overline\optvar^{(j)}$ and the full gradient $\overline{\gamma}^{(j)}$, resulting in a memory of \replaced{$\mathcal O(|\optspace|)$}{$\mathcal O(\dim \optspace)$} which is independent of the number of subsets and comparable to GD. However, this comes at the cost of computing two stochastic gradients per iteration, one at the previous iterate and another at the reference iterate. In some applications, the memory overhead of storing $\nabla h_i(\overline\optvar^{(j)})$ might be preferable over the extra gradient computation $\nabla h_i(\overline\optvar^{(j)})$. As for SAG(A), we can derive a modified version of SVRG which has the same memory footprint \replaced{$\mathcal O(|\optspace|+ |\dualspace|)$}{$\mathcal O(\dim \optspace + \dim \dualspace)$} and the same computational cost as the modified SAGA.
    
SPDHG has essentially the same memory footprint as PDHG: \replaced{$\mathcal O(|\optspace|+ |\dualspace|)$}{$\mathcal O(\dim \optspace + \dim \dualspace)$}. Depending on the dual space, the memory footprint for primal-dual methods can be significantly larger compared to GD. Note that the computational effort and memory footprint is the same as the modified version of SAGA.

\section{Numerical Examples} \label{sec:numerics}

In this section we present a brief numerical study, showcasing representative stochastic optimisation algorithms discussed in this review, and comparing them against their deterministic counterparts. 
Note that this numerical study is not intended to provide a comprehensive and thorough comparison of algorithmic performance.
Instead, our aim is to show that stochastic algorithms can indeed be beneficial for certain problems and to illustrate their behaviour. 
{\color{red}The code to reproduce all results in this paper will be made available upon acceptance of this manuscript.}

\subsection{Sparse Spikes Deblurring}\label{subsec:sparse_spikes}

We consider a toy deblurring problem where the ground truth are sparse spikes and the forward operator uses a uniform blurring kernel.
The goal of this example is to compare PGD and SAGA and investigate how is the performance of those two algorithms affected by the size of the blurring kernel and step{\removed{-}}size regime. The numerical results were produced with MATLAB.

Let $\reconvar^\dagger \in \R^\recondim,\ \recondim= 1,000$ be a sparse vector with $20$ evenly distributed \replaced{non-zeros}{nonzero} entries, representing the ground truth, $\ipop \in \R^{\recondim\times \recondim}$ be the matrix representation of a blur with uniform kernel of size $\kappa$, with periodic boundary conditions, which is normalised such that $\|\ipop\| = 1$. 
We consider the recovery of $\reconvar^\dagger$ from noisy data $\datavar = \ipop\reconvar + \eta $, where $\eta \in \mathbb{R}^\recondim$ is random Gaussian noise. 
We consider kernels of size $\kappa_{\new{j}} = 2\replaced{i}{j}-1$ for 
$\replaced{i}{j}=1,\ldots,100$, and random Gaussian noise $\eta\sim\mathcal{N}(0,\sigma^2 I)$, with $10$ different noise levels, taken at equidistant points in $[10^{-4}, 10^{-1}]$. Note that for $\replaced{i}{j}=1$ the kernel matrix $\ipop$ is the identity matrix, and as $\replaced{i}{j}$ increases so does the correlation between rows of $\ipop$.

To recover $\reconvar^\dagger$ from the observation $\datavar$, we minimise $\Phi(\optvar) = g(\optvar) + h(\optvar)$ with $\optvar = \reconvar$, $g(\optvar) = \mu\norm{\optvar}_1$ and $h(\optvar) = \frac{1}{2} \norm{\ipop\optvar - v}_2^2$. 
The objective $\Phi$ is for SAGA rewritten into the finite-sum form $\Phi(\optvar)=g(\optvar)+ \sum_{i=1}^{\numh} h_i(\optvar)$ where $\numh=\recondim$, $ h_i(\optvar):= \frac{1}{2} (\ipop_i \optvar - \datavar_i)^2$ and $\ipop_i\in\R^{1\times \numh},i=1,...,\numh$ is the $i$-th row of $\ipop$. 
Note that $L = 1$, $\norm{\ipop_i} = 1/\sqrt{\kappa}, i=1,...,\numh$ and thus $L_{\max} = 1/\kappa$\new{, for a given kernel size $\kappa$}. 
The regularisation parameter $\mu$ is chosen as 
$\mu = \frac 12\norm{\ipop^* v}_{\infty} $. We use two step{\removed{-}}size choices for both methods. 
For PGD we use $\tau = 1/L$ (denoted PGD1) and $\tau = 1.9/L$ (denoted PGD2). For SAGA, we use $\tau = 1/(3\numh \Lmax)$ (denoted SAGA1) and $\tau = 1/(1.75\numh \Lmax)$ (denoted SAGA2). 

\begin{figure}[htbp]
    \centering
       \includegraphics[width=0.6\textwidth]{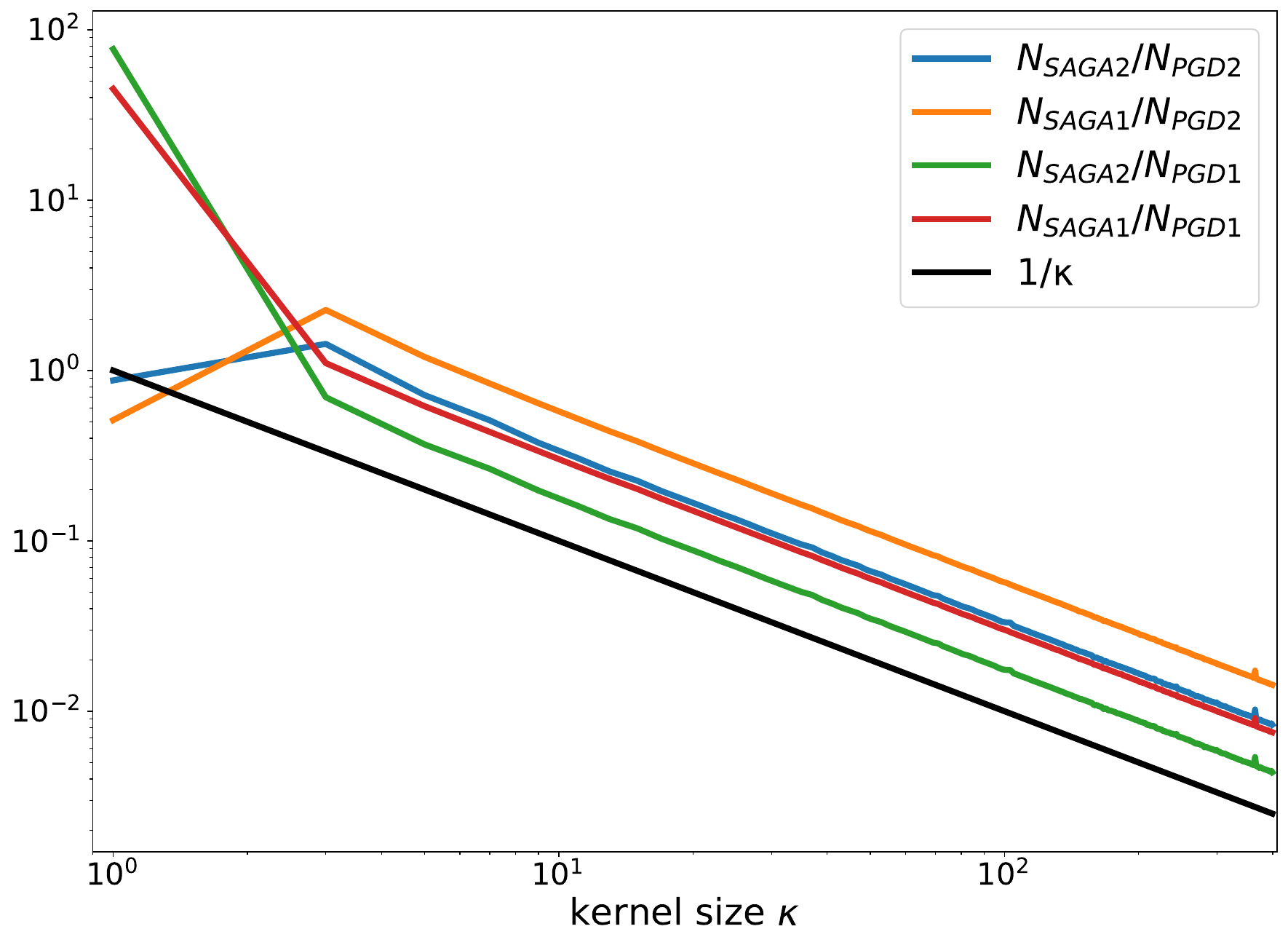}
    \label{subfig:cmp_ssd_averaged}
    \caption{The lines represent the ratio of the number of data passes for SAGA and PGD, i.e.\ $N_{\text{SAGA}} / N_{\text{PGD}}$, to reach a target reconstruction accuracy for the deblurring of sparse spikes, for different step{\removed{-}}size choices. See accompanying text in Section \ref{subsec:sparse_spikes} for more details.}  
    \label{fig:cmp_ssd}
\end{figure}

For each kernel size and noise level we compute a high-precision solution~$\optvar^\ast$ and then run PGD and SAGA using the stopping criterion 
$
{\norm{\xk-x^\ast}} / {\norm{x^\ast}} \leq 10^{-4}.
$
We record the number of data passes needed to reach that accuracy, denote it as $N_{\text{PGD}}$ and $N_{\text{SAGA}}$ for the corresponding choices of PGD and SAGA, respectively, and compute the ratio $N_{\text{SAGA}} / N_{\text{PGD}}$, which is averaged over the 10 noise levels.
The resulting $4$ ratios (for each step{\removed{-}}size pairing), and a line $1/\kappa_{\replaced{ i }{ j }}$, for 
$i=1,\ldots,100$ are shown in Figure \ref{fig:cmp_ssd}.  

We can reach the following observations from the results. First, when $\ipop$ is the identity matrix (for $\kappa=1$\new{)}, PGD with $\tau=1/L$ converges in one step. Hence, the ratio $N_{\text{SAGA}} / N_{\text{PGD}}$ is very large. Moreover, SAGA1 is faster than SAGA2. 
This behaviour can be attributed to the properties of $\ipop$.
When \replaced{$\kappa$}{the kernel size} is small the gradients $\nabla h_i$ and rows of $\ipop$ are less correlated,  decreasing the data redundancy. In such a setting stochastic methods can lose some of their competitive edge. On the other side of the coin is the case of large \removed{$\kappa$}{kernel sizes}. In this setting the the rows of $\ipop$ are more correlated, which benefits stochastic methods since each $\ipop_i$ carries a much larger proportion of the informational content of $\ipop$. In such a setting stochastic methods are expected to perform better.
Moreover, larger step{\removed{-}}size leads to a larger initial variance and slower performance. 
Second, for kernel sizes between $5$ and $199$, we observe that the ratio $N_{\text{SAGA}} / N_{\text{PGD}}$ is proportional to $1/\kappa$. 
Moreover, larger step{\removed{-}}size leads to faster performance for both PGD and SAGA. To see this, notice that the orange line $N_\text{SAGA1} / N_\text{PGD2}$ is above the red line $N_\text{SAGA1} / N_\text{PGD1}$, indicating that PGD1 needs more data passes to reach stopping criterion. Similarly, the green line stays below the red line, indicating that compared to SAGA2, SAGA1 requires more data passes to reach the stopping criterion. 

\subsection{Simulated Parallel-Beam CT}\label{subsec:parallel_CT}
We reproduce here the numerical results from an upcoming publication on a software framework for stochastic optimisation methods for large-scale imaging problems within the Core Imaging Library (CIL) \cite{CIL_I, CIL_II,CILZenodo}. This is being developed as a community effort around the CCPi (\url{www.ccpi.ac.uk}) and SyneRBI (\url{www.ccpsynerbi.ac.uk}) projects. The reproduced results deal with the CT reconstruction of a Shepp--Logan phantom of size $128\times128$ px. 
The measurements follow the forward problem $\datavar=\ipop\reconvar^\dagger+\eta$, where $\eta$ is a random sample from $\mathcal{N}(0,\sigma^2 I)$, and the forward operator $\ipop$ is defined by parallel-beam CT with $240$ angles from the $[0, \pi)$ angle range.

We compare the performance of a deterministic algorithm PGD (Algorithm~\ref{alg:pgd}) with SGD (Algorithm~\ref{alg:sgd}) and SAGA (Algorithm~\ref{alg:saga}).
As discussed in Section~\ref{sec:stoch}, when solving inverse problems it is not common to utilise fully stochastic approaches, in the sense that instead of using randomly constructed batches the  measurement data is pre-partitioned into batches. One element of the partition is then sampled at random in each iteration.

To recover the phantom $\reconvar^\dagger$ we minimise $\Phi(\optvar)=g(\optvar)+h(\optvar)$, where $\optvar=\reconvar$, $g(\optvar)=\lambda\text{TV}(\optvar)$, and $h(\optvar)=\frac{1}{2}\|\ipop\optvar-\datavar\|^2$ for PGD.
SGD and SAGA rewrite the objective as $\Phi(\optvar)=g(\optvar)+\sum_{i=1}^{\numh} h_i(\optvar)$, where $\numh\in\mathbb{N}$, and $h_i(\optvar)=\frac{1}{2}\|\ipop_i\optvar-\datavar_i\|^2$ is such that $\sum_{i=1}^\numh h_i(\optvar)=h(\optvar)$. 
Operator $\ipop_i$, for $1\leq i\leq \numh$, is defined as the CT operator that uses $240/\numh$ angles from the range $[(i-1)/\numh, \pi)$ by taking every $\numh$-th angle starting with $i/\numh$ with increments of size $\pi\numh/240$. This is sometimes called the \emph{equidistant} or \emph{staggered} partition, and results in all the operator norms $\|\ipop_i\|$ to be of a nearly identical value. Moreover, the number of batches is typically a divisor of the number of angles, which ensures that all the batches are of the same size.

\begin{figure}[htbp]
    \centering
    \subfloat[Objective value over 200 data passes]{\includegraphics[width=0.48\textwidth]{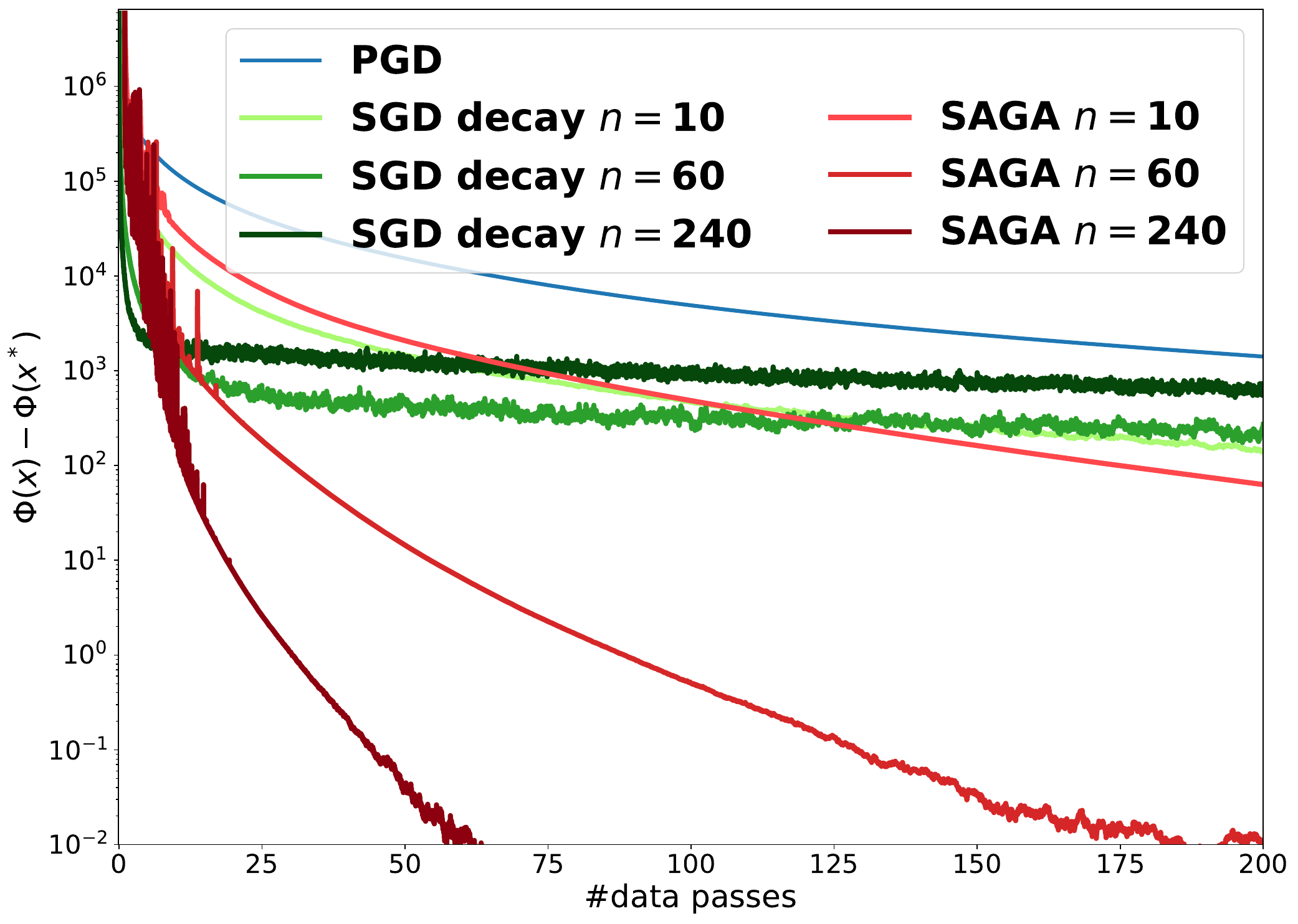}\label{fig:SL_objectives_200}}
    \hfill
    \subfloat[Objective value over the first 20 data passes]{\includegraphics[width=0.48\textwidth]{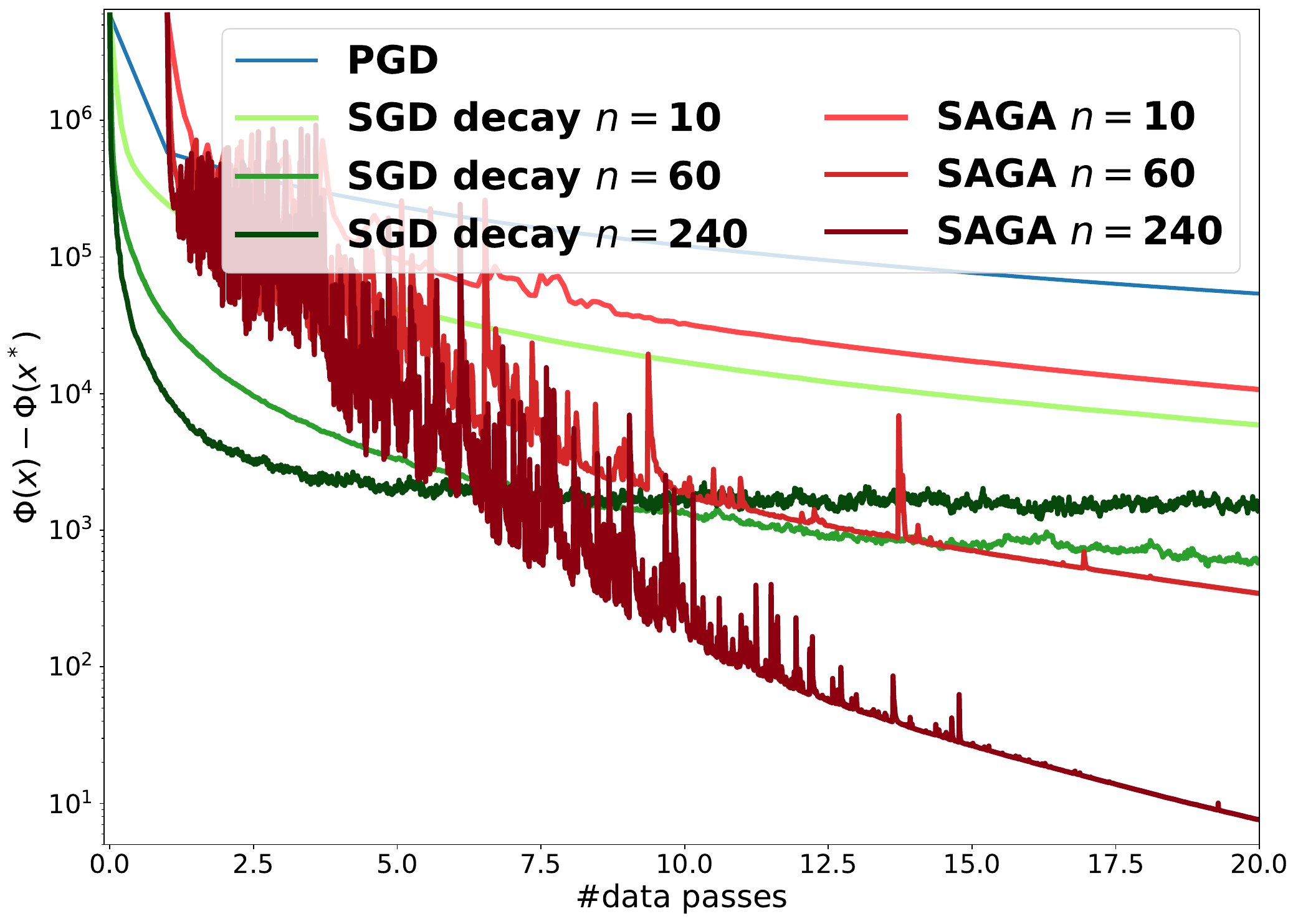}\label{fig:SL_objectives_20}}
    \\
    \includegraphics[width=\textwidth]{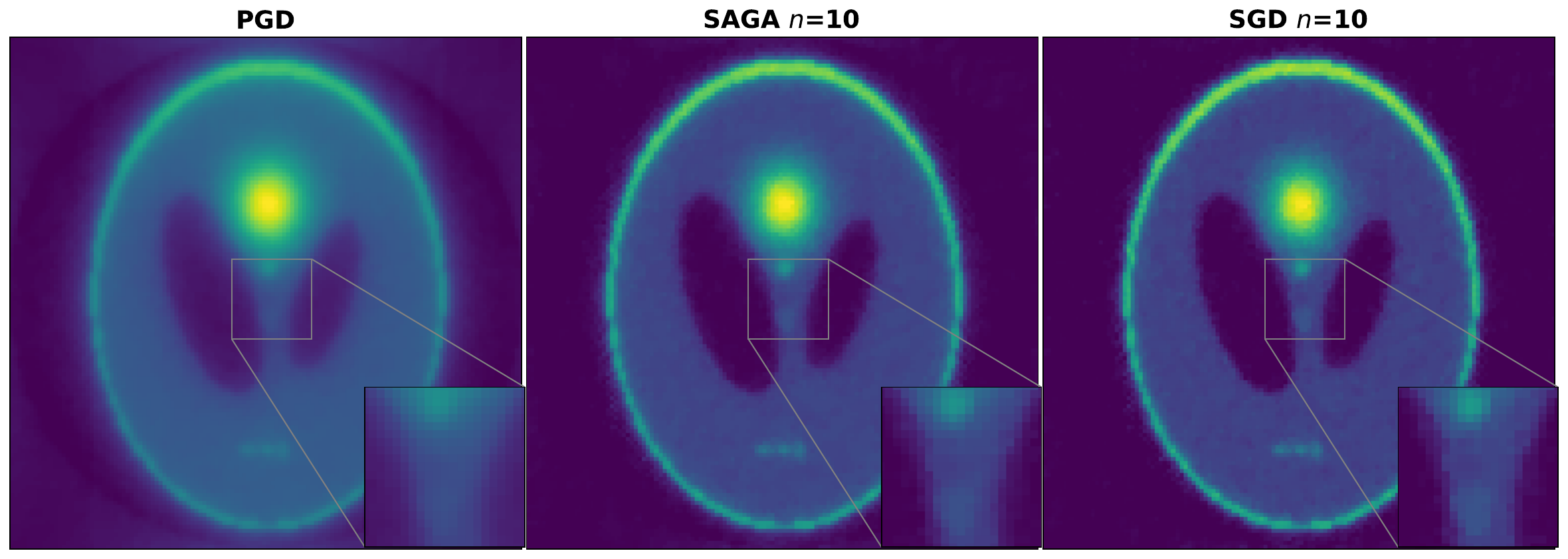}
    \\
    \subfloat[Reconstructions with PGD, SGD, and SAGA after 10 epochs]{\includegraphics[width=\textwidth]{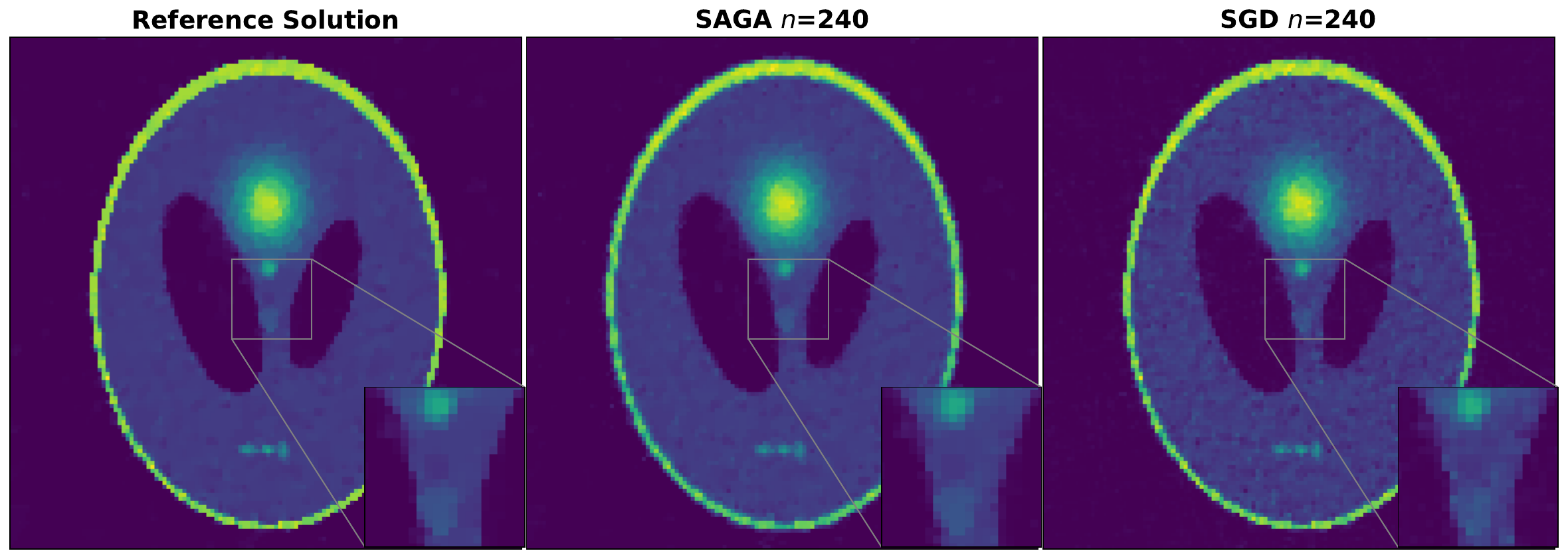}\label{fig:SL_reconstructions}}
    \caption{Performance comparison for a Shepp--Logan phantom from CT measurements corrupted by Gaussian noise with $\sigma=1$. Note that the objective value for SAGA starts at the $1$-st data pass due to initialising the table of gradients at $\optvar^{(0)}$.}    
    \label{fig:SL_results}
\end{figure}

The initial image is $\optvar^{(0)}=0$.
The step{\removed{-}}size for PGD is $\tau=1/\|\ipop\|^2$. 
For SGD we use a decaying step{\removed{-}}size schedule, $\tau_k= 1/(2\numh \Lmax (1+0.01k/\numh))$, where $\Lmax=\max_i \|\ipop_i\|^2$, following \cite{Twyman2023rdp}, which satisfies condition \eqref{eqn:sgd_stepsize_condition}.
In this regime the step{\removed{-}}sizes decay from $1/(2\numh \Lmax)$  to roughly \new{$1/(6\numh \Lmax)$}, over the studied \new{$200$} data passes, giving a behaviour similar to a constant step{\removed{-}}size regime, albeit convergent.
For SAGA we use $\tau=1/(3\numh \Lmax)$ as step{\removed{-}}size. 

In Figure \ref{fig:SL_results} we show the performance of these 3 algorithms in terms of the sub-optimality of the objective value: $\Phi(\xk)-\Phi(\optvar^\ast)$, where $\Phi(\optvar^\ast)$ is a reference solution computed with FISTA (Algorithm \ref{alg:fista}).
Figure \ref{fig:SL_results} shows a significant speed up with SGD and SAGA over PGD. 
SGD is the fastest method in early iterations, but it then slows down after $5$-$10$ data passes.
On the other hand, SAGA is steady and fast throughout the optimisation trajectory and
the increased variance it presents in initial iterations disappears already after around $20$ data passes. 
Figure \ref{fig:SL_objectives_200} and Figure \ref{fig:SL_objectives_20} also investigate the effect of the number of subsets: the lightest shade shows a low number of subsets ($\numh=10$), normal shade shows a higher number of subsets ($\numh=60$), and the darkest shade shows the highest number of subsets ($\numh=240$).
We observe that increasing the number of subsets has a dual effect: it increases the variance of the optimisation trajectory for stochastic algorithms, but it greatly increases the speed. 
In Figure \ref{fig:SL_reconstructions} we show the reconstructions of the 4 studied algorithms obtained after $10$ data passes. 
The reconstruction with PGD presents as oversmooth, and fails to adequately reconstruct the finer image detail. For example, consider the zoomed in region where stochastic algorithms manage to recover two distinct objects (one of a lower and the other of a higher intensity), whereas PGD does not. 

\begin{figure}[htbp]
    \centering
    \subfloat[SGD with $\numh=30$]{\includegraphics[width=0.33\textwidth]{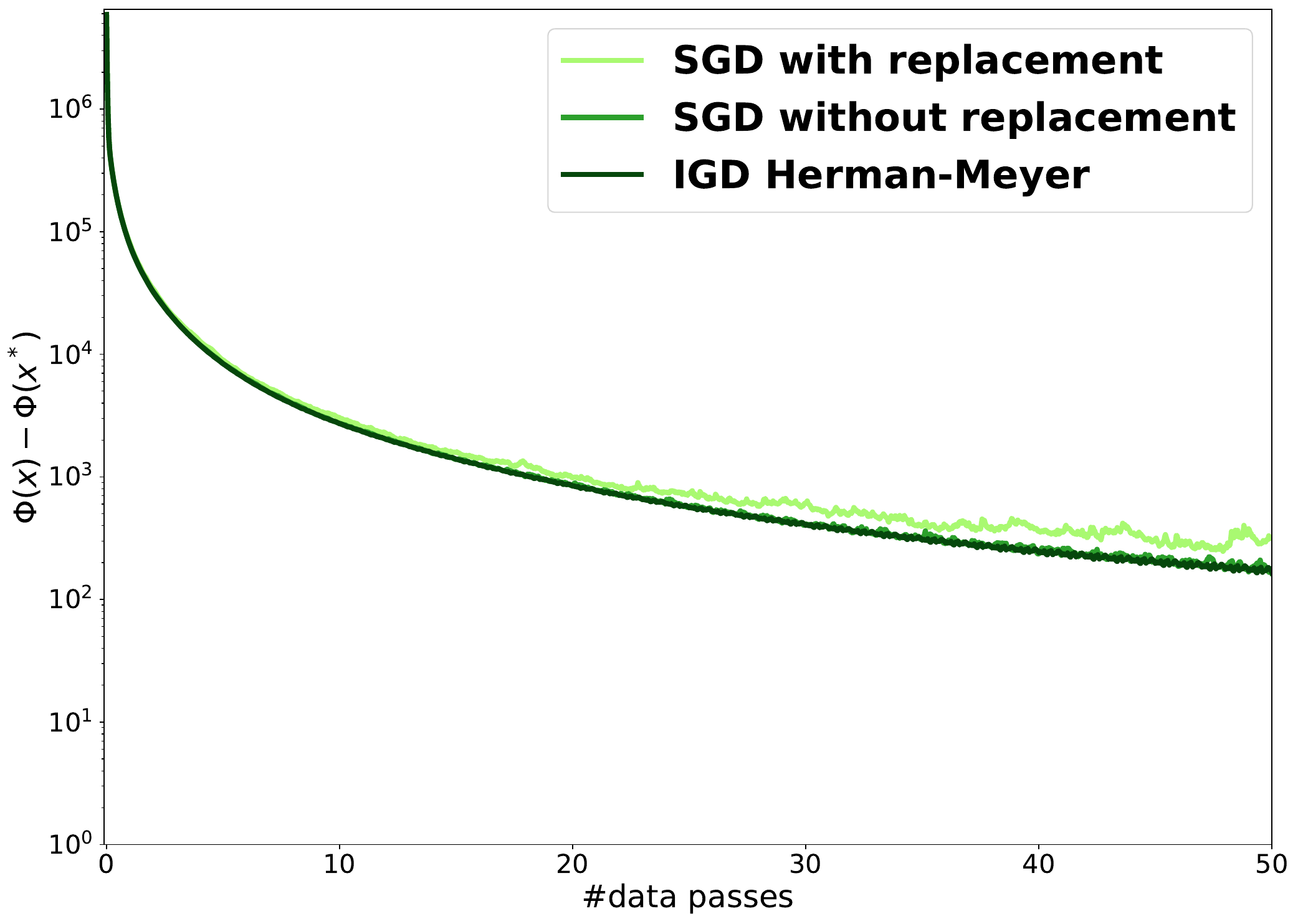}\label{fig:SL_figure1}}
    \hfill
    \subfloat[SAGA with $\numh=30$]{\includegraphics[width=0.33\textwidth]{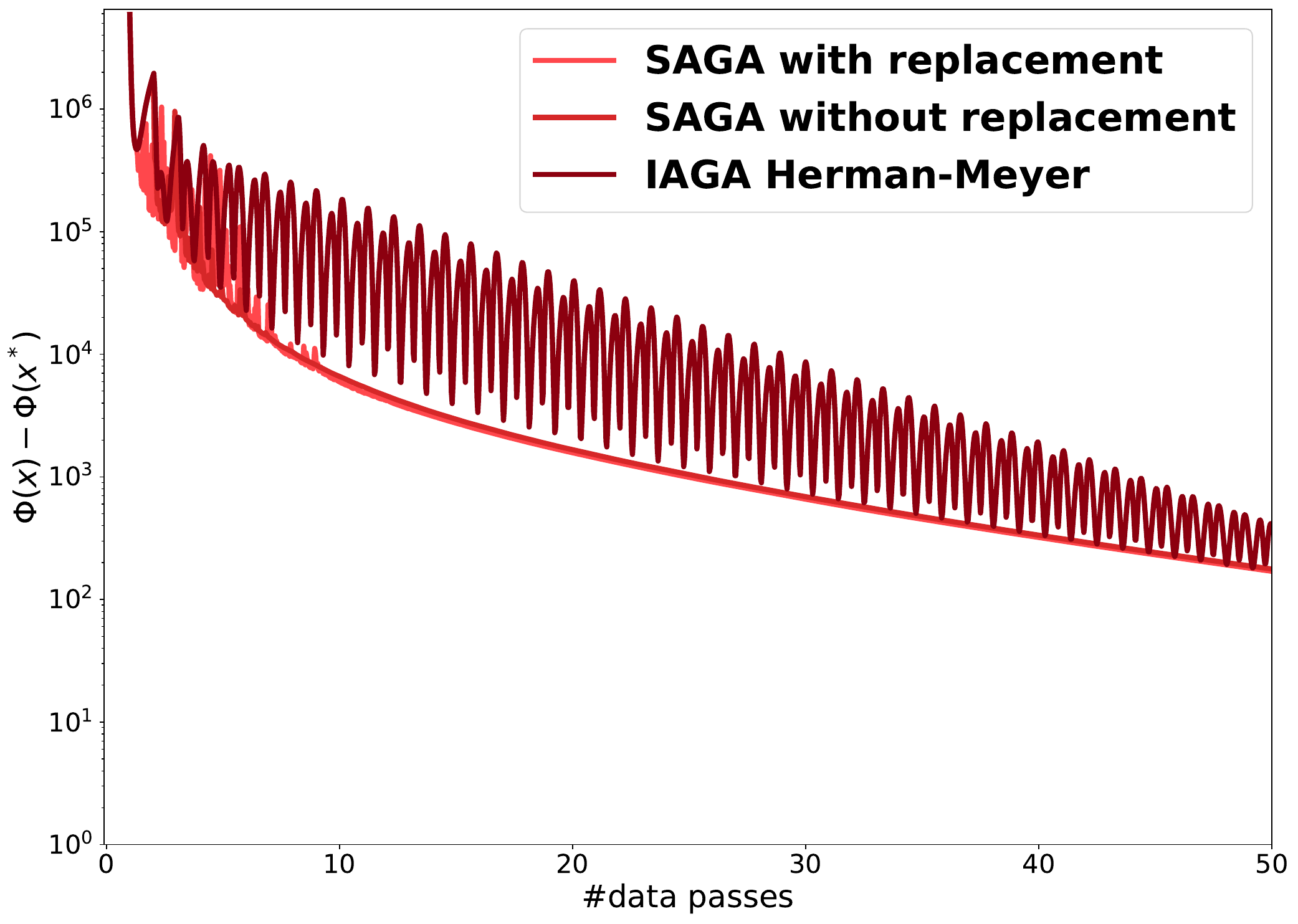}\label{fig:SL_figure2}}
    \hfill
    \subfloat[SAGA with $\numh=240$]{\includegraphics[width=0.33\textwidth]{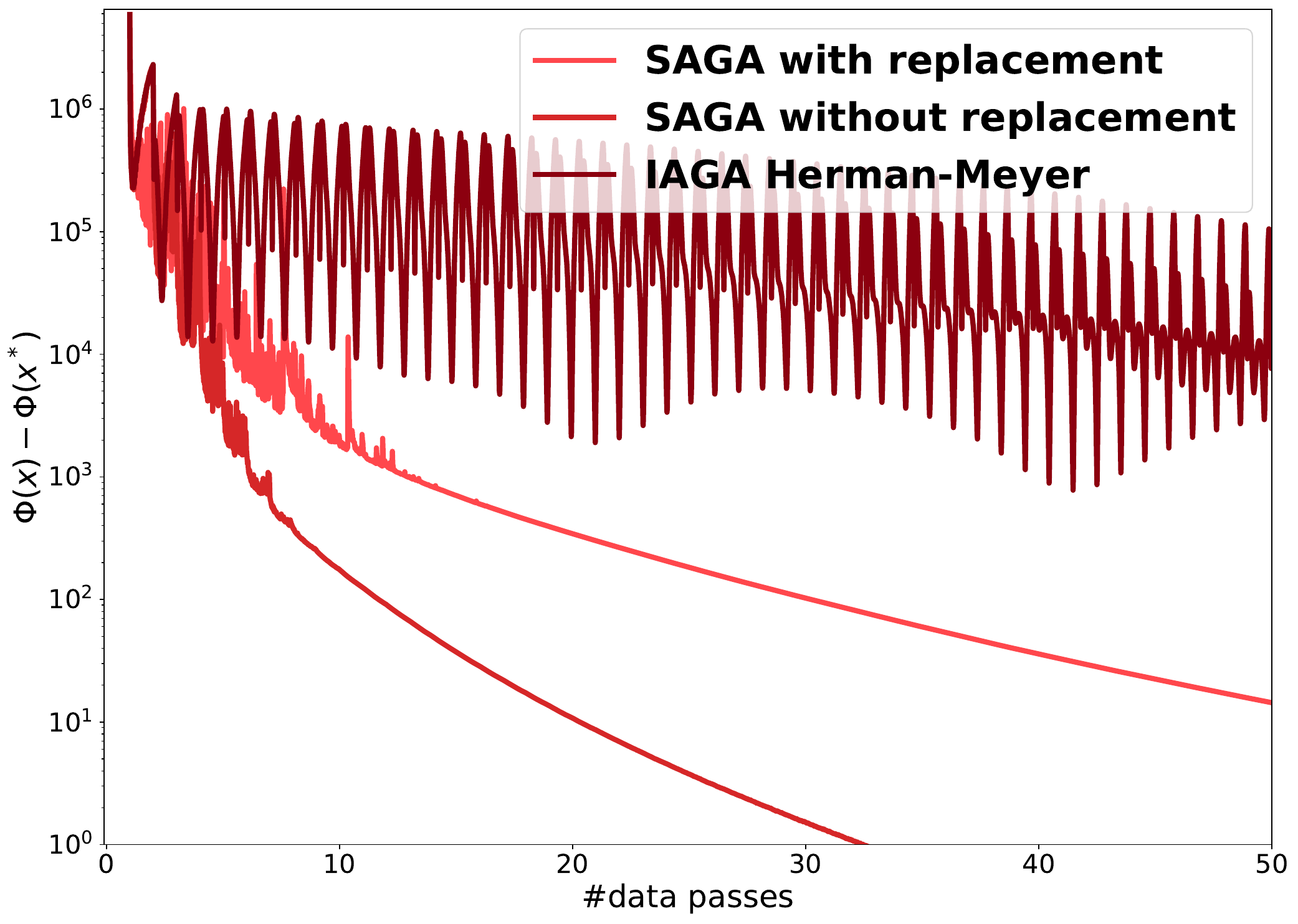}\label{fig:SL_figure3}
    }
    \caption{Comparison of the performance of stochastic and incremental methods (using the Herman--Meyer order) on parallel-beam  CT reconstruction of the Shepp--Logan}
    \label{fig:SL_sampling}
\end{figure}

In Figure \ref{fig:SL_sampling} we compare the effect of the strategy for accessing the function indices. 
Specifically, we compare SGD and SAGA, with indices sampled uniformly at random (with \new{and without} replacement), and IGD and IAGA, which are incremental gradient formulations of SGD and SAGA that access indices according to a deterministic Herman--Meyer order~\cite{herman1993hmorder}.
The results show that the behaviour highly depends on the number $\numh$ of subsets and the used algorithm. 
Comparing SGD and IGD reveals that the deterministic strategy can outperform stochastic \new{sampling strategies}, though the differences are not significant.
However, for variance reduced algorithms the differences are significantly more severe. 
When $\numh$ is small, SAGA outperforms IAGA and presents a much smaller variance, though the variance with IAGA reduces significantly within the first $50$ data passes.
When $\numh$ is large, SAGA is significantly better than IAGA, \new{and sampling with replacement can outperform sampling without replacement}.

\begin{figure}[htbp]
    \centering
    \subfloat[SGD and ADAM with $\numh=10$]{\includegraphics[width=0.5\textwidth]{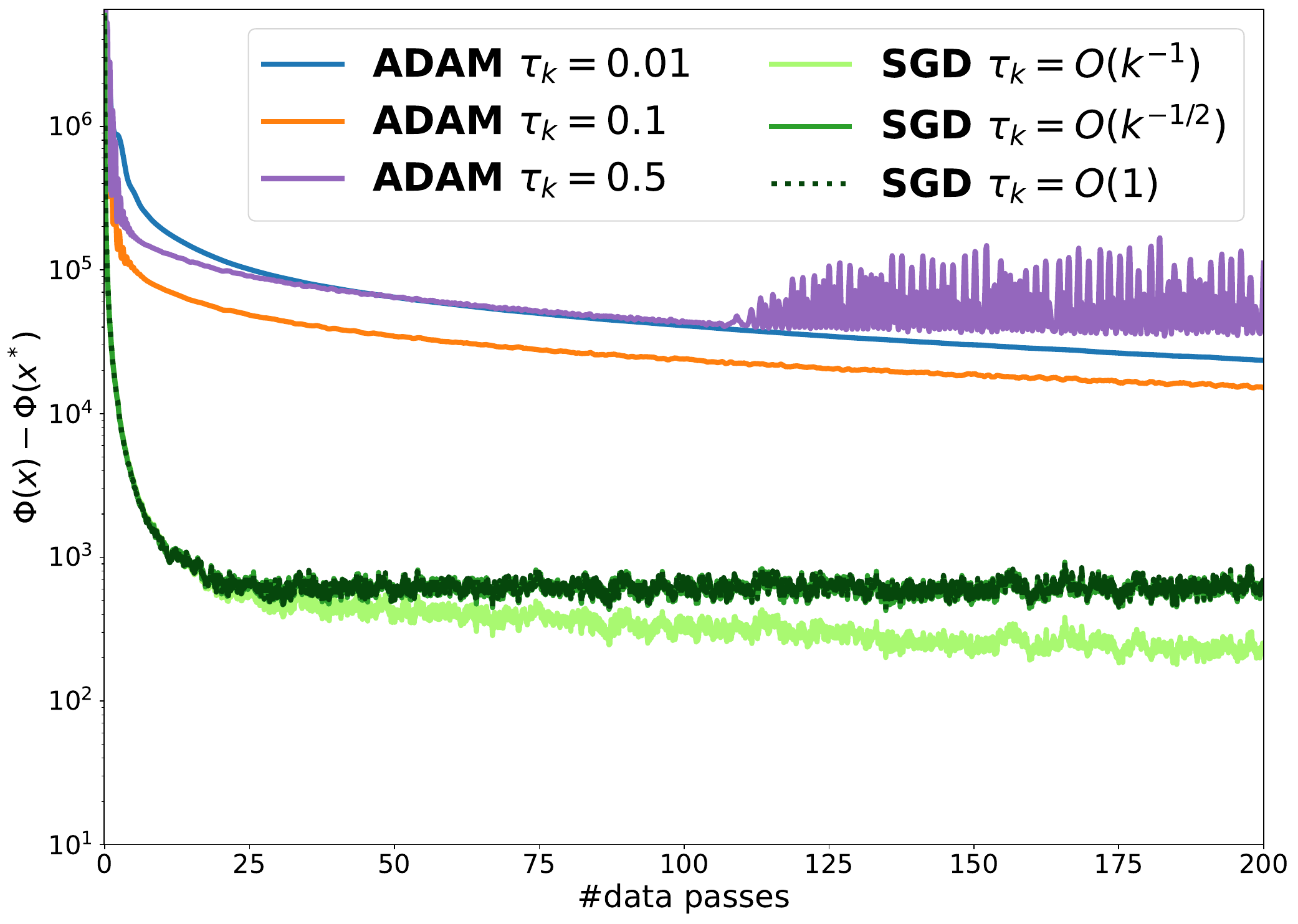}\label{fig:SL_stepsize1}}
    \hfill
    \subfloat[SGD and ADAM with $\numh=60$]{\includegraphics[width=0.5\textwidth]{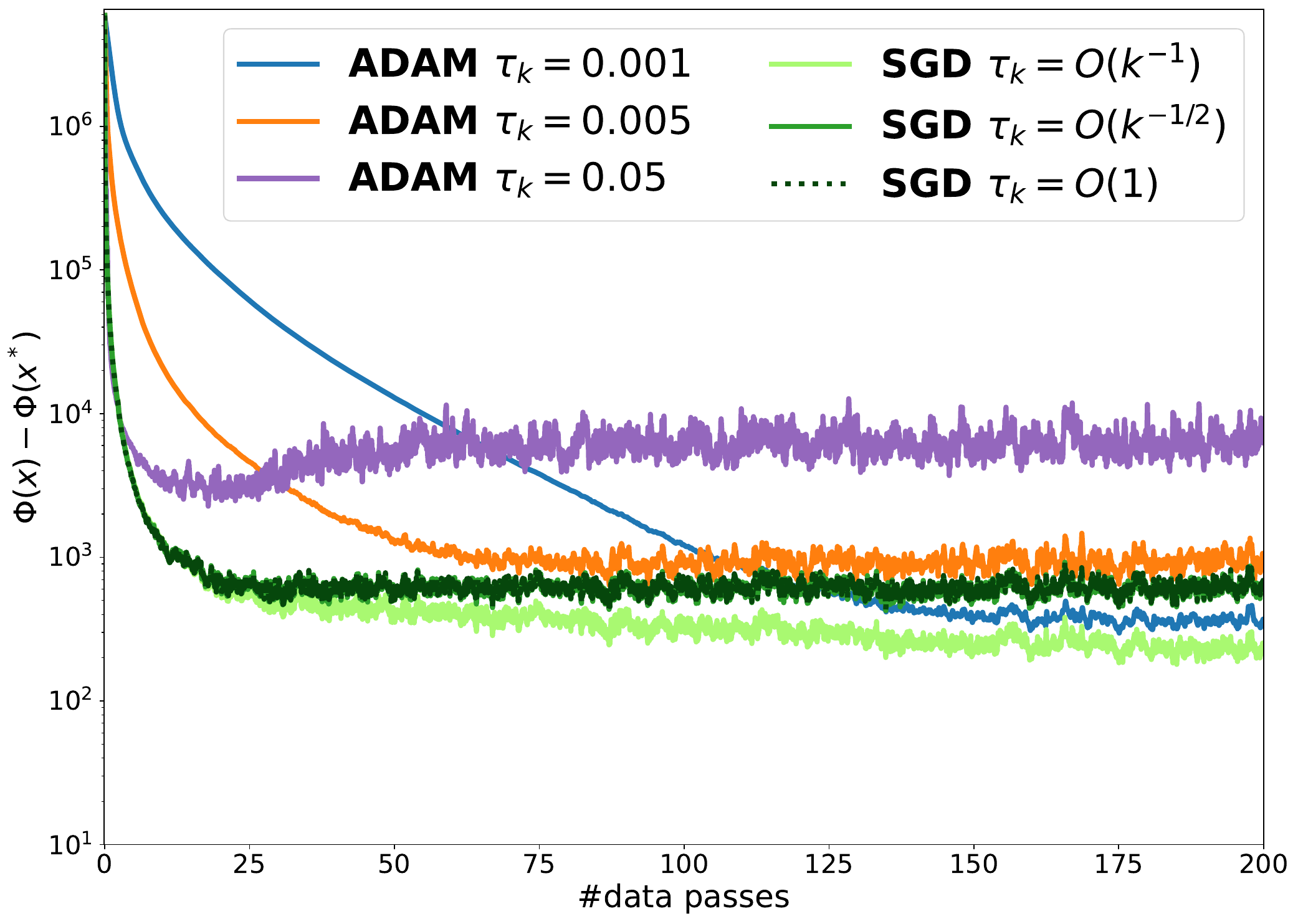}\label{fig:SL_stepsize2}}
    \hfill
    \caption{\new{Comparison of the performance of SGD with respect to the stepsize regime, and a comparison with ADAM on parallel-beam CT reconstruction of the Shepp-Logan}}
    \label{fig:SL_stepsize}
\end{figure}
\new{
In Figure \ref{fig:SL_stepsize} we investigate the effect of the stepsize schedule for SGD with respect to the number of subsets. We examine three options, each with a progressively more aggressive stepsize: constant stepsize $\tau_k=1/(2n\Lmax)$, sublinearly decaying stepsize $\tau_k= 1/(2\numh \Lmax (1+0.01\sqrt{k}/\numh))$, and linearly decaying stepsize $\tau_k= 1/(2\numh \Lmax (1+0.01k/\numh))$.
The results indicate that in a regime with fewer subsets (and consequently lower variance), larger stepsizes can be utilised. However, variance increases with the number of iterations (as it is not properly kept in check by stepsizes). 
This effect is amplified as the number of subsets (and thus variance) grows. Further hyperparameter tuning could enhance method performance, but such an investigation is beyond the scope of this review.
Note that constant stepsize and sublinear stepsize SGD are nearly identical, which is due to a very slow decay of the stepsizes over 200 epochs.

Lastly, we compare SGD with ADAM, a popular optimisation algorithm in machine learning that extends SGD by incorporating momentum, exponentially decaying averages, and estimators of second-order terms. ADAM has become the go-to optimisation tool due to its effectiveness in smooth optimisation tasks. 
ADAM is defined for smooth optimisation task and it is typically not studied for problems with non-smooth terms. Consequently, tuning its hyperparameters for inverse problems applications can be difficult.
We use ADAM with default momentum parameters, and investigated two standard stepsize regimes: constant stepsize and sublinearly decaying stepsizes. The results shown in Figure \ref{fig:SL_stepsize} were found by parameter search and were selected for balance between convergence and speed\footnote{For smaller stepsize the convergence was slower, whereas for larger stepsizes the variance increased in early epochs}. 
As for SGD, constant and decaying stepsize ADAM have a very similar behaviour due to a slow decrease of the stepsizes, and thus we only show the results for different values of the constant stepsize.
A study of ADAM and similar adaptive step-size, and the reasons why it underperforms on large scale inverse problems, has not been fully investigated in the literature.
}

\subsection{Fan-Beam CT Data of a Walnut}
In this set of experiments we take the central slice of a cone beam CT scan of a walnut~\cite{zenodoXrayDatasetWalnut}.
As in Section \ref{subsec:parallel_CT}, we reproduce the numerical results from an upcoming publication on the stochastic framework in the CIL software library.
Full sinogram data has 721 projections in the $[0, 2\pi)$ range.
To reduce the computational costs we consider only $360$ projections in the $[0, 2\pi)$ range, with a $1^\circ$ angle separation. 
The reconstructed images are of size $280\times 280$ px.

\begin{figure}
    \centering
    \includegraphics[width=\textwidth]{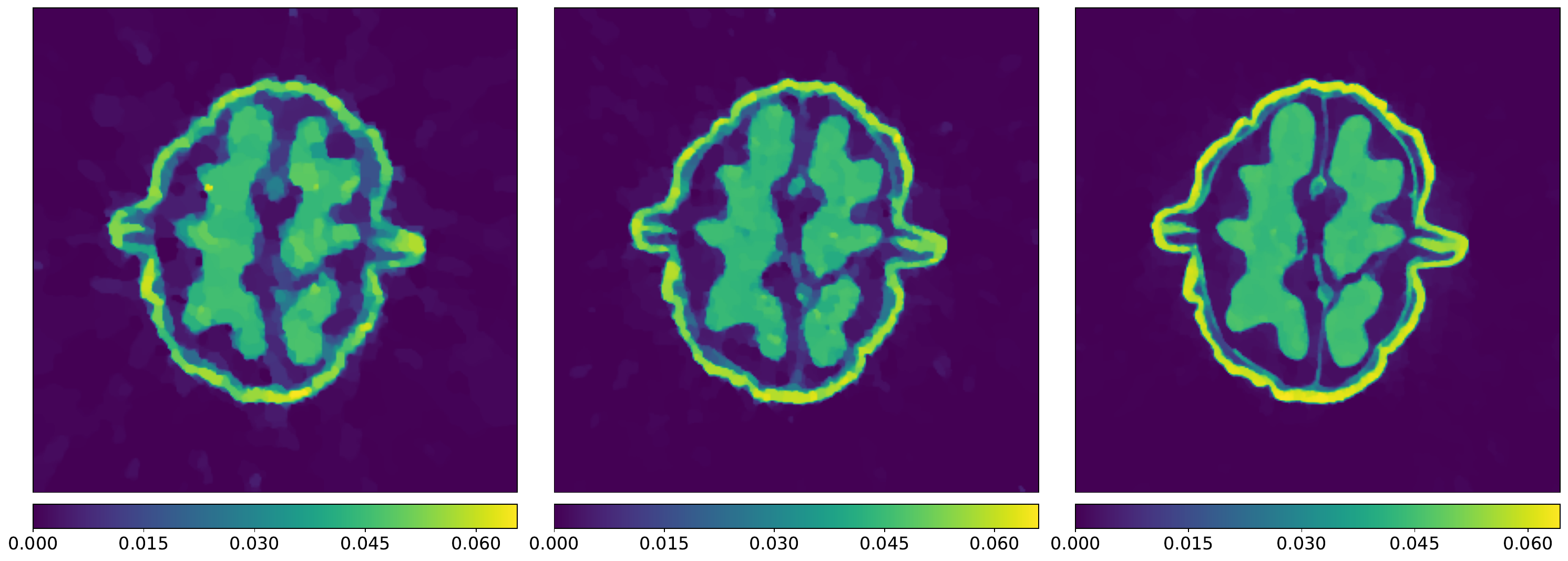}
    \caption{Reference solutions for the walnut with high, medium and low noise levels}
    \label{fig:walnut_referencesolns}
\end{figure}

Noisy measurements are generated by corrupting the sinogram (the measurement data) following the Beer--Lamberts law. 
We first compute the expected number of counts via $N_I=\text{Pois}(I_0\exp(-\datavar))$, where $\datavar$ is the loaded (noiseless) sinogram and $I_0$ the intensity of the initial beam. 
The noisy measurement is then computed by first correcting all $0$ entries in $N_I$ to $1$, and then applying linearisation (post-log measurements) as $\datavar^\delta = -\ln(N_I/I_0)$.
We consider three intensity levels: low intensity ($I_0=50$), medium intensity ($I_0=250$), and high intensity ($I_0=5000$), simulating high to low noise.  
The objective is defined as 
\begin{align}\label{eqn:walnut_objective}
\Phi(\optvar)=\frac{1}{2} \|\ipop\optvar-\datavar\|_2^2 +\alpha_I \|\nabla\optvar\|_1 + \iota_{[0,\infty)^d}(\optvar),
\end{align}
where $\alpha_I\in\{0.124, 0.045, 0.01\}$ is the regularisation parameter for high, medium, and low noise levels, respectively.
Reference solutions are shown in Figure \ref{fig:walnut_referencesolns}.
We study gradient based methods (PGD, SGD, SAGA, SVRG) and primal-dual methods (PDHG and SPDHG).
The terms in the objective \eqref{eqn:walnut_objective} are identified with optimisation templates~\eqref{eq:opttemplate} and~\eqref{eq:opttemplate_finite_sum}, for gradient based and primal-dual methods as follows
\begin{description}
\item[Gradient based methods] $g(\optvar)=\alpha_I \|\nabla \optvar\|_1+\iota_{[0,\infty)^d}(\optvar)$, $f\equiv0$, and $h(\optvar)=\frac{1}{2}\|\ipop\optvar-\datavar\|^2$ for PGD and  $h_i(\optvar)=\frac{1}{2}\|\ipop_i\optvar-\datavar_i\|^2$ for SGD, SAGA and SVRG;
\item[Primal-dual methods] $g(\optvar)=\alpha_I \|\nabla \optvar\|_1+\iota_{[0,\infty)^d}(\optvar)$, $h\equiv0$, $\optop=\ipop$ with $f(\dualvar)=\frac{1}{2}\|\dualvar-\datavar\|^2$ for PDHG and $\optop_i=\ipop_i$ with $f_i(\dualvar)=\frac{1}{2}\|\dualvar-\datavar_i\|^2$ for SPDHG.
\end{description}
The proximal operator of $g$ is approximated by 100 FGP iterations \cite{beck2009fast}.
Stochastic methods use $\numh=60$ (for SGD, SAGA and SVRG), and $\numf=60$ (for SPDHG), subsets of the data. 
Operators $\ipop_i$ are defined analogously to the parallel-beam CT case in Section~\ref{subsec:parallel_CT}: starting with the angle $i/\numh$ (or $i/\numf$) we take every $\numh$-th (or $\numf$-th) angle from the equidistant partition of $[0,2\pi)$ into $360$ angles. 
As a result, the norms $\|\ipop_i\|$ are all nearly the same.
\replaced{SGD uses a constant step{\removed{-}}size $\tau=1/(2\numh\Lmax)$, and SAGA and SVRG use $\tau=1/(3\numh\Lmax)$.}{Here we compare to SGD with a constant step{\removed{-}}size $\tau=1/(2\numh\Lmax)$ as is commonly used in applications. SAGA and SVRG use a constant stepsize of $\tau=1/(3\numh\Lmax)$.}
SPDHG is implemented following \cite[Algorithm 2]{Ehrhardt2019pmb}, with $\sigma=\gamma\rho / K_{\text{max}}$, $\tau= 1/(\numh\gamma K_{\text{max}})$, with $\rho=0.99$, $\gamma=1$ and $K_{\text{max}} = \max_i \|\ipop_i\|$.

\begin{figure}[htbp]
    \centering
    \subfloat[High noise]{\includegraphics[width=0.33\textwidth]{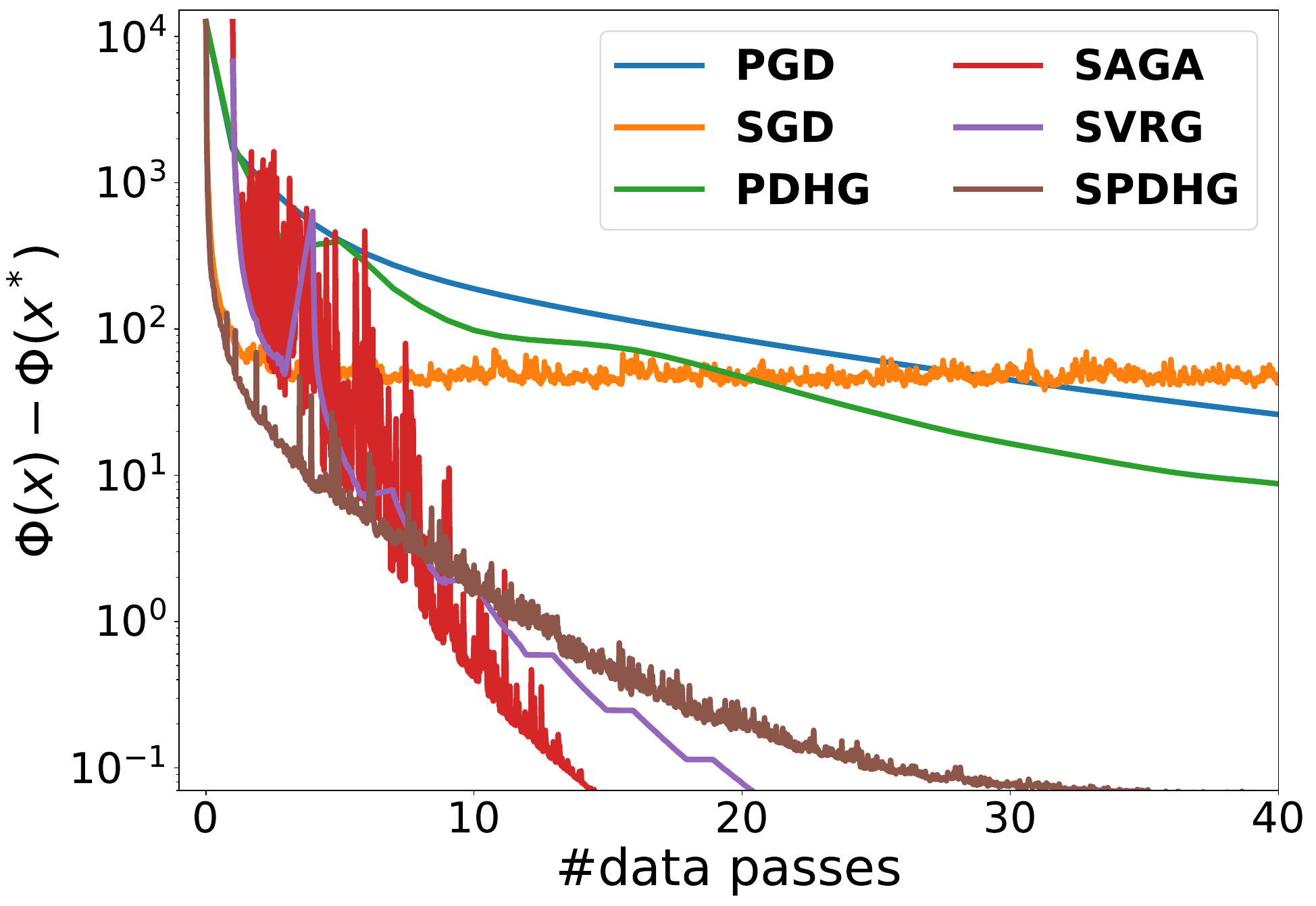}\label{fig:CW_figure1}}
    \hfill
    \subfloat[Medium noise]{\includegraphics[width=0.33\textwidth]{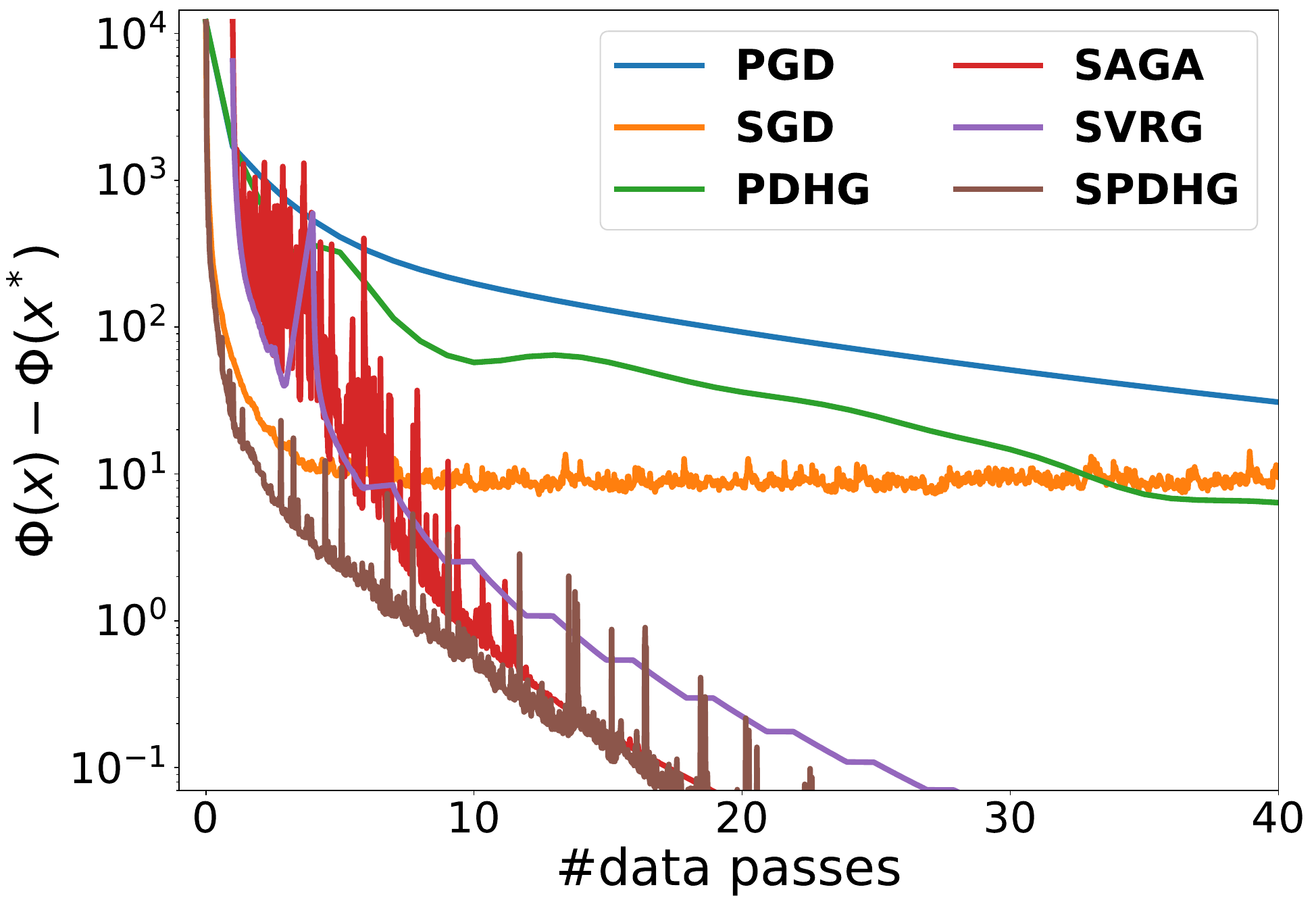}\label{fig:CW_figure2}}
    \hfill
    \subfloat[Low noise]{\includegraphics[width=0.33\textwidth]{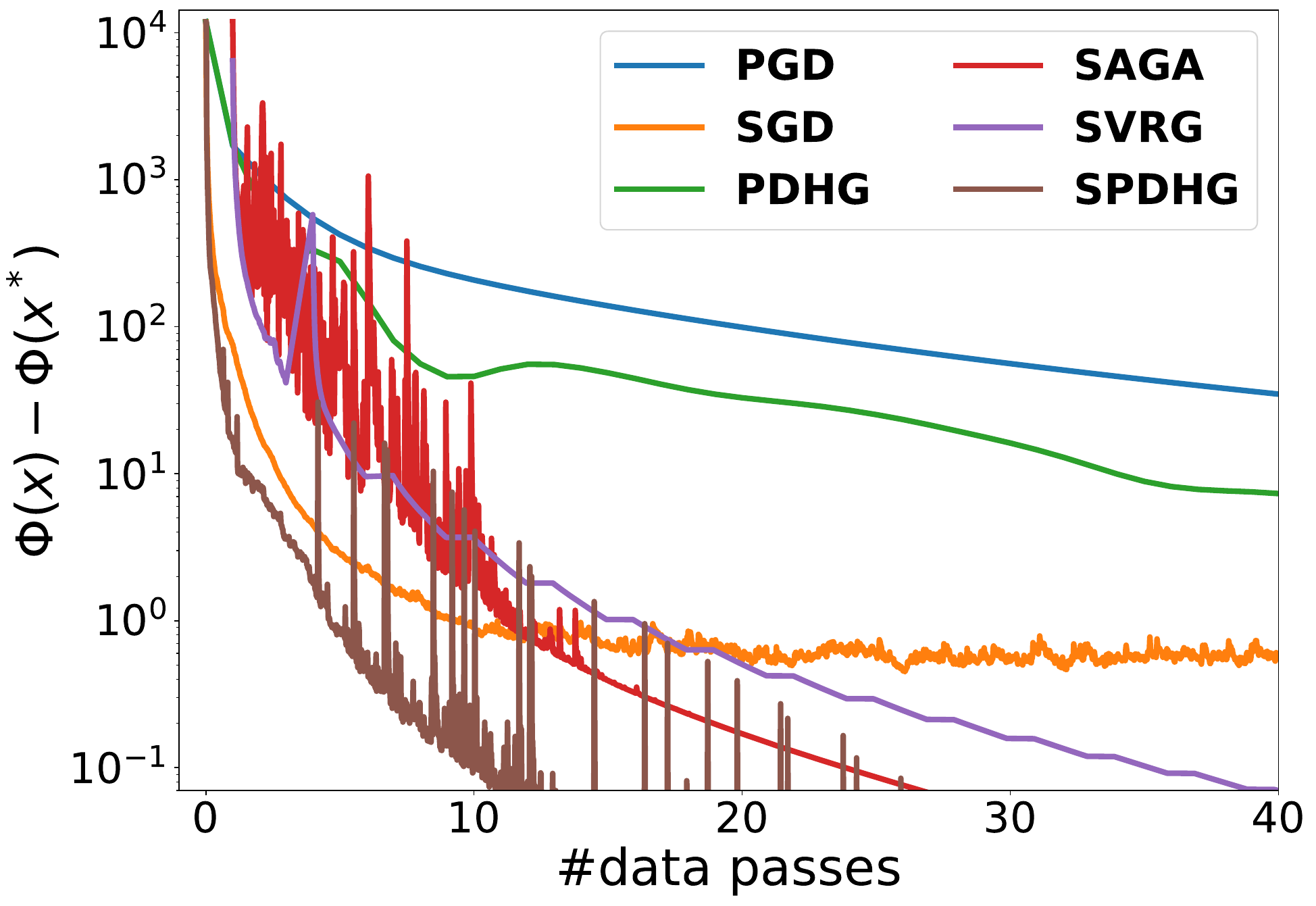}\label{fig:CW_figure3}
    }
    \caption{Comparison of reconstruction performance for a walnut from fan-beam CT measurements with three noise levels}
    \label{fig:CW_basic}
\end{figure}

The results in Figure \ref{fig:CW_basic} show that the noise level can have a noticeable effect on the convergence behaviour. 
First, for lower noise levels, SGD is competitive with PGD within the given optimisation window, whereas for higher levels it is less competitive. However, SAGA and SVRG are clearly ahead of PGD regardless of the noise level.
Second, SPDHG is the fastest method for lower noise, but its performance decays (relative to SVRG and SAGA) as the noise increases.
In spite of the superior convergence behaviour of SAGA, SVRG and SPDHG, running a large number of iterations might not be desirable due to computational or time constraints or even needed if the sought features of the image are recovered early.
It is therefore of interest to investigate which algorithms recover a good quality reconstruction within a limited iteration budget, which image details are missing, and whether an increased computational or memory budget is justified or whether SGD would suffice.
To investigate this, we look at reconstructions with SGD, SVRG, SAGA, and SPDHG within the first $10$ data passes for medium noise.

\begin{figure}
    \centering
    \subfloat[SGD]{\includegraphics[width=\textwidth]{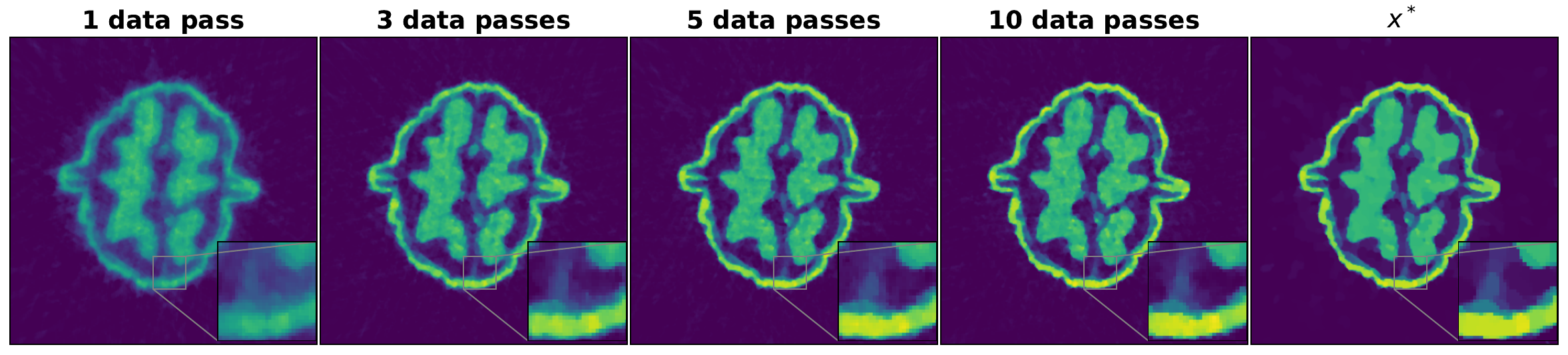}    \label{fig:walnut_sgd_iter_recons}}\newline
    \subfloat[SVRG]{\includegraphics[width=\textwidth]{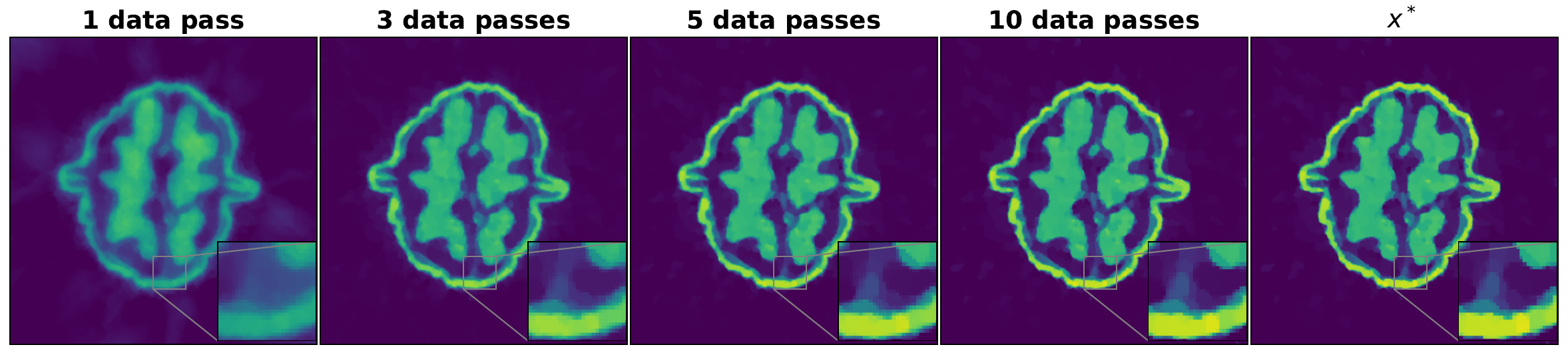}\label{fig:walnut_svrg_iter_recons}}\newline
    \subfloat[SAGA]{\includegraphics[width=\textwidth]{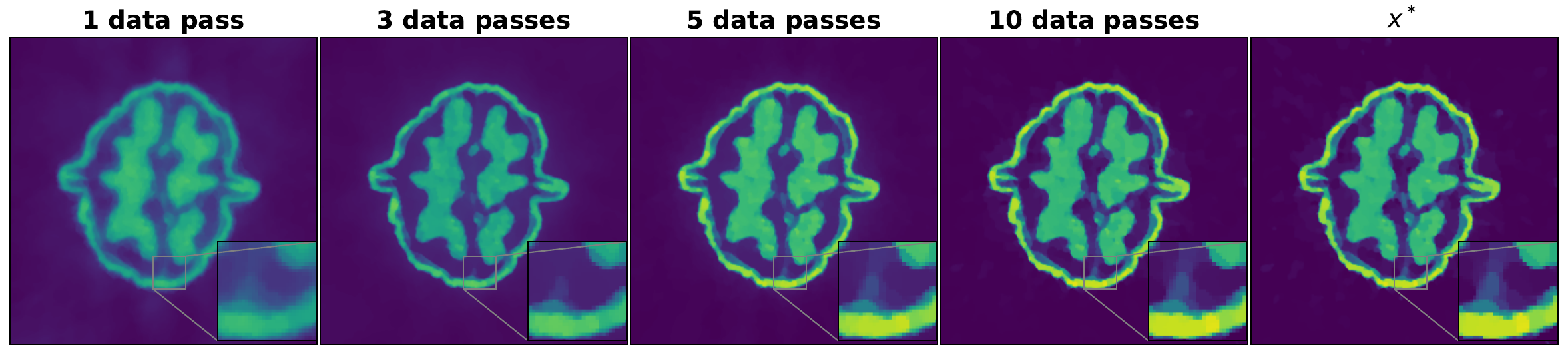}\label{fig:walnut_saga_iter_recons}} \newline
    \subfloat[SPDHG]{\includegraphics[width=\textwidth]{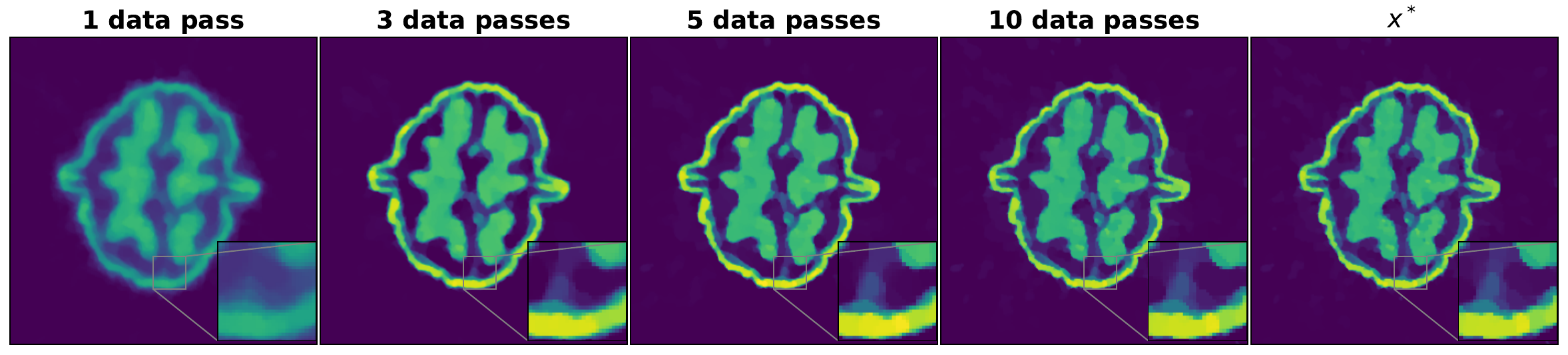}\label{fig:walnut_spdhg_iter_recons}}
    \caption{Reconstruction of the walnut from fan-beam CT with medium noise using SGD, SVRG, SAGA, and SPDHG. The colour scheme in all the images is such that the smallest and largest value corresponds to the maximum and minimum values in the reference solution.}\label{fig:walnut_iter_recons}
\end{figure}

The results in Figure \ref{fig:walnut_iter_recons} confirm that all the stochastic algorithms perform well and recover a visually satisfying solution quickly.
Moreover, variance reduced algorithms show a steady convergence behaviour and better quality of detail.
SGD performs well but fails to reconstruct the desired image structures promoted by the TV regulariser and retains more reconstruction noise, within this time frame. 
However, depending on the problem, only a handful of data passes might be needed to recover general features of the solution.
Thus, a sound strategy to improve the performance of an image reconstruction method would be to warm start the reconstruction with a few data passes with standard SGD, before employing a more advanced method.
This would also help reduce the initially high variance seen with SAGA, which updates a gradient stored in the table only once the corresponding batch (or subset) index is selected (which introduces an initial bias towards selected subset gradients) and may result in a reconstruction of unreliable quality within the first $10$ data passes.

\begin{figure}
    \centering
    \includegraphics[width=0.4\textwidth]{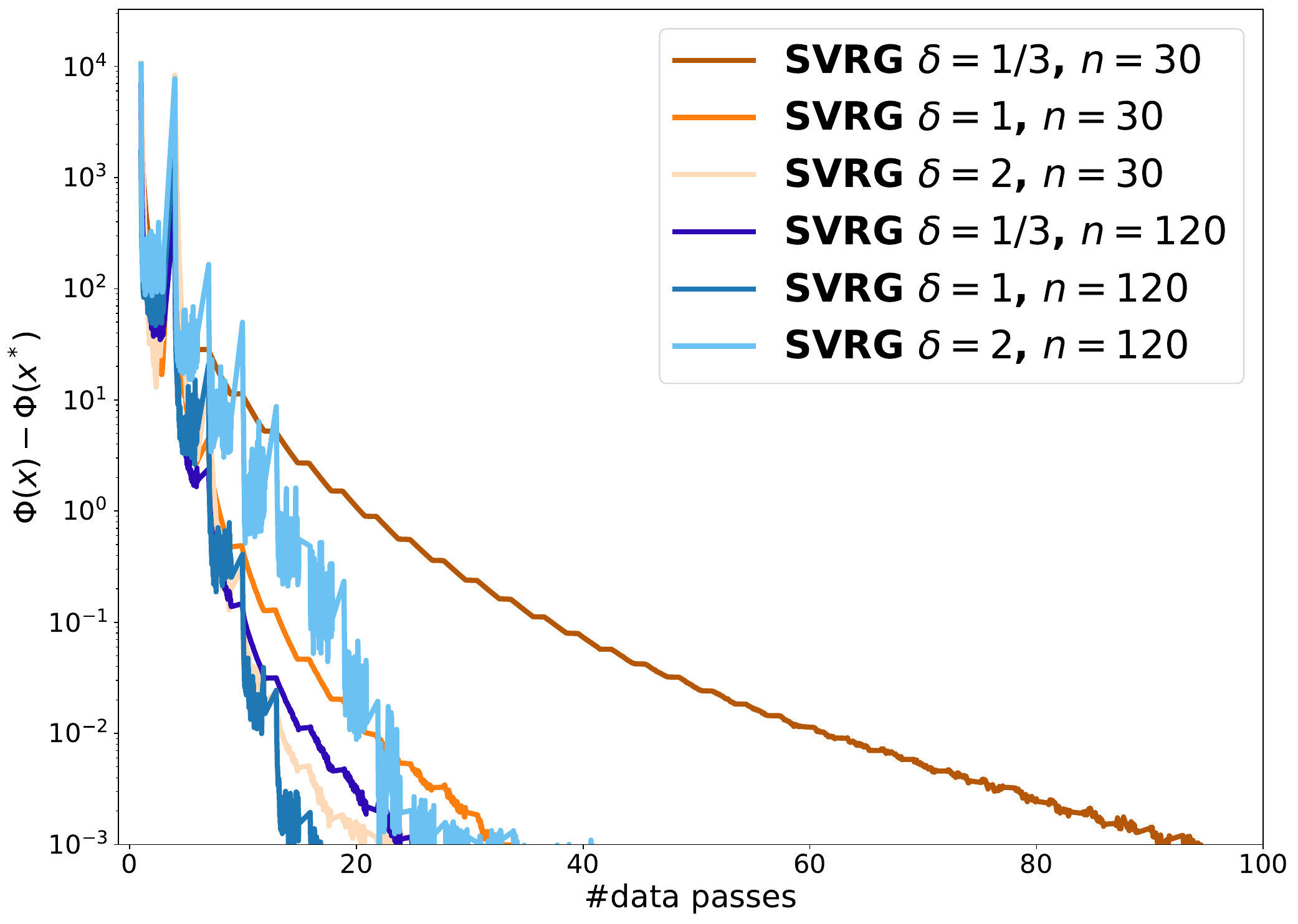}
    \caption{Robustness of SVRG with respect to mis-specification of the step{\removed{-}}size for the fan-beam reconstruction of the walnut.}
    \label{fig:svrg_stepsizes}
  \end{figure}

As argued in Section \ref{sec:challenges}, theoretical step{\removed{-}}size regimes are often pessimistic, in the sense that they are intended to give the largest value that ensures convergence towards a minimum, with a given method, for all functions and data in the given function class. 
In Figure \ref{fig:svrg_stepsizes} we take a closer look at this behaviour by looking at the convergence trajectory of SVRG for medium noise ($I_0=250$), with $\numh\in\{30,60, 120\}$ and constant step{\removed{-}}sizes $\tau = \delta /(\numh \Lmax)$ for $\delta \in \{1/3, 1, 2\}$. 
SVRG with step{\removed{-}}size $\delta=7/3$ no longer showed a convergent behaviour within the studied optimisation window, for all the investigated batch sizes.
The results show that SVRG is stable with respect to the step{\removed{-}}size choice.
The variance progressively increases as the number of subsets $\numh$ increases: when $\numh=120$, SVRG performs best for $\delta = 1$ and the worst for $\delta=2$, whereas for $\numh=30$ it performs the best for $\delta = 2$ and the worst for $\delta = 1/3$. 
This behaviour is (roughly) consistent across the studied noise levels. 
The same experiment was carried out for SAGA, which showed divergent behaviour already for $\delta = 1$.

\subsection{Image Reconstruction of Simulated PET Data}

The sensitivity of the convergence behaviour of an iterative method on the Lipschitz constant can have severe consequence in reconstruction tasks where the Lipschitz constant is not available, hard to estimate, or is highly varying. 
A prototypical case is PET, where the data fidelity term is the Poisson log-likelihood, or equivalently, the Kullback--Leibler (KL) divergence.
This data fidelity term is in general not globally Lipschitz smooth. However, in case of PET it is Lipschitz smooth for any realistic PET scan which has \replaced{non-zero}{nonzero} scatter and background events. 
However, the corresponding Lipschitz constant is pessimistic and usually not used in practice.
Stochastic (and \replaced{non-stochastic}{nonstochastic}) gradient methods are nevertheless used for PET, and despite the fact that their application is in principle more heuristic, stochastic methods can show fast and reliable performance \cite{Twyman2023rdp,kereta2021svrem,Chambolle2018spdhg}. 
To provide further insight into this behaviour, we report the numerical results from \cite{papoutsellis2024psmr} that use the stochastic framework from  the CIL library. 
The results concern the reconstruction of a simulated PET thorax dataset with simulated Poisson noise.

In terms of the optimisation template \eqref{eq:opttemplate} we write the objective as follows, where $r$ denotes the background and scatter events.
\begin{description}
\item[Gradient based methods] $g(\optvar)=0.1 \|\nabla \optvar\|_1+\iota_{[0,\infty)^d}(\optvar)$, $f\equiv0$, and $h(\optvar)=\text{KL}(\datavar|\ipop\optvar+r)$ for FISTA and  $h_i(\optvar)=\text{KL}(\datavar_i|\ipop_i\optvar+r_i)$ for SGD, SAGA and SVRG;
\item[SPDHG] $g(\optvar)=0.1 \|\nabla \optvar\|_1+\iota_{[0,\infty)^d}(\optvar)$, $\optop_i=\ipop_i$ with $f_i(\dualvar)=\text{KL}(\datavar_i|\dualvar+r_i)$, $h\equiv0$.
\end{description}
The numerical study compares SGD, SAGA, SVRG, Loopless-SVRG (L-SVRG) and SPDHG, therein termed Prox-SGD, Prox-SAGA, Prox-SVRG and Prox-LSVRG, respectively, against FISTA, in terms of the relative objective error and the distance to the reference solution.
Stochastic methods use $32$ equidistant subsets of data. Stepsize regimes were chosen via a rough grid search, aiming to optimise the performance in terms of the distance to the solution. FISTA uses stepsize $\tau=0.0455$. SGD uses $\tau=0.1$, SAGA uses $\tau=0.025$, SVRG and L-SVRG use $\tau=0.05$, and L-SVRG uses the probability parameter $p=0.0078$.  
SPDHG uses 
$\sigma=\gamma\rho / K_{\text{max}}$, $\tau= 1/(\numh\gamma^2 K_{\text{max}})$, with $\rho=0.99$, $\gamma=0.25$ and $K_{\text{max}} = \max_i \|\ipop_i\|$.
 
\begin{figure}
    \centering
    \subfloat[Convergence of iterates]{
    \includegraphics[width=0.49\textwidth]{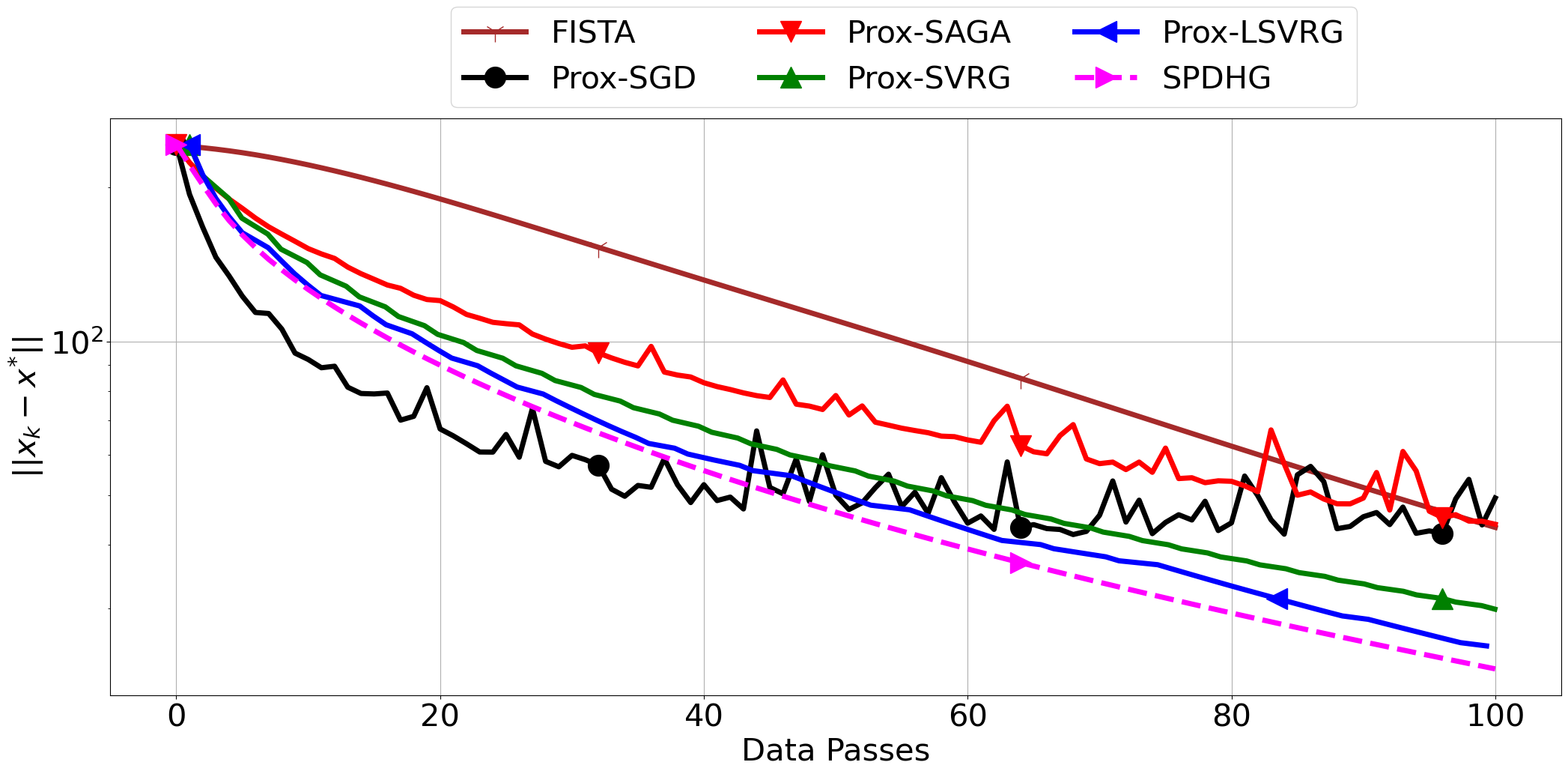} \label{fig:pet_convergence}}
    \subfloat[Convergence of objective value. Function $F$ is equal to the function $\Psi$ from template \eqref{eq:opttemplate_finite_sum}]{
    \includegraphics[width=0.49\textwidth]{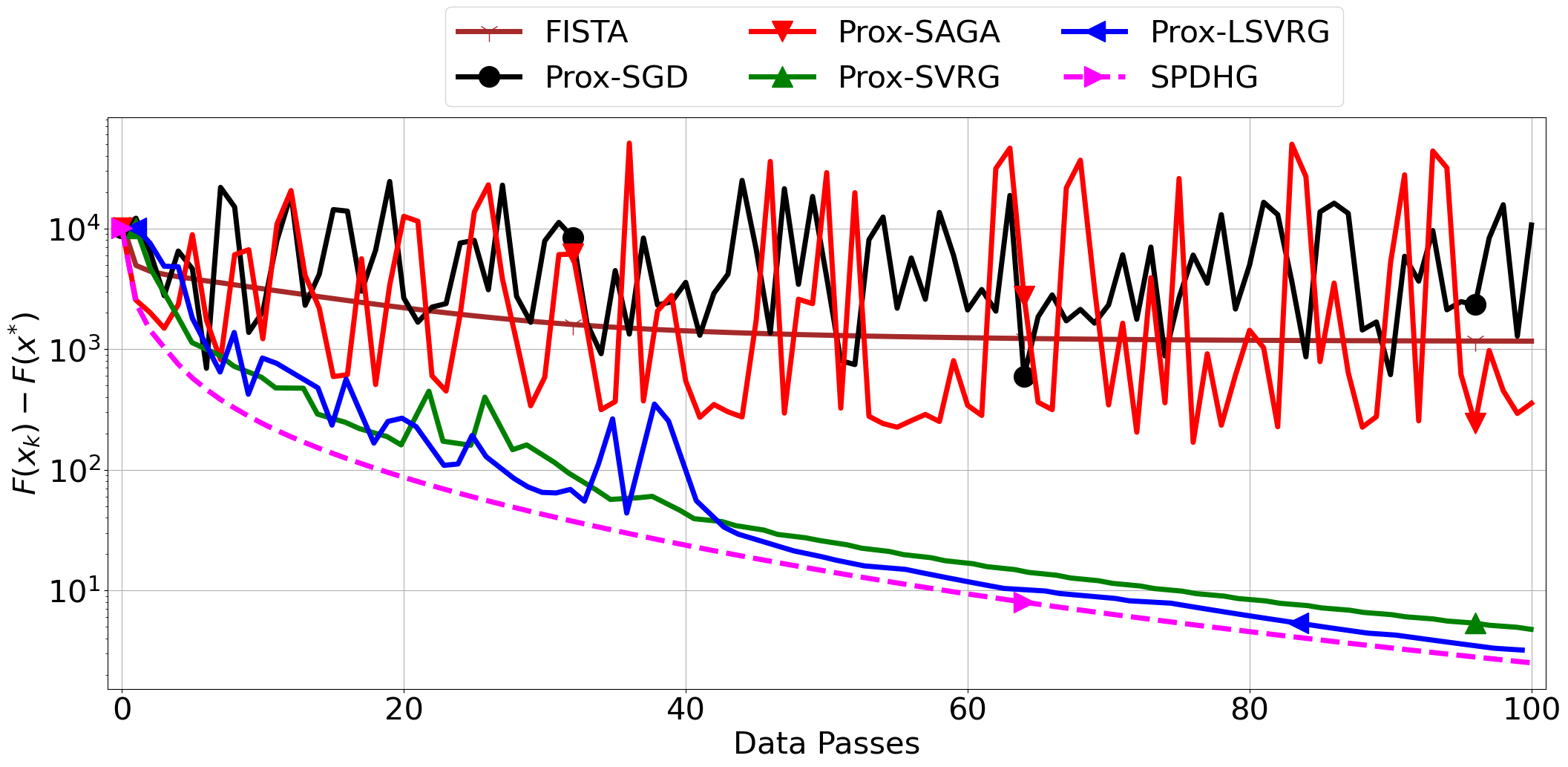} \label{fig:pet_obj}}
    \caption{PET reconstruction of a 2D thorax, from \cite{papoutsellis2024psmr}. Here Prox-SGD, Prox-SAGA, Prox-SVRG and Prox-LSVRG refer to SGD, SAGA, SVRG, Loopless-SVRG.}
    \label{fig:stepsizes}
\end{figure}

The results, reported in Figure \ref{fig:pet_convergence}, show that all the methods show good behaviour in terms of the distance to the optimal solution, with SAGA being different compared to its performance in CT. However, in terms of the objective function error we see a somewhat different picture. 
Namely, SVRG, L-SVRG and SPDHG perform well, but SAGA is (qualitatively) indistinguishable from SGD and fails to show variance-reducing behaviour (which can be attributed to the fact that the variance is only expected to be reduced as we get closer to the solution).
This diverging optimisation behaviour is particularly interesting due to the fact that SVRG and SAGA have the same step{\removed{-}}size conditions in many problems, but in practice often show a different behaviour.

\section{Outlook}\label{sec:outlook}

Optimisation algorithms are a backbone for solving inverse problems and stochastic approaches becoming increasingly relevant for large-scale applications. We have discussed state-of-the-art algorithms for inverse problems with a particular focus on the toolbox nature of the topic. Virtually all of the relevant algorithms are assembled from a fairly small range of tools from both the deterministic and the stochastic optimisation toolkit. Particularly, the stochastic toolkit has been developed with machine learning in mind but it turns out to be of excellent use for inverse problems, too.

Over the next decade, we expect this trend to continue, with more and more applications realising the potential that stochastic optimisation offers. Furthermore, we expect stochastic optimisation tools to be developed specifically for inverse problems which may not be relevant to machine learning. Potential examples are sketching in the image domain and nonconvexity from nonlinear forward operators. We also expect that theory could be developed specifically for inverse problems. For instance, it has been shown that stochastic gradient descent can be much faster in the so called interpolation regime \cite{Vaswani2019}. In the context of inverse problems this regime implies that there is no noise and no regularisation which of course is unrealistic. Instead, there could be theory tailored to inverse problems explaining which algorithms are fastest in a specific noise regime, thereby potentially explaining the observations we made in Section~\ref{sec:numerics}.

Finally, as we indicated throughout the manuscript in various further reading sections, an exciting path forward is the interplay between inverse problems, optimisation and machine learning. This could be in the context of modelling, such as learning regularisation models or learning forward operators or their sampling, but also in the context of learning optimal algorithms to solve inverse problems. We believe that the current optimisation-based paradigm with learned explicit regularisation will become a dominating computational paradigm in inverse problems, while developing efficient optimisation algorithms tailored for these applications will continue to be important to practitioners. To the best of our knowledge this has not been attempted in the context of stochastic algorithms. In our view, this development is far from over and more research is needed to fully understand and exploit the connections in the context of inverse problems.

\ack

We acknowledge support from the EPSRC: MJE (EP/Y037286/1, EP/S026045/1, EP/T026693/1, EP/V026259/1) and ZK (EP/X010740/1).
JL is supported by the ``Fundamental Research Funds for the Central Universities’’, the National Science Foundation of China (BC4190065) and the Shanghai Municipal Science and Technology Major Project (2021SHZDZX0102).

\appendix

\addtocontents{toc}{\protect\setcounter{tocdepth}{-1}}

\section{Mathematical Preliminaries and Notation}\label{app:prelims}

\addtocontents{toc}{\protect\setcounter{tocdepth}{+1}}

\addcontentsline{toc}{section}{{Appendix A}\hspace{1em}Mathematical Preliminaries and Notation}

Throughout this paper, $\optspace$ and $\dualspace$ denote generic real Hilbert spaces with inner products $\langle \cdot,\cdot\rangle$ and associated norms $\|x\| = \sqrt{\langle x, x\rangle}$. It will be clear from the context which inner product is being used. The (Hermitian) adjoint of a linear and bounded operator $A: \optspace \to \dualspace$ is denoted as $A^\ast$ and $I$ is the identity operator.
We use $\|A\|:=\sup_{\|x\|=1} \|Ax\|$ to denote the induced operator norm. 

By $(a, b]$ we denote the interval $a < x \leq b$. If $\optspace = \R^d$, we write $x_i$ for the $i$-th entry of a vector $x \in \R^d$. For $1\leq p<\infty$ we denote the $p$-norm of a vector $x\in\R^d$ by $\|x\|_p=(\sum_{i=1}^d |x_i|^p)^{1/p}$. 
\new{$I$ denotes the identity operator.}

\paragraph{Convexity} 
A function $f\new{:\optspace \to \RI}$ is called \emph{convex} if  
\begin{equation*}
f(\alpha x_1+(1-\alpha)x_2)\leq \alpha f(x_1)+(1-\alpha) f(x_2),~ \forall \alpha \in [0, 1] \new{,~ \forall x_1,x_2\in\optspace,}
\end{equation*}
and it is \emph{strictly convex} if the above equality is strict. 
It is \emph{$\mu$-strongly convex} if 
\begin{align*}
f(x) - \frac{\mu}{2} \|x\|^2\new{,~\forall x\in\optspace}
\end{align*}
is convex. 
It is \emph{$L$-Lipschitz continuous} if for some $L \geq 0$ it holds that
\begin{equation}\label{eqn:lipschitz_cont}
\|f(x_1) - f(x_2)\| \leq L \|x_1-x_2\|\new{,~ \forall x_1,x_2\in\optspace} . 
\end{equation} 
A function $f$ is called \emph{$L$-smooth} if it is differentiable and has $L$-Lipschitz continuous gradient, that is
\begin{equation}\label{eqn:L_smooth}
\|\nabla f(x_1) - \nabla  f(x_2)\| \leq L \|x_1-x_2\|,~ \forall x_1,x_2 \in \optspace . 
\end{equation} 
If $f$ is twice differentiable and $\mu$-strongly convex, then all the eigenvalues of its Hessian are all larger than $\mu$ (i.e.\ $\nabla^2 f(x)\succeq
 \mu I) $, and it is $L$-smooth if they are all bounded by $L$ (i.e.\ $\nabla^2 f(x)\preceq L I$).

A function $f$ is called \emph{proper} if its effective domain $\{x\in \optspace:f(x) < \infty\}$ is nonempty. It is also called closed, if its sublevel set $\{x\in \optspace: f(x)\leq \gamma \}$ is closed for each $\gamma\in\R$.

\paragraph{Nonsmooth analysis}
For a proper closed and convex function $f$,  its \emph{subdifferential} is defined by
\begin{equation}\label{eqn:subdifferential}
    \partial f(x) =\{p\in \optspace: f(z)\ge f(x)+\langle p, z-x\rangle, \text{ for all } z \in \optspace\}.
\end{equation}
\new{An element from the subdifferential is called a subgradient.}
The \emph{Fenchel (or convex) conjugate} of a function $f$ is defined as
\begin{equation}\label{eqn:convex_conjugate}
f^\ast(y) = \sup_{x\in \optspace} \left\{\langle x, y\rangle - f(x)\right\}.
\end{equation}
Note that $f^\ast$ is closed convex regardless the convexity of $f$. 
The bi-conjugate of $f$ is defined by
\begin{equation}\label{eqn:bi_conjugate}
f^{\ast\ast}(x) = \sup_{y\in \optspace} \left\{\langle x, y\rangle - f^{\ast}(y)\right\}.
\end{equation}
The \emph{indicator function} $\iota_{\mathcal C}$ of a set $\mathcal C\subset \optspace$ is defined as 
\begin{equation}\label{eqn:indicator_function} 
\iota_{\mathcal C}(x)=\begin{cases}0,&\text{ if } x\in \mathcal C , \\ \infty,&\text{otherwise} . \end{cases}
\end{equation}
The \emph{proximal operator} with parameter $\tau > 0$ of a proper closed convex function $f:\optspace\rightarrow \R_\infty$ is the mapping $\prox_{\tau f}: \optspace \rightarrow \optspace$ defined by
\begin{equation}\label{eqn:prox_operator} 
\prox_{\tau f} (z)=\argmin_{x \in \optspace} \left\{ f(x)+\frac{1}{2\tau}\|x-z\|^2\right\}.
\end{equation}
Due to the convexity of $f$ (and thus the strong convexity of the objective), the minimiser always exists and is unique. 
If $\optspace = \R^\recondim$ with the Euclidean norm, then the proximal operator of many convex and \replaced{non-convex}{nonconvex} functions, e.g.\ norms and indicator functions, has a closed-form expression \cite{combettes2011proximal,beck2017first}, and is known for many other functions \cite{proximityoperatorProxRepository}. \new{Such functions are called prox-friendly.} For instance, the proximal operator of the indicator function of the \replaced{non-negative}{nonnegative} orthant $[0,\infty)^d$ has a particularly simple expression
$\left[\prox_{\tau \iota_{[0, \infty)^d}}(x)\right]_i = \max(x_i, 0)$. 

\new{An important class of cases when the proximal operator of the composition term $f\circ \optop$ can be explicitly computed is when $\optop\optop^\ast=\alpha I,~\alpha>0$ and the proximal operator of $f$ can be computed.} In this case the proximal operator of $f \circ \optop$ can be computed as $\prox_{f\circ \optop}(x) = x + \tfrac{1}{\alpha} \optop^\ast \big( \prox_{\alpha f}(\optop x) - \optop x \big)$ \cite{beck2017first}. 

\section*{References}

\bibliographystyle{alpha}
\bibliography{references}

\end{document}